\documentclass[11pt,a4paper]{amsart}
\usepackage{iftex}
\ifPDFTeX
\usepackage[utf8]{inputenc}
\usepackage[T1]{fontenc}
\usepackage{lmodern}
\else
\usepackage{fontspec}
\fi
\usepackage[margin=1in]{geometry}
\usepackage{mathtools}
\usepackage{amssymb}
\usepackage{amsthm}
\usepackage{mathrsfs}
\usepackage{stmaryrd}
\usepackage{enumitem}
\usepackage[numbers,sort]{natbib}
\usepackage[colorlinks=true,linkcolor=blue,citecolor=blue,urlcolor=black]{hyperref}

\numberwithin{equation}{section}
\title{On the regularity of branching loci for $\mathbb Z_2$-harmonic spinors}
\author[Yang Li]{Yang Li}
\address{Department of Pure Mathematics and Mathematical Statistics, Cambridge University, Wilberforce Road, Cambridge CB3 0WB, UK}
\email{yl454@cam.ac.uk}

\author[Davide Parise]{Davide Parise}
\address{Department of Mathematics, Zeeman Building, University of Warwick, Gibbet Hill Road, Coventry CV4 7AL, UK}
\email{Davide.Parise@warwick.ac.uk}

\author[Yichen Zhang]{Yichen Zhang}
\address{Department of Mathematics, University College London, Gower Street, London WC1E 6BT, UK}
\email{yichen.zhang.24@ucl.ac.uk}

\date{\today}

\newtheorem{maintheorem}{Theorem}
\newtheorem{theorem}{Theorem}[section]
\newtheorem{lemma}[theorem]{Lemma}
\newtheorem{proposition}[theorem]{Proposition}
\newtheorem{corollary}[theorem]{Corollary}
\newtheorem{example}[theorem]{Example}
\theoremstyle{definition}
\newtheorem{definition}[theorem]{Definition}
\newtheorem{remark}[theorem]{Remark}

\newcommand{\ie}{\emph{i.e.} }
\newcommand{\cf}{\emph{cf.} }
\newcommand{\R}{\mathbb{R}}
\newcommand{\C}{\mathbb{C}}
\newcommand{\gradG}{\nabla^{G}}

\makeatletter
\def\@part[#1]#2{%
  \refstepcounter{part}%
  \addcontentsline{toc}{part}{\thepart\hspace{1em}#1}%
  {\parindent\z@\raggedright\normalfont\Large\bfseries\boldmath
   \partname~\thepart.\quad#2\par}%
  \nobreak\vskip 3ex\@afterheading}
\makeatother

\begin{document}
\maketitle
\begin{abstract}
We show that the branching locus of a $\mathbb Z_2$-harmonic spinor is
smooth (resp. real-analytic) near every point of frequency $\frac{1}{2}$ for smooth (resp. real-analytic) ambient defining data. This partially strengthens the previous regularity results obtained by Clifford Taubes and Boyu Zhang.
\end{abstract}

\section{Introduction}
	\label{sec:front-introduction}
	Multivalued harmonic functions, $1$-forms, and spinors arise naturally in geometric measure theory \cite{DLS11}, gauge theory \cite{Tau13}, calibrated geometry \cite{He23}, and special holonomy \cite{Don16}. In this paper we study the regularity of the branching locus of a $\mathbb Z_2$-harmonic spinor.
	
	Let $M$ be a Riemannian manifold of dimension $n\geqslant2$, equipped with a real or complex spinor bundle $E$ carrying a fibre metric $h$, Clifford multiplication $\mathfrak c$, and a compatible connection $\nabla^E$, and write $\mathcal D=\sum_i\mathfrak c(e_i)\nabla^E_{e_i}$ for the associated Dirac operator in some local orthonormal frame $\{e_i\}$. We call a symmetric two-valued section $u$ a $\mathbb Z_2$-harmonic spinor if $u\in W^{1,2}_{\operatorname{loc}}(M,\mathcal{A}_2^{\operatorname{sym}}(E))$ and
	\[
	\mathcal{D}u=0\quad\text{almost everywhere}.
	\]
	Here $\mathcal A^{\mathrm{sym}}_2(E_x):=\{\llbracket v\rrbracket+\llbracket-v\rrbracket:v\in E_x\}$ denotes the space of symmetric unordered pairs in the fibre $E_x$ for $x\in M$, and $\mathcal Du$ is well-defined almost everywhere for $W^{1,2}$ two-valued sections \cite{DLS11}.
	
	It turns out that $u$ admits a H\"older continuous representative (see Section~\ref{sec:holder-regularity}). The branching set $\mathcal{B}_u$ and the zero locus $Z_u$ are defined by
	\[
	\mathcal{B}_u:=\{x:u\text{ has no \textit{continuous} local lift near }x\}
	\subset Z_u:=\{x:u(x)=2\llbracket0\rrbracket\}.
	\]
	The branching set $\mathcal B_u$ is closed.
	
	Our main result concerns branch points of frequency $\frac{1}{2}$.
	\begin{maintheorem}
		\label{thm:main-regularity}
		Let $u$ be a nonzero $\mathbb Z_2$-harmonic spinor and let $p_0\in\mathcal B_u$ have frequency $N_u(p_0,0)=1/2$. Then $\mathcal B_u$ is a smooth (resp. real-analytic) embedded submanifold of codimension two near $p_0$ for smooth (resp. real-analytic) ambient defining data $(M,E,h,\mathfrak{c},\nabla^E)$.
	\end{maintheorem}
	
	\begin{remark}
		We place this result in context.
		\begin{enumerate}
			\item (\textit{Earlier regularity results})
			Clifford Taubes \cite{Tau14} developed the foundational theory of $\mathbb Z_2$-harmonic spinors and their zero loci, assuming continuity of $|u|$ and the local growth bound
			$\int_{B_r(p)}|u|^2\,dV_g\lesssim_p r^{n+\varepsilon}$
			for $p\in Z_u$, $r\ll1$, and some $\varepsilon>0$.
			Under the same hypotheses, Boyu Zhang \cite{Zha17} proved rectifiability and locally finite Minkowski content of $Z_u$ in dimension four.
			In our setting, continuity and the growth bound are proved \textit{a priori} in Part I.
			
			\item (\textit{Sharpness for smooth ambient data})
			Jiahuang Chen and Siqi He \cite{CH26} construct $\mathbb Z_2$-harmonic $1$-forms on $\mathbb R^3$ with a suitable smooth metric for which
			$Z_u=\mathcal B_u=\mathcal B_u^{sm}\cup
			(\overline{\mathcal B_u^{sm}}\setminus\mathcal B_u^{sm})$.
			The points of the smooth part $\mathcal B_u^{sm}$ have frequency $\frac{1}{2}$, while its complement in $Z_u$ is the standard Cantor set;
			the authors explained to us in private communication that the construction also applies to fat Cantor sets of Hausdorff dimension $1$ and positive $\mathcal H^1$-measure.
			These examples illustrate the sharpness of Theorem~\ref{thm:main-regularity} for smooth ambient data. The real-analytic case remains under investigation.
			
			\item (\textit{Frequency-$\frac{1}{2}$ and geometric applications})
			Certain gluing constructions in gauge theory and calibrated geometry \cite{HPW26,Par24} involve \textit{non-degenerate} $\mathbb Z_2$-harmonic spinors with a smooth branching locus and
			$|u(x)|\sim\operatorname{dist}(x,\mathcal B_u)^{1/2}$
			near the branching locus.
			To address the loss of derivatives in the deformation theory, one uses Nash--Moser methods, which require sufficiently high regularity of $\mathcal B_u$ \cite{Don21,Par23,HPW26}.
			Frequency measures the local vanishing order, so the frequency-$\frac{1}{2}$ assumption gives a weak formulation of non-degeneracy. Theorem~\ref{thm:main-regularity} then supplies the regularity of the branching locus needed for this deformation theory.
			
			\item (\textit{Regularity of singular sets})
			The branching set is itself part of the unknown, and no equation is prescribed for it \textit{a priori}. Its regularity is therefore a codimension-two free boundary problem.
			Related regularity questions arise in the classical obstacle problem \cite{Caf77,FS19,FZ25} and the Signorini (thin obstacle) problem \cite{KPS15}. They also arise in the study of cylindrical tangent cones in minimal submanifolds \cite{Sim93}, multivalued harmonic functions \cite{KW13,KW17,Kru14}, singular sets of area-minimizing integral currents \cite{DMS23,KW26}, and cylindrical singular sets in mean curvature flow \cite{SWX25}.
		\end{enumerate}
	\end{remark}
	
	\subsection{Proof strategy}
	\label{sec:front-proof-strategy}
	The proof has two parts. We first show that the branching locus is a $C^{1,\alpha}$ graph and obtain quantitative estimates for approximation by its tangent models. We then use certain components of the spinor as normal coordinates and prove higher regularity in the tangential directions.
	
	\textbf{Part I: $C^{1,\alpha}$-regularity of the branching locus and quantitative estimates.}
	We first show that a $\mathbb Z_2$-harmonic spinor locally minimises a modified Dirichlet energy $\mathscr J$. The key input is an integrated Weitzenb\"ock formula, obtained from a transgression form on the total space $E$. This observation was motivated by work on branched deformations of special Lagrangians \cite{He23} and by the compactness problem for Fueter sections with values in moduli spaces of monopoles of charge $\geqslant2$ \cite{EL24, LO}.
	
	Local minimality yields stationarity and the variational estimates needed to establish almost monotonicity of frequency, H\"older continuity, and strong $W^{1,2}$-compactness for blow-ups \cite{DLS11,Stu19}. With these inputs, we adapt the excess-decay method of Leon Simon \cite{Sim93}, as developed for two-valued harmonic functions by Krummel--Wickramasekera \cite{KW13}.
	
	The main output of Part I is the following.
	\begin{maintheorem}[$C^{1,\alpha}$-regularity with quantitative estimates]
		\label{thm:main-tangent-estimates}
		For smooth ambient data, there is a ball $B(p_0,r_0)$ such that the set $\mathcal B_u\cap B(p_0,r_0)$ is a single $C^{1,\alpha}$-graph of codimension two over the spine at $p_0$, for some $\alpha\in(0,1)$. For each $p\in\mathcal{B}_u\cap B(p_0,r_0)$, there is a unique non-degenerate flat Dirac tangent model $\Phi_p$ with spine $T_p\mathcal{B}_u$. Moreover, the following quantitative estimates hold: 
		\begin{equation}
			\sup_{B_{r/2}(p)}\mathcal G(u,\Phi_p)\leqslant Cr^{1/2+\alpha},\qquad 
            r^{1-n}\int_{B_{r/2}(p)}\mathcal G(\nabla u,\nabla\Phi_p)^2\leqslant Cr^{2\alpha},\qquad 0<r\leqslant r_0,
			\label{eqn:main-tangent-height}
		\end{equation}
		where $\mathcal{G}$ denotes the distance between unordered two-valued sections.
	\end{maintheorem}
	
	\textbf{Part II: The partial Legendre transform and Campanato iterations.}
	Choose local coordinates $(z,y)$ in which the spine at $p_0$ is $\{0\}\times\mathbb R^{n-2}$. By Part I, the branching locus is the graph $\{(g(y),y)\}$ of a $C^{1,\alpha}$ function, and the tangent model at $p_0$ has the form $\Phi_{p_0}=a_0z^{1/2}+b_0\bar z^{1/2}$, with $|a_0|^2+|b_0|^2>0$. The quantitative approximation by these models allows us to use two real components of the spinor as new normal coordinates $\xi$. On the branched double cover, the resulting partial Legendre transform is described by
	\[
	X(\xi,y)=(\varphi(\xi,y),y),\qquad\widehat u(\xi,y)=(\xi,V(\xi,y)),\qquad U=(\varphi,V),
	\]
	where $\varphi$ is the inverse normal coordinate map and $V$ contains the remaining spinor components. The map $\varphi$ is even in $\xi$, while $V$ is odd. At $p_0=p(y_0)=(g(y_0),y_0)$, the transformed tangent model has the form
	\begin{equation}
		P=(a+B(y-y_0)+q(\xi),M\xi),
		\label{eqn:intro-tangent-model}
	\end{equation}
	where $a\in\mathbb C$, $B$ and $M$ are linear maps acting on $y-y_0$ and $\xi$, respectively, and $q$ is a homogeneous quadratic polynomial. The simplest example is the square model $U_0=(\xi^2,0)$.
	
	In these coordinates, the branching set is the fixed plane $\{\xi=0\}$, and its geometry is encoded in the equations for $U$. These form a quasilinear degenerate elliptic system with the inherited variational structure. The identity $g(y)=\varphi(0,y)$ reduces regularity of the original branching locus to tangential regularity of $\varphi$. We write $D_y^\beta=D_{y_1}^{\beta_1}\cdots D_{y_{n-2}}^{\beta_{n-2}}$ for a tangential derivative of order $|\beta|=\sum_{k=1}^{n-2}\beta_k$.
	
	The main result of Part II is the following.
	\begin{maintheorem}
		\label{thm:main-higher-regularity}
		For smooth (resp. real-analytic) ambient data, on every sufficiently small cylinder $Q_r(y_0)=B_r^2(0)\times B_{r^2}^{n-2}(y_0)$ centred at a frequency-$\frac{1}{2}$ branch point $(0,y_0)$ in the transformed coordinates, the fields $D_y^\beta U$ have a unique tangent polynomial of the form \eqref{eqn:intro-tangent-model} and satisfy
		\begin{equation}
			\|D_y^\beta U\|_{L^\infty(Q_r(y_0))}\leqslant C_\beta,
			\qquad(\textit{resp.}\ \|D_y^\beta U\|_{L^\infty(Q_r(y_0))} \leqslant C^{|\beta|+1}|\beta|!),
			\label{eqn:front-tangential-factorial-bound}
		\end{equation}
		where $C$ is independent of the multi-index $\beta$ in the real-analytic case.
		
		In particular, the branching set near $p_0=p(y_0)$ is smooth (resp. real-analytic). In Legendre coordinates, the $\mathbb Z_2$-harmonic spinor and the model coefficients $a=g(y)$, $B=D_yg$, $q$, and $M$ depend smoothly (resp. real-analytically) on $y$.
	\end{maintheorem}
	
	The proof proceeds in three steps.
	\begin{enumerate}
		\item We first study the square model $U_0=(\xi^2,0)$ in flat ambient space. Its linearised Euler--Lagrange equation $\mathscr B_0(Z,Y)=0$ is the weak formulation of the system
		\begin{equation}
			\left\{\begin{aligned}
				&\Delta_\xi Z_\varphi-\frac{4}{\xi}\partial_{\bar\xi}Z_\varphi+4|\xi|^2\Delta_yZ_\varphi=0,
                \\&\Delta_\xi Z_V+4|\xi|^2\Delta_yZ_V=0,
			\end{aligned}\right.
			\label{eqn:front-square-model-linearisation}
		\end{equation}
		where $Z_\varphi$ is even and $Z_V$ is odd in $\xi$. In weighted Sobolev and H\"older spaces adapted to Grushin Laplacian $\Delta_\xi+4|\xi|^2\Delta_y$, we prove a weighted Caccioppoli inequality and excess decay for weak solutions of \eqref{eqn:front-square-model-linearisation}.
		
		\item The Grushin Laplacian commutes with tangential derivatives. Differentiating the nonlinear system in the $y$-directions gives an equation of the schematic form
		\[
		\mathscr B_U(D_y^\beta U,Y)=\mathscr F_\beta(Y),
		\qquad \textit{for every admissible test field }Y,
		\]
		where $\mathscr B_U$ is the second variation of the transformed $\mathscr J$-energy and $\mathscr F_\beta$ contains only field derivatives with lower tangential order. The estimates from Part I control the difference between $\mathscr B_U$ and the form obtained by freezing at a transformed tangent model. A Campanato iteration then transfers the frozen excess decay to $D_y^\beta U$.
		
		\item To continue the induction, we need derivative estimates for the fields appearing in the source $\mathscr F_\beta$. We therefore upgrade the $L^2$-height estimate to a weighted $C^{1,\alpha}$-estimate. In particular, we obtain
		\[
		|D_\xi D_y^\beta\varphi(\xi,y)|\lesssim_\beta|\xi|,
		\]
		in accordance with the quadratic normal term and the evenness of $\varphi$. These estimates close the induction for smooth ambient data. For real-analytic data, we track their dependence on the derivative order to obtain the factorial bounds in \eqref{eqn:front-tangential-factorial-bound}, and hence real analyticity of the branching locus.
	\end{enumerate}
	
	\begin{remark}
			Krummel \cite{Kru14} studied real analyticity of the branch set for two-valued Dirichlet minimisers and $C^{1,\alpha}$ two-valued harmonic functions under the assumption that the frequency is locally constant along the branch set. His approach uses quantitative asymptotics at branch points to construct a partial Legendre transform, following Kinderlehrer--Nirenberg--Spruck \cite{KNS78}. The intended route to analyticity is then through Schauder estimates for the resulting singular elliptic system and a majorant argument. Partial Legendre transforms have also been used for the Signorini problem \cite{KPS15}. Our treatment of the transformed problem instead uses excess decay and Campanato iteration. This approach is of independent interest.
	\end{remark}

    \begin{remark}[Frequency jumping set]
        We expect that our method can be extended to show that $\mathcal B_u$ is smooth (resp. real-analytic) near every branch point at which $N_u(\cdot,0)$ is locally constant on $\mathcal B_u$, with value $\frac{1}{2}+k$ for some $k\in\mathbb Z_{\geqslant0}$ (\cf \cite{Kru14}). 
        How large is the set of branch points near which $N_u(\cdot,0)$ is not constant on $\mathcal B_u$? For real-analytic ambient data, does it have vanishing $\mathcal H^{n-2}$-measure, or even Hausdorff dimension at most $n-3$? By contrast, for smooth data, it may have positive $\mathcal H^{n-2}$-measure (\cf the fat Cantor set example by extending the construction of \cite{CH26}).
    \end{remark}
	\subsection{Organisation} Part I begins in Section~\ref{sec:minimality-compactness} with $\mathscr J$-minimality, H\"older regularity, strong blow-up compactness and the frequency-$\frac{1}{2}$ models. Section~\ref{sec:cylindrical-spine-estimates} establishes sheet matching, the no-hole property and the non-concentration estimates. Section~\ref{sec:excess-branch-regularity} establishes linear excess decay and carries out the Campanato iteration, completing the proof of Theorem~\ref{thm:main-tangent-estimates}. In Part II, Section~\ref{sec:legendre-coordinates} introduces the partial Legendre transform and the weighted spaces. Sections~\ref{sec:frozen-operator} and \ref{sec:linearised-polynomial-approximation} prove the frozen and variable-coefficient linear estimates, corresponding to the first two steps above. Section~\ref{sec:tangential-regularity} carries out the tangential regularity induction and proves real analyticity, completing the proof of Theorem~\ref{thm:main-regularity}.
	
	\subsection*{Acknowledgements.} Part I develops and refines material from the second author's Ph.D. thesis \cite{Par22}. Y. L. is supported by the Royal Society URF. D. P. was partially supported by the Cecilia Tanner Research Impulse Fund from Imperial College London. Y. Z. is supported by the London School of Geometry and Number Theory Centre for Doctoral Training (CDT), a joint venture between University College London, Imperial College London and King's College London.
	
	We thank Simon Donaldson, Siqi He, Gregory Parker, Luca Spolaor, Thomas Walpuski, and Neshan Wickramasekera for their interest in this work and for useful discussions. Parts of this work were discussed during the 2026 Clay workshop \textit{Multi-valued Harmonicity} at Oxford. We are also grateful to the organisers and all participants for an inspiring meeting.
	
	\subsection*{AI disclosure.} The main ideas of this paper were developed by the authors. Large language models (LLMs) were used to check routine calculations and assist with manuscript preparation. The authors take full responsibility for the content and correctness of the paper. 

\part{\texorpdfstring{\(C^{1,\alpha}\)}{C1-alpha}-regularity of the branching locus}
	\label{part:branch-regularity}
	
	This part establishes \(C^{1,\alpha}\)-regularity of the branching locus together with the quantitative estimates, which are later used in Part II (\cf Proposition \ref{prop:legendre-initial-estimates}). The principal new ingredient is the integrated Weitzenb\"ock formula in Section~\ref{sec:integrated-weitzenbock}, proved without using the Laplacian. It leads to the local minimising property in Section~\ref{sec:dirac-minimality}. H\"older regularity is then essentially a consequence of Stuvard's work \cite{Stu19}; we briefly discuss the argument for the benefit of the reader. The remainder of Part I is adapted from the second author's thesis \cite{Par22}, which follows the strategy of Simon \cite{Sim93} and Krummel--Wickramasekera \cite{KW13}, and uses the Dirac equation to simplify several steps. 
	
	\section{\texorpdfstring{$\mathscr J$-energy minimality, H\"older regularity and frequency-$\frac{1}{2}$ models}{J-energy minimality, Holder regularity and frequency-one-half models}}
	\label{sec:minimality-compactness}
	
	\subsection{Geometric setting}
	\label{sec:front-geometric-setting}
	
	Recall that $E\to M$ is a real or complex spinor bundle over a Riemannian manifold $(M,g)$, equipped with a fibre metric $h$, Clifford multiplication $\mathfrak c$, and a compatible connection $\nabla^E$, which we denote by $\nabla$ below. The associated Dirac operator is $\mathcal D=\sum_i\mathfrak c(e_i)\nabla_{e_i}$ in a local orthonormal frame $\{e_i\}$. We use the Clifford conventions
	\[
	\mathfrak c(v)\mathfrak c(w)+\mathfrak c(w)\mathfrak c(v)=-2g(v,w), \qquad \mathfrak c(v)^*=-\mathfrak c(v).\]
	The Weitzenb\"ock formula on $C^2$ sections of $E$ is
	\[\mathcal D^*\mathcal D=\nabla^*\nabla+\mathscr R,
    \qquad \mathscr R=\sum_{i<j}\mathfrak c(e_i)\mathfrak c(e_j)R^E_{e_i,e_j},
	\]
	where $R^E$ is the curvature of $\nabla$. We write $\langle v,w\rangle=\operatorname{Re}h(v,w)$. For a two-valued $W^{1,2}$-section $W=\llbracket w_1\rrbracket+\llbracket w_2\rrbracket$, set
	\[
	|W|^2=\frac{1}{2}\sum_{\ell=1}^2|w_\ell|^2,
	\qquad |\nabla W|^2=\frac{1}{2}\sum_{\ell=1}^2|\nabla w_\ell|^2,
	\qquad |\mathcal DW|^2=\frac{1}{2}\sum_{\ell=1}^2|\mathcal Dw_\ell|^2.
	\]
	Here the covariant gradient is defined almost everywhere using the approximate derivative of $W$ in local trivialisations. We measure the pointwise distance between two sections of $\mathcal A_2(E)$ using the Almgren metric \cite{DLS11}, normalised by
	\[
	\mathcal G\!\left(\sum_{\ell=1}^2\llbracket v_\ell\rrbracket, \sum_{\ell=1}^2\llbracket w_\ell\rrbracket\right)^2=\frac{1}{2}\min_{\pi\in S_2}\sum_{\ell=1}^2|v_\ell-w_{\pi(\ell)}|^2,
	\]
	so $|W|^2=\mathcal G(W,2\llbracket0\rrbracket)^2$. The same convention applies to covariant derivatives in $T^*M\otimes E$ and higher covariant tensors, with one permutation matching all components of each sheet.
	
	On a geodesic ball $B_r(p)\Subset M$, define the energy of $W\in W^{1,2}(B_r(p);\mathcal A_2(E))$ by
	\[
	\mathscr J_{B_r(p)}(W)=\frac{1}{2}\sum_{\ell=1}^2\int_{B_r(p)}(|\nabla w_\ell|^2+\langle\mathscr Rw_\ell,w_\ell\rangle)\,dV_g.
	\]
	This is independent of the ordering of the sheets. For $p\in Z_u$, set
	\[
	H_u(p,r)=\int_{\partial B_r(p)}|u|^2\,d\mathrm{Area}_g,\qquad N_u(p,r)=\frac{r\mathscr J_{B_r(p)}(u)}{H_u(p,r)},\qquad N_u(p,0)=\lim_{r\downarrow0}N_u(p,r).
	\]
	The existence of this limit is proved in Section~\ref{sec:frequency-control}.
	
	The tangent models are described in Subsection~\ref{sec:half-frequency-models}. We compare them in geodesic coordinates and the \emph{exponential gauge centred at $p$}: an orthonormal basis of $E_p$ is extended by parallel transport along radial geodesics. This identifies $E_x$ with $E_p$. We write $X=\exp_p^{-1}x$ and suppress the parallel transport maps below.

	\subsection{The integrated Weitzenböck formula}
	\label{sec:integrated-weitzenbock}
    Although our applications concern two-valued sections, the following integrated Weitzenböck formula holds for arbitrary $Q$-valued $W^{1,2}$-sections. The conventions in the previous Section \ref{sec:front-geometric-setting} extend to $\mathcal A_Q(E)$, with quantities averaged over the $Q$ sheets. 
  
	Let $\pi:E\to M$ be the bundle projection. Then the pullback bundle $\pi^*E\to E$ admits a tautological section $\tau(e)=e$; we write
	\[
    \vartheta=\nabla^{\pi^*E}\tau\in\Omega^1(E;\pi^*E).
	\]
	The form $\vartheta$ vanishes on horizontal vectors and identifies each vertical tangent space $T_e(E_{\pi(e)})$ with $E_{\pi(e)}$. For a smooth section $v:\Omega\to E$ over an open set $\Omega\subset M$, the natural identification $v^*(\pi^*E)=E|_\Omega$ implies
    \[
    v^*\vartheta=\nabla v\in\Omega^1(\Omega;E).
	\]
	In a local oriented orthonormal frame $\{e_i\}_{i=1}^n$ of $TM$, set
	\[
	\Theta=\sum_{i<j}\langle\tau,\mathfrak c(e_i)\mathfrak c(e_j)\vartheta\rangle\wedge\pi^*(\iota_{e_j}\iota_{e_i}dV_g)\in\Omega^{n-1}(E;\mathbb R).
	\]
	Here the Clifford endomorphisms act on the pullback bundle. Thus $\Theta$ is a real differential form on the total space $E$, with horizontal degree $n-2$ and vertical degree one relative to the connection.
	
	\begin{proposition}[Integrated Weitzenböck formula]
		\label{prop:integrated-weitzenbock}
		For every \(W=\sum_{\ell=1}^Q\llbracket w_\ell\rrbracket\in W^{1,2}(\Omega;\mathcal A_Q(E))\) and every \(\chi\in C_c^\infty(\Omega)\),
		\begin{equation}
			\begin{aligned}
				\sum_{\ell=1}^Q\int_\Omega\chi(|\nabla w_\ell|^2+\langle\mathscr Rw_\ell,w_\ell\rangle-|\mathcal Dw_\ell|^2)\,dV_g =\sum_{\ell=1}^Q\int_\Omega(\langle\mathcal Dw_\ell,\mathfrak c(d\chi)w_\ell\rangle-\langle\nabla_{\nabla \chi}w_\ell,w_\ell\rangle)\,dV_g.
			\end{aligned}
			\label{eqn:integrated-weitzenbock}
		\end{equation}
	\end{proposition}
	
	\begin{proof} We establish \eqref{eqn:integrated-weitzenbock} successively for smooth, Lipschitz, and $W^{1,2}$-sections. Note that the definition of $\Theta$ is independent of the orthonormal frame. For a smooth section $v$, the Clifford relation implies the pointwise identity
		\[
		|\mathcal Dv|^2=|\nabla v|^2-2\sum_{i<j}\langle\nabla_{e_i}v,\mathfrak c(e_i)\mathfrak c(e_j)\nabla_{e_j}v\rangle.
		\]
		Computing in local normal coordinates,
		$$d\langle v,\mathfrak c(e_i)\mathfrak c(e_j)\nabla v\rangle(e_i,e_j)=2\langle\nabla_{e_i}v, \mathfrak c(e_i)\mathfrak c(e_j)\nabla_{e_j}v\rangle +\langle v,\mathfrak c(e_i)\mathfrak c(e_j)R^E_{e_i,e_j}v\rangle.$$
		Summing over $i<j$ yields the smooth transgression formula 
		\[
		v^*d\Theta=(|\nabla v|^2+\langle\mathscr Rv,v\rangle-|\mathcal Dv|^2)\,dV_g.
		\]
		Thus $v^*d\Theta$ is the divergence term in the Weitzenböck formula. 
		
		Now, let $W$ be a Lipschitz $Q$-valued section, and let $\Gamma_W$ be the current associated with its graph in $E|_\Omega$, counting each sheet with multiplicity one. The boundary theorem \cite[Theorem 2.1]{DLS15} gives $\partial\Gamma_W=0$ in $E|_\Omega$. Since $W$ is bounded on $\operatorname{spt}\chi$, the support of $\pi^*\chi\,\Theta$ meets the graph in a compact set. Stokes' theorem therefore yields
		\[
		0=\partial\Gamma_W(\pi^*\chi\,\Theta)=\Gamma_W(d(\pi^*\chi\,\Theta))=\sum_{\ell=1}^Q\int_\Omega(\chi\,w_\ell^*d\Theta+d\chi\wedge w_\ell^*\Theta).
		\]
		The pullbacks $w_\ell^*\Theta$ and $w_\ell^*d\Theta$ are interpreted using the approximate differentials of $W$ and do not require globally defined Lipschitz sheets. It follows from the Clifford relation that
		\[
		-d\chi\wedge w_\ell^*\Theta=(\langle\mathcal Dw_\ell,\mathfrak c(d\chi)w_\ell\rangle-\langle\nabla_{\nabla\chi}w_\ell,w_\ell\rangle)dV_g, \quad\text{for almost every point in }\Omega.
		\]
		Dividing the sum over the sheets by $Q$ proves \eqref{eqn:integrated-weitzenbock} for Lipschitz $W$.
		
		Finally, for a general $W^{1,2}$-section, choose Lipschitz approximations $W_k$ converging strongly to $W$ in $W^{1,2}$, as in \cite[Proposition 4.4]{DLS11}. The approximations may be chosen with $\lambda_k\to\infty$, $\operatorname{Lip}W_k\leqslant C\lambda_k$, and
		\[
		A_k=\{W_k\ne W\},\qquad |A_k|=o(\lambda_k^{-2}),\qquad\|\mathcal G(W_k,W)\|_{L^2(\Omega)}\rightarrow0.
		\]
		Absolute continuity of the integral and the Lipschitz bound imply
		\[
		\varepsilon_k:=\int_{A_k}(|\nabla W_k|^2+|\nabla W|^2+|W_k|^2+|W|^2)\,dV_g
		\rightarrow0.
		\]
		Outside $A_k$, the sections and their covariant derivatives agree almost everywhere. With the $L^1$-norms taken over $\operatorname{spt}\chi$, Clifford contraction and the boundedness of the curvature give
		\[
		\||\nabla W_k|^2-|\nabla W|^2\|_{L^1}+\||\mathcal DW_k|^2-|\mathcal DW|^2\|_{L^1}+\|\langle\mathscr RW_k,W_k\rangle-\langle\mathscr RW,W\rangle\|_{L^1}\lesssim\varepsilon_k\rightarrow0.
		\]
		Denote the integrand on the right-hand side of \eqref{eqn:integrated-weitzenbock} by $b_\chi$. Over $A_k$, the Cauchy--Schwarz inequality implies that 
		\[
		\begin{aligned}
			\|b_\chi(W_k)-b_\chi(W)\|_{L^1}&\lesssim\|d\chi\|_\infty\left(\int_{A_k}(|\nabla W_k|^2+|\nabla W|^2)\,dV_g\right)^{1/2}\left(\int_{A_k}(|W_k|^2+|W|^2)\,dV_g\right)^{1/2}
            \\&\lesssim_\chi \varepsilon_k\rightarrow0,\quad k\to \infty.
		\end{aligned}
		\]
		Both sides of \eqref{eqn:integrated-weitzenbock} therefore converge in $L^1$. The integrated Weitzenböck formula follows. 
	\end{proof}

	\subsection{\texorpdfstring{Local minimality of $\mathscr J$-energy for $\mathbb Z_2$-harmonic spinors}{Local minimality of J-energy for Z2-harmonic spinors}}
	\label{sec:dirac-minimality}
	
	\begin{corollary}[Dirac solutions are $\mathscr J$-energy local minimisers]
		\label{cor:dirac-minimality}
		Let $u$ be a $\mathbb Z_2$-harmonic spinor.  For every ball $B\Subset M$ and every $V\in W^{1,2}(B;\mathcal A_2(E))$ with the same trace as $u$ on $\partial B$, we have
		\[
		\mathscr J_B(V)-\mathscr J_B(u)=\|\mathcal DV\|_{L^2(B)}^2\geqslant0.
		\]
		
	\end{corollary}

	\begin{proof}
		Choose a slightly larger ball $\widetilde B$ with $B\Subset\widetilde B\Subset M$, and set
		\[
		\widetilde V=
		\begin{cases}
			V&\text{in }B,\\
			u&\text{in }\widetilde B\setminus B.
		\end{cases}
		\]
		The equality of the traces implies $\widetilde V\in W^{1,2}(\widetilde B;\mathcal A_2(E))$: in local trivialisations, the traces agree after the Almgren embedding \cite[Proposition 2.10]{DLS11}, so the ordinary Sobolev gluing theorem applies; using the Lipschitz inverse on the embedded target, we then obtain the asserted two-valued Sobolev section. Choose a cutoff 
		\[
		\chi\in C_c^\infty(\widetilde B),\qquad 0\leqslant\chi\leqslant1, \qquad \chi=1\quad\text{on a neighbourhood of }\overline B.
		\]
		In particular, $\operatorname{spt}d\chi\subset\widetilde B\setminus \overline B$, where $\widetilde V=u$ and their approximate differentials agree almost everywhere. Applying \eqref{eqn:integrated-weitzenbock} to $\widetilde V$ and $u$ and taking differences, we find
		\begin{equation*}
			0=\mathscr J_B(V)-\mathscr J_B(u)-\|\mathcal DV\|_{L^2(B)}^2+\|\mathcal Du\|_{L^2(B)}^2\leqslant \mathscr J_B(V)-\mathscr J_B(u),
		\end{equation*}
		as $\mathcal{D}u=0$. 
	\end{proof}
	
	Minimality implies the outer and inner variation identities recorded below.
	
	\begin{remark}[Deviation from flat spaces]\label{rmk:rescaled-geometry}
		In geodesic coordinates and the exponential gauge, rescale a ball of radius $r$ to $B_1$. The rescaled $\mathscr J$-energy density is
		\[
		\frac{1}{2}|dv|^2+\frac{1}{2}a^{ij}_{AB}\partial_i v^A\partial_j v^B+b^i_{AB}v^A\partial_i v^B+\frac{1}{2}c_{AB}v^Av^B,
		\]
		where
		\[
		|a|\lesssim r^2|x|^2,\qquad|da|+|b|\lesssim r^2|x|,\qquad|db|+|c|+|dc|\lesssim r^2.
		\]
		For a local minimiser $u$ of this energy, the outer variation $(1+t\zeta)u$, with $\zeta\in C_c^\infty(B_1)$, gives
		\[
		\int_{B_1}(\zeta|du|^2+\langle du,d\zeta\otimes u\rangle)
		=O\!\left(r^2\|\zeta\|_{C^1}\int_{B_1}(|\nabla u|^2+|u|^2)\right).
		\]
		For $Y\in C_c^\infty(B_1;\mathbb R^n)$, the inner variation $u\circ(\operatorname{id}-tY)$ similarly gives
		\[
		\int_{B_1}\left(\tfrac{1}{2}|du|^2\operatorname{div}Y-\langle\partial_i u,\partial_j u\rangle\partial_iY_j\right)=O\!\left(r^2\|Y\|_{C^1}\int_{B_1}(|\nabla u|^2+|u|^2)\right).
		\]
	\end{remark}
	
	\begin{corollary}[Caccioppoli estimate]\label{cor:caccioppoli-estimate}
		Let $W$ be a local minimiser of the $\mathscr J$-energy in $B_t(x_0)$. For concentric balls $B_s(x_0)\subset B_t(x_0)$, where $0<s<t<1$, we have
		\begin{equation}
			\label{eqn:caccioppoli-estimate}
			\int_{B_s(x_0)}|\nabla W|^2
			\leqslant\frac{C}{(t-s)^2}\int_{B_t(x_0)}|W|^2.
		\end{equation}
	\end{corollary}
	\begin{proof}
    We simply sketch the proof as it is standard. Choose $\zeta\in C_c^\infty(B_t(x_0))$ with $\zeta=1$ on $B_s(x_0)$,$0\leqslant\zeta\leqslant1$, and $|d\zeta|\lesssim(t-s)^{-1}$. The outer variation $W_\tau=(1+\tau\zeta^2)W$ and Young's inequality give
		\[
		\int_{B_t}\zeta^2|\nabla W|^2 \lesssim\int_{B_t}(|d\zeta|^2+\zeta^2)|W|^2 \lesssim\frac{1}{(t-s)^2}\int_{B_t}|W|^2.
		\]
	\end{proof}
	
	Following \cite[Theorem 3.9 and Proposition 3.10]{DLS11} and \cite[Theorem 5.1 and Proposition 5.2]{Stu19}, we first obtain a uniform Hölder estimate and frequency control for $\mathbb Z_2$-harmonic spinors, then prove strong \(W^{1,2}\)-convergence of blow-ups. The former is a standard consequence of $\mathscr{J}$-minimality and Campanato–Morrey iteration, while the latter follows by applying the integrated Weitzenböck formula to the rescaled spinors.

	\subsection{A uniform Hölder estimate}
	\label{sec:holder-regularity}
	
	In a smooth metric trivialisation over $B_s(x_0)$, the $\mathscr J$-energy has the form
	\begin{equation}
		\mathscr J_{B_s(x_0)}(W)=\frac{1}{2}\sum_{\ell=1}^2\int_{B_s(x_0)}
		\left(
		a^{ij}_{\alpha\beta}\partial_iw_\ell^\alpha \partial_jw_\ell^\beta
		+2b^i_{\alpha\beta}w_\ell^\alpha\partial_iw_\ell^\beta
		+c_{\alpha\beta}w_\ell^\alpha w_\ell^\beta
		\right)dx,
		\label{eqn:jacobi-functional}
	\end{equation}
	where $a^{ij}_{\alpha\beta}=\sqrt{\det g}\,g^{ij}h_{\alpha\beta}$ is uniformly positive, and $a,b,c$ are smooth and bounded. The Euclidean Dirichlet energy is $
	\operatorname{Dir}(W,B_s(x_0))
	=\int_{B_s(x_0)}|dW|^2
	=\frac{1}{2}\sum_{\ell=1}^2\int_{B_s(x_0)}|dw_\ell|^2,$
	with $\operatorname{Dir}(W,\partial B_s(x_0))$ denoting the tangential energy on the boundary.
	
	\begin{proposition}[Hölder regularity for local minimisers]
		\label{prop:holder-regularity}
		
		Let \(W\) be a local minimiser of the functional $\mathscr J$. There exist $\gamma>0$, $\Lambda,C<\infty$, and $s_0>0$, depending only on the ellipticity and coefficient bounds, such that, if $W$ is defined in $B_{s_0}(x_0)$, then for $0<s<t<s_0$,
		\begin{equation}
			\operatorname{Dir}(W,B_s)
			\leqslant C\left(\frac{s}{t}\right)^{n-2+\gamma}
			\left(
			\operatorname{Dir}(W,B_t)+\Lambda\int_{B_t}|W|^2
			\right).
			\label{eqn:energy-morrey-decay}
		\end{equation}
		In particular, $W\in C^{\gamma/2}_{\mathrm{loc}}$. 
		
	\end{proposition}
	
	\begin{proof}
		We use geodesic coordinates centred at $x_0$ with radius $s_0\ll1$ and an orthonormal bundle frame. The following estimates are uniform for centres in a fixed compact set of $M$. As a competitor, take a Euclidean Dirichlet minimiser $f$ with the same boundary data as $W$ on $\partial B_s$. Since $|a-I|\lesssim s$ and the lower-order coefficients are bounded, minimality and Young's inequality give, for any $\varepsilon>0$ and sufficiently small $s_0$,
		\begin{equation}
			\operatorname{Dir}(f,B_s)\leqslant\operatorname{Dir}(W,B_s)
			\leqslant(1+\varepsilon)\operatorname{Dir}(f,B_s)
			+C_\varepsilon\int_{B_s}(|W|^2+|f|^2),\qquad 0<s<s_0.
			\label{eqn:harmonic-competitor-comparison}
		\end{equation}

		By \cite[Proposition 3.10]{DLS11} and \cite[Corollary 5.4]{Stu19}, we may fix $\gamma\in(0,1)$ sufficiently small so that
		\begin{equation}
			\operatorname{Dir}(f,B_s)
			\leqslant\frac{s}{n-2+2\gamma}\operatorname{Dir}(W,\partial B_s),
			\label{eqn:harmonic-boundary-energy}
		\end{equation}
		\begin{equation}
			\int_{B_s}|v|^2
			\leqslant\frac{s}{n-2+2\gamma}\int_{\partial B_s}|v|^2
			+Cs^2\operatorname{Dir}(v,B_s),\qquad v=W,f.
			\label{eqn:scalar-trace-estimate}
		\end{equation}
		To absorb the lower-order terms, we follow \cite[Proposition 5.2]{Stu19} and add $\Lambda\int_{B_s}|W|^2$ to the Dirichlet energy. Choose $\varepsilon$ small, then $\Lambda$ large, and finally decrease $s_0$. Combining these estimates and absorbing the terms proportional to $s^2\operatorname{Dir}(W,B_s)$ gives, for almost every $0<s<s_0$,
		\begin{equation}
			\operatorname{Dir}(W,B_s)+\Lambda\int_{B_s}|W|^2
			\leqslant\frac{s}{n-2+\gamma}\left(
			\operatorname{Dir}(W,\partial B_s)
			+\Lambda\int_{\partial B_s}|W|^2\right).
			\label{eqn:energy-boundary-comparison}
		\end{equation}
		With $\varphi(r)=\operatorname{Dir}(W,B_r)+\Lambda\int_{B_r}|W|^2$, the coarea formula and \eqref{eqn:energy-boundary-comparison} imply
		\[
		r\varphi'(r)\geqslant (n-2+\gamma)\varphi(r) \quad\text{for almost every }r \in (0,s_0).
		\]
		Integrating from $s$ to $t$ gives \eqref{eqn:energy-morrey-decay}. The local $C^{\gamma/2}$-regularity follows from the Campanato--Morrey estimate.
	\end{proof}

	\subsection{Frequency control}
	\label{sec:frequency-control}
	
	We collect the frequency estimates needed to construct nonzero homogeneous blow-ups. For $\mathbb Z_2$-harmonic spinors, the underlying identities and curvature estimates appear in \cite[Lemmas 3.1--3.3]{Tau14} and \cite[Section 3.2]{Par22}; the corresponding variational argument for multiple-valued Jacobi fields is developed in \cite[Lemma 6.13 and Propositions 6.9--6.10]{Stu19}.
	
	Let $u$ be a symmetric local $\mathscr J$-minimiser and $p\in\mathcal B_u$. Suppressing the centre $p$, write
	\[
	H(r)=H_u(p,r),\qquad
	\mathscr J(r)=\mathscr J_{B_r(p)}(u),\qquad
	N(r)=\frac{r\mathscr J(r)}{H(r)},\qquad
	h(r)=r^{1-n}H(r).
	\]
	For sufficiently small $r_0>0$, we have $H(r)>0$ and $N(r)\geqslant-Cr^2$ for $0<r\leqslant r_0$. The key ingredients are the inner and outer variation identities in Remark~\ref{rmk:rescaled-geometry}. Applying the outer variation to radial cutoffs, we get
	\begin{equation}
		\mathscr J(r)=\int_{\partial B_r(p)}\langle u,\nabla_\nu u\rangle.
		\label{eqn:dirac-boundary-flux}
	\end{equation}
	Differentiating the height and using radial inner variations give, respectively,
	\begin{align}
		H'(r)&=\frac{n-1}{r}H(r)+2\mathscr J(r)+O(rH(r)),\label{eqn:height-variation}\\
		\mathscr J'(r)&=\frac{n-2}{r}\mathscr J(r)+2\int_{\partial B_r(p)}|\nabla_\nu u|^2+O\!\left(r\int_{B_r(p)}|\nabla u|^2+H(r)\right).
		\label{eqn:energy-radial-variation}
	\end{align}
	These identities hold for almost every $0<r<r_0$ and imply
	\begin{equation}
		N'(r)\geqslant\frac{2r}{H(r)}\int_{\partial B_r(p)}\left|\nabla_\nu u-\frac{N(r)}{r}u\right|^2-Cr(1+N(r)).
		\label{eqn:frequency-derivative}
	\end{equation}
	Thus $e^{Ar^2}(1+N(r))$ is nondecreasing for sufficiently large $A$, and the limit $\nu_p:=N_u(p,0)=\lim_{r\downarrow0}N(r)\in[0,\infty)$ exists. Integrating \eqref{eqn:height-variation} gives, for $0<s<t\leqslant r_0$,
	\begin{equation}
		\left|\log\frac{h(t)}{h(s)}-2\int_s^t\frac{N(\tau)}{\tau}\,d\tau\right|
		\leqslant C(t^2-s^2).
		\label{eqn:frequency-height-growth}
	\end{equation}
	Here $r_0$, $A$ and $C$ depend only on the local geometry and the coefficients of $\mathscr J$.

	\subsection{Strong \texorpdfstring{\(W^{1,2}\)}{W1,2}-convergence of blow-ups}
	\label{sec:strong-blowup-convergence}
	In the exponential gauge centred at a branch point $p$, define the normalised blow-ups by $u_r(X)=h(r)^{-1/2}u(\exp_p(rX))$. 
	
	\begin{theorem}[Blow-ups of $\mathbb Z_2$-harmonic spinors]
		\label{thm:spinor-blowups}
		Let $u$ be a $\mathbb Z_2$-harmonic spinor, and let $p\in\mathcal B_u$. For any sequence $r_j\downarrow0$, the normalised blow-ups $u_{r_j}$ admit a subsequence converging strongly in $W^{1,2}_{\mathrm{loc}}$ to a nonzero limit $\Phi$ on $T_pM$.
		
		The limit is locally Dirichlet-minimising, is homogeneous of degree $N_u(p,0)$, and satisfies $\mathcal D_0\Phi=0$, where $\mathcal D_0$ is the constant-coefficient Dirac operator obtained by freezing the principal symbol of $\mathcal{D}$ at $p$. 
	\end{theorem}
	
	\begin{proof}
		Write $u_j=u_{r_j}$ and $g_j=r_j^{-2}(\exp_p(r_j\,\cdot))^*g$. We first establish the $W^{1,2}$-weak convergence using the frequency function and then apply the integrated Weitzenböck formula \eqref{prop:integrated-weitzenbock} to derive minimality of the limit and energy convergence. 
		
		 For each fixed $R\geqslant1$, the frequency bounds and \eqref{eqn:frequency-height-growth} give $h(\rho r_j)\leqslant C_Rh(r_j)$ for $0<\rho\leqslant2R$ and $r_j<r_0/(2R)$. Here $C_R$ depends also on $N_u(p,r_0)$, but not on $j$. By the coarea formula,
		\[
		\int_{B_{2R}}|u_j|^2\,dV_{g_j}
        =\int_0^{2R}\rho^{n-1}\frac{h(\rho r_j)}{h(r_j)}\,d\rho \leqslant C_R.
		\]
		The Caccioppoli estimate \eqref{eqn:caccioppoli-estimate} and uniform equivalence of the rescaled metrics then yield
		\begin{equation}
			\int_{B_R}(|u_j|^2+|du_j|^2)\,dX\leqslant C_R.
			\label{eqn:normalised-blowup-bound}
		\end{equation}
		Proposition~\ref{prop:holder-regularity} also gives uniform H\"older bounds on each fixed ball. Passing to a diagonal subsequence, we can assume
		\begin{equation}\label{eqn:weak-blowup-convergence}
			u_j\to\Phi\quad\text{strongly in }L^2_{\mathrm{loc}}\cap C^{\gamma'/2}_{\mathrm{loc}},
			\qquad
			u_j\rightharpoonup\Phi\quad\text{in }W^{1,2}_{\mathrm{loc}},\quad \gamma'<\gamma,
		\end{equation}
		with $\Phi(0)=2\llbracket0\rrbracket$. Passing to the limit in the weak formulation of the rescaled Dirac equations, using \eqref{eqn:weak-blowup-convergence} and convergence of the coefficients, gives $\mathcal D_0\Phi=0$. The integrated Weitzenb\"ock formula, as in Corollary~\ref{cor:dirac-minimality}, then shows that $\Phi$ is locally Dirichlet-minimising. The usual homogeneity argument from \eqref{eqn:frequency-derivative}, using strong $L^2_{\mathrm{loc}}$-convergence, shows that $\Phi$ is homogeneous of degree $N_u(p,0)$. The normalisation then gives
		\begin{equation}
			\int_{\partial B_\rho}|\Phi|^2\,d\mathrm{Area}=\rho^{n-1+2N_u(p,0)}>0,\qquad \rho>0.
			\label{eqn:normalised-blowup-height}
		\end{equation}
		
		Now, for $\chi\in C_c^\infty(\mathbb R^n)$, apply \eqref{eqn:integrated-weitzenbock} to $u_j$ in the rescaled geometry and use its Dirac equation. Smooth convergence of the rescaled coefficients and the bound \eqref{eqn:normalised-blowup-bound} give
		\[
		\begin{aligned}
		\int\chi|du_j|^2\,dX&=\frac{1}{2}\int(\Delta\chi)|u_j|^2\,dX+o(1)
        \\&\rightarrow\frac{1}{2}\int(\Delta\chi)|\Phi|^2\,dX=\int\chi|d\Phi|^2\,dX.
		\end{aligned}
		\]
		Here $\Delta=\sum_i\partial_i^2$, the limit follows from strong $L^2_{\mathrm{loc}}$-convergence, and the last equality follows from the same formula and $\mathcal D_0\Phi=0$. Approximating ball indicators by smooth cutoffs yields
		\begin{equation}\label{eqn:blowup-energy-convergence}
		\lim_{j\to\infty}\operatorname{Dir}(u_j,B_R)=\operatorname{Dir}(\Phi,B_R),\qquad R>0,
		\end{equation}
		and hence strong $W^{1,2}_{\mathrm{loc}}$-convergence.
	\end{proof}
	
	\subsection{\texorpdfstring{Frequency-\(\frac{1}{2}\)}{Frequency-one-half} models}
	\label{sec:half-frequency-models}
	
	For any codimension-two plane $P\subset\mathbb R^n$, we decompose $\R^n$ as $\mathbb R^n=P^\perp\oplus P\cong \mathbb C_z\oplus P$. If $E$ is real, we complexify $E$ and extend Clifford multiplication complex linearly.
	
	\begin{definition}[Frequency-\(\frac{1}{2}\) models]
		\label{def:half-frequency-models}
		
		Let \(\mathcal C_{\mathcal D}\) be the class of nonzero maps
		\[\Phi_{P,a,b}(z):=\llbracket az^{1/2}+b\bar z^{1/2}\rrbracket
		+\llbracket-az^{1/2}-b\bar z^{1/2}\rrbracket,\quad z\in \mathbb C,
		\]
		where
		\[Ja=-ia,\qquad Jb=ib,\qquad J=\mathfrak c(e_1)\mathfrak c(e_2),
		\]
		for an oriented orthonormal basis $e_1,e_2$ of $P^\perp$. The model is translation invariant along $P$. The coefficients lie in $E_p$ in the complex case and in its complexification in the real case. For a real spinor bundle $E$, we impose in addition $b=\bar a$, where the bar denotes the natural conjugation on the complexification of $E_p$. 
	\end{definition}
	
	\begin{proposition}[Minimum frequency and tangent classification]
		\label{prop:minimum-frequency-classification}
		
		Let $u$ be a nonzero $\mathbb Z_2$-harmonic spinor and $p\in \mathcal{B}_u$. Then
		\[
		N_u(p,0)\geqslant\frac{1}{2}.
		\]
		If $N_u(p,0)=\frac{1}{2}$, every tangent map at $p$ belongs to $\mathcal C_{\mathcal D}$. 
		
	\end{proposition}
	
	\begin{proof}
		The dimension-reduction argument in Taubes \cite[Lemmas 5.1, 5.4, 5.5 and 6.4]{Tau14} applies in arbitrary dimension. We recall the main points. By Theorem~\ref{thm:spinor-blowups}, every tangent map $\Phi$ is a nonzero homogeneous Dirichlet minimiser satisfying $\mathcal D_0\Phi=0$, of positive degree $\alpha=N_u(p,0)$. In dimension two, the angular Dirac equation gives $\alpha\in\frac{1}{2}\mathbb Z_{\geqslant1}$. For $n\geqslant3$, if $Z_\Phi=\{0\}$, a sheet extends smoothly across the origin, so $\alpha$ is a positive integer. Otherwise, at a zero $q\ne0$, the cylindrical blow-up lemma \cite[Lemma 3.24]{DLS11} gives a tangent map invariant in the direction $q$. Induction on the dimension, together with frequency monotonicity and $N_\Phi(q,R)\to\alpha$ as $R\to\infty$, yields
		\[
		\frac{1}{2}\leqslant N_\Phi(q,0)\leqslant N_\Phi(q,R)\leqslant\alpha.
		\]
		If $\alpha=\frac{1}{2}$, equality holds throughout. Thus $\Phi$ is homogeneous about both $0$ and $q$, hence invariant in the direction $q$. Iterating reduces to the two-dimensional angular equation and gives a sheet
		\begin{equation}
			\phi(z)=az^{1/2}+b\bar z^{1/2},\qquad Ja=-ia,\quad Jb=ib,\quad |\Phi(z)|^2=(|a|^2+|b|^2)|z|,
			\label{eqn:half-frequency-model}
		\end{equation}
		translation invariant along an $(n-2)$-plane. In the real case $b=\bar a$. Hence $\Phi\in\mathcal C_{\mathcal D}$.
	\end{proof}
	
	\begin{remark}[Model gauge]
		\label{rmk:model-gauge}
		On the branched double cover $z=w^2$, a tangent model is $\Phi=\pm T(w)$ (\cf Definition~\ref{def:half-frequency-models}), where $T:\C \to E_p$, $w\mapsto aw+b\bar w$, is real-linear. This parametrisation is unique \textit{up to rotating the oriented orthonormal frame $(e_1,e_2)$ of $P^\perp$} (which changes $z=x_1+ix_2$) \textit{and the sign of the square root $w$}: namely the model is invariant under $w\mapsto e^{i\phi}w$, $(a,b)\mapsto(e^{-i\phi}a,e^{i\phi}b)$.
		
		We fix this freedom near a given model $\Phi_0=\pm T_0(w)$ as follows. Choose a real-linear projection $\pi:E_p\to\C$ such that $\pi T_0$, regarded as a real $2\times2$ matrix, is invertible with positive determinant. By continuity, so is $\pi T$ for every model $\Phi$ near $\Phi_0$ with the same centre, and by the polar decomposition there is a unique $\phi\in [0,2\pi)$ such that $\pi T\circ e^{-i\phi}$ is symmetric positive definite. We call the corresponding normal frame and square root \textit{the model gauge}. In this gauge the parametrisation is unique, and we can compare nearby models by
		\[
		 d_{\mathcal C}(\Phi,\Phi')^2=\sum_{j=1}^{2}|e_j-e_j'|^2+|a-a'|^2+|b-b'|^2 ,
		\]
		where $(e_1,e_2)$ and $(e_1',e_2')$ are the normal frames of $\Phi$ and $\Phi'$ in the model gauge; the first term also measures the distance between their spines. A compact family of models in such a neighbourhood is called a \emph{non-degenerate model chart}.
	\end{remark}
	
	\medskip
	We record the infinitesimal variations needed in the excess-decay argument.
	\begin{lemma}[Variations of Dirac models]
		\label{lem:dirac-model-variations}
		Fix $\Phi=\{\pm\varphi\}\in\mathcal C_{\mathcal D}$ with $\varphi=az^{1/2}+b\bar z^{1/2}$ and spine $P$. Write $z=x_1+ix_2$ and choose orthonormal coordinates $(y_1,\ldots,y_{n-2})$ on $P$.
		The tangent space $T_\Phi\mathcal C_{\mathcal D}$ consists of the fields
		\begin{equation}
		\ell=Az^{1/2}+B\bar z^{1/2}+\sum_{\alpha=1}^{n-2}\sum_{i=1}^2c_{\alpha i}y_\alpha\partial_i\varphi,
		\label{eqn:dirac-model-tangent-field}
		\end{equation}
		where $A,B$ are fibre constants and $c_{\alpha i}\in\mathbb R$, satisfying $\mathcal D_0\ell=0$ on $\mathbb R^n\setminus P$; in the real case, $B=\bar A$. Every such field is realised by a smooth one-parameter family $\Phi_s\in\mathcal C_{\mathcal D}$ with $\Phi_0=\Phi$ and
	   \begin{equation}
		\|\mathcal G(\Phi_s,\{\pm(\varphi+s\ell)\})\|_{L^2(B_1)}=O(|s|^{3/2}).
		\label{eqn:model-variation-expansion}
		\end{equation}
	\end{lemma}
	\begin{proof}[Sketch of the proof]
		Varying the coefficients gives the first two terms in \eqref{eqn:dirac-model-tangent-field}; tilting the normal plane gives the terms $c_{\alpha i}y_\alpha\partial_i\varphi$.
		The Dirac constraints are $(\mathfrak c(e_1)+i\mathfrak c(e_2))a=0$ and $(\mathfrak c(e_1)-i\mathfrak c(e_2))b=0$. Differentiating these constraints restricts the fields in \eqref{eqn:dirac-model-tangent-field} to those satisfying $\mathcal D_0\ell=0$. Since the corresponding eigenspaces of $J$ vary smoothly with the normal plane, every such infinitesimal variation is realised by a smooth one-parameter family in $\mathcal C_{\mathcal D}$. The error estimate \eqref{eqn:model-variation-expansion} follows from Taylor expansion.
	\end{proof}
	\section{Estimates near a cylindrical spine}
	\label{sec:cylindrical-spine-estimates}
	This section gives estimates for spinors close to cylindrical models of frequency $1/2$, beginning with sheet matching and the topological exclusion of holes. We adapt Simon's cylindrical cone technique, as developed further by Krummel--Wickramasekera \cite{Sim93,KW13}, taking advantage of some simplifications from the Dirac equation. See also the special case of triple junctions \cite{Min24}. 
	
	Given a $\mathbb Z_2$-harmonic spinor on a geodesic ball $B_r(0)$, dilate the ball to $B_1$ and multiply the spinor values by $r^{-1/2}$. Throughout this section, write $u$ for the resulting spinor on $B_1$ and $g$, $\nabla$, $\mathscr R$ for the rescaled ambient data, with $\delta=Cr^2$ as in Remark~\ref{rmk:rescaled-geometry}.
	We use geodesic coordinates centred at $0$, identifying fibres by parallel transport along radial geodesics.
    
	Write
	\[
	X=(z,y)\in\mathbb C\times\mathbb R^{n-2},\qquad P=\{z=0\},\qquad \rho=|z|,\qquad R=|X|.
	\]
	Unless indicated otherwise, balls and volume integrals use $g$. We write $D$ for coordinate differentiation and $dX$ for coordinate measure, and set $z=x_1+ix_2$. Throughout this section, $\Phi=\{\pm\varphi\}$ belongs to a fixed compact nondegenerate Dirac model chart $\mathcal C_0\Subset\mathcal C_{\mathcal D}$, with spine $P$. Its $1/2$-homogeneity makes $\Phi$ invariant under the dilation. Small $L^2$-excess relative to this family gives a uniform $L^2$-bound on $u$.
	
	\subsection{Sheet matching away from the branching locus}
	\label{sec:sheet-matching}
	Away from the model spine, the H\"older estimate and small $L^2$-excess allow us to match the local sheets of $u$ and $\Phi$ and apply elliptic estimates to their difference.
	
	\begin{lemma}[Sheet matching]
		\label{lem:sheet-matching}
		Let $\Phi\in\mathcal C_0$ have spine $P$, and let $K$ be a compact set in $B_{3/4}$ at distance at least $c_0>0$ from both $P$ and $\partial B_{3/4}$. There exist $\varepsilon,C>0$ and $\kappa_0\in(0,1)$ such that, if
		\[
		\int_{B_1}\mathcal G(u,\Phi)^2+\delta<\varepsilon^2,
		\]
		then $u$ has no zero near $K$ and has a unique matching with $\Phi$ such that $\mathcal G(u,\Phi)^2=|u-\varphi|^2$ and
		\begin{equation}
			\|u-\varphi\|_{C^1(K)}\leqslant C\left(\|\mathcal G(u,\Phi)\|_{L^2(B_1)}^{\kappa_0}+\delta\right).
			\label{eqn:sheet-matching-estimate}
		\end{equation}
		The constants are uniform for $\Phi\in\mathcal C_0$.
	\end{lemma}
	\begin{proof}
	Corollary~\ref{cor:caccioppoli-estimate} and Proposition~\ref{prop:holder-regularity} give a uniform H\"older bound for $\mathcal G(u,\Phi)$. Interpolation with its $L^2$-norm yields
	\[
	\mathcal G(u,\Phi)\leqslant C\|\mathcal G(u,\Phi)\|_{L^2(B_1)}^{\kappa_0}
	\]
	on a fixed neighbourhood of $K$, for some $\kappa_0\in(0,1)$.

	Choose a local model branch $\varphi$. Its norm has a uniform positive lower bound; choose $\varepsilon$ so that $C\varepsilon^{\kappa_0}<\frac{1}{2}|\varphi|$ on this neighbourhood. The hypothesis then gives $\mathcal G(u,\Phi)<|\varphi|/2$. The disjoint fibre balls of radius $|\varphi|/2$ around $\pm\varphi$ therefore contain one value of $u$ each. This uniquely selects continuous $W^{1,2}$-branches, which satisfy the Dirac equation weakly.

	Writing $u$ for this branch, we have $\mathcal G(u,\Phi)=|u-\varphi|$ and $\mathcal D(u-\varphi)=(\mathcal D_0-\mathcal D)\varphi$. The derivatives of $\varphi$ are uniformly bounded away from $P$, while the coefficients of $\mathcal D-\mathcal D_0$ and their derivatives are $O(\delta)$. Interior elliptic estimates therefore give \eqref{eqn:sheet-matching-estimate}.
	\end{proof}
	\subsection{Topological exclusion of holes}
	\label{sec:no-hole-principle}
	
	The matching in Lemma~\ref{lem:sheet-matching} identifies the $\mathbb Z_2$-monodromy of $u$ and $\Phi$: their sheets have the same sign change around closed curves. In particular, $u$ has nontrivial monodromy around each small normal circle where the sheets are matched.

	\begin{lemma}[No-hole principle]
	\label{lem:no-hole-principle}
	For each $0<\tau<1/4$, there exists $\varepsilon(\tau)>0$, uniform for $\Phi\in\mathcal C_0$, such that
	\[
	\int_{B_1}\mathcal G(u,\Phi)^2\,dX+\delta<\varepsilon(\tau)
	\]
	implies that every normal disk of radius $\tau/4$ centred at a point of $P\cap B_{1/2}$ contains a branch point.
	\end{lemma}

	\begin{proof}
	The boundary circles of these disks stay a fixed positive distance from $P$ and $\partial B_{3/4}$, depending only on $\tau$. For sufficiently small $\varepsilon(\tau)$, Lemma~\ref{lem:sheet-matching} gives nonzero matched sheets near all these circles. Thus $u$ has the same nontrivial $\mathbb Z_2$-monodromy as $\Phi$ around each circle.

	If a disk contained no branch point, choose a contractible neighbourhood of the disk disjoint from $\mathcal B_u$, using that $\mathcal B_u$ is closed. Thus we can choose local continuous sheets of $u$. They satisfy the Dirac equation weakly and hence are smooth. By unique continuation and contractibility, $u$ admits a single-valued smooth lift on the disk. The monodromy around the boundary circle is therefore trivial, a contradiction.
	\end{proof}

	\subsection{The estimates at branch points}
	\label{sec:branch-point-estimates}
	Assume now that $0\in\mathcal B_u$, and define
	\[
	\mathcal E(u,\Phi;p,s)=s^{-n-1}\int_{B_s(p)}\mathcal G(u,\Phi)^2\,dV_g, \qquad E:=\mathcal E(u,\Phi;0,1).
	\]
	For $p=(\xi,\zeta)$, write $R_p=d_g(p,\cdot)$ and let $\partial_{R_p}$ denote radial differentiation in geodesic coordinates centred at $p$. 
	
	\begin{proposition}[$L^2$-estimates at branch points]
		\label{prop:branch-point-estimates}
		There exists $\varepsilon>0$, depending only on the fixed model family and ambient bounds, such that if $E+\delta<\varepsilon$, then for every $p=(\xi,\zeta)\in\mathcal B_u\cap B_{9/16}$,
		\begin{equation}
			|\xi|^2+\int_{B_{1/2}}R^{2-n}\left|\partial_R\left(\frac{u}{R^{1/2}}\right)\right|^2+\int_{B_{1/2}}|D_yu|^2\leqslant C(E+\delta).
			\label{eqn:branch-radial-tilt}
		\end{equation}
		For every \(0<\sigma<1\), we have the non-concentration estimate, with translation taken in the original coordinates,
		\begin{equation}
			\int_{B_{1/4}(p)}\frac{\mathcal G(u(X),\Phi(X-p))^2}{R_p(X)^{n+1-\sigma}}\,dV_g(X)+\int_{B_{1/2}}\frac{\mathcal G(u(X),\Phi(X))^2}{R_p(X)^{n-1-\sigma}}\,dV_g(X)\leqslant C_\sigma(E+\delta).
			\label{eqn:branch-weighted-excess}
		\end{equation}
	\end{proposition}
	\begin{remark}\label{rmk:branch-estimate-heuristics}
		Estimate \eqref{eqn:branch-radial-tilt} controls the normal displacement, homogeneity defect and tangential energy. The bound for the second integral in \eqref{eqn:branch-weighted-excess}, combined with the no-hole principle, yields $L^2$ non-concentration of the normalised errors near the spine.
	\end{remark}
	
	\begin{proof}
		\textbf{Step 1: the Weiss energy and variational inequalities.}
		Let $H$ and $\mathscr J$ be as in Subsection~\ref{sec:frequency-control}, with centre $0$. Define the Weiss energy
		\[
		W(s)=s^{1-n}\mathscr J(s)-\tfrac{1}{2}s^{-n}H(s),\qquad 0<s<3/4.
		\]
		By Proposition~\ref{prop:minimum-frequency-classification} and the frequency estimates, $N(s)\geqslant\frac{1}{2}-C\delta s^2$. The uniform $L^2$-bound and \eqref{eqn:caccioppoli-estimate} then imply
		\[
		H(s)\leqslant Cs^n,\qquad \int_{B_s}|\nabla u|_g^2\,dV_g\leqslant Cs^{n-1}.
		\]
		Consequently $W(s)=s^{-n}H(s)(N(s)-\frac{1}{2})\geqslant-C\delta s^2$, so $\liminf_{s\downarrow0}W(s)\geqslant0$. Using \eqref{eqn:dirac-boundary-flux}--\eqref{eqn:energy-radial-variation} for the rescaled data, we obtain
		\[
		W'(s)=2s^{1-n}\int_{\partial B_s}\left|\nabla_{\partial_R}u-\frac{u}{2s}\right|^2\,d\mathrm{Area}_g
		+O\!\left(\delta s^{2-n}\int_{B_s}|\nabla u|_g^2\,dV_g+\delta s^{1-n}H(s)\right).
		\]
		The preceding bounds give an $O(\delta s)$ error. Integrating yields
		\begin{equation}
			W(s)-W(t)=2\int_{B_s\setminus B_t}R^{2-n}\left|\partial_R\left(\frac{u}{R^{1/2}}\right)\right|^2\,dV_g
            +O(\delta s^2),\qquad 0<t<s<3/4.
			\label{eqn:weiss-variation}
		\end{equation}
		Letting $t\downarrow0$ along a sequence attaining the lower limit, we obtain
		\[
		W(s)\geqslant2\int_{B_s}R^{2-n}\left|\partial_R\left(\frac{u}{R^{1/2}}\right)\right|^2\,dV_g-C\delta s^2.
		\]

		Choose a smooth nonincreasing function $\psi:[0,\infty)\to[0,1]$ with $\psi=1$ on $[0,1/2]$ and $\psi=0$ on $[3/4,\infty)$. Throughout the following calculations, we write $\psi$ and $\psi'$ for $\psi(R)$ and $\psi'(R)$. By Remark~\ref{rmk:rescaled-geometry} and the Caccioppoli estimate \eqref{eqn:caccioppoli-estimate}, replacing $|\nabla u|_g^2+\langle\mathscr Ru,u\rangle$ and $dV_g$ by $|Du|^2$ and $dX$, respectively, introduces an $O(\delta)$ error. Combined this with the coarea formula and integration by parts, we find	
\begin{equation}
    \begin{aligned}
    &\int_{B_{3/4}}(|Du|^2\psi^2+R^{-1}|u|^2\psi\psi')\,dX\\
    ={}&\int_{B_{3/4}}[(|\nabla u|_g^2+\langle\mathscr Ru,u\rangle)\psi^2+R^{-1}|u|^2\psi\psi']\,dV_g+O(\delta)\\
    ={}&\int_{1/2}^{3/4}(-2\psi(s)\psi'(s))s^{n-1}W(s)\,ds+O(\delta)\\
    \geqslant{}&2\int_{1/2}^{3/4}(-2\psi\psi')s^{n-1}\left(\int_{B_s}R^{2-n}\left|\partial_R\left(\frac{u}{R^{1/2}}\right)\right|^2\,dV_g\right)ds-C\delta\\
    \geqslant{}&c\int_{B_{1/2}}R^{2-n}\left|\partial_R\left(\frac{u}{R^{1/2}}\right)\right|^2\,dX-C\delta.
    \end{aligned}
    \label{eqn:weiss-cutoff-identity}
\end{equation}
		where $c>0$ depends only on $n$.
		
		To control the tangential energy as well, apply the outer variation $\psi^2u$ and the inner variation $\psi^2(z,0)$ to the rescaled $\mathscr J$-energy, with the same $O(\delta)$ control of the coordinate errors. Thus

		\[
		\begin{aligned}
			\frac{1}{2}\int|Du|^2\psi^2\,dX
			&=-\int\langle u,\rho\partial_\rho u+y\cdot D_yu\rangle
			R^{-1}\psi\psi'\,dX+O(\delta),\\
			\int|D_yu|^2\psi^2\,dX
			&=-2\int\left(\tfrac{1}{2}\rho^2|Du|^2
			-\langle\rho\partial_\rho u,\rho\partial_\rho u+y\cdot D_yu\rangle\right)
			R^{-1}\psi\psi'\,dX+O(\delta).
		\end{aligned}
		\]
		To express the error relative to $\Phi$, set
		\[
		\mathrm{Err}(u,\Phi)=\tfrac{1}{2}\rho^2(|Du|^2-|D\Phi|^2)-\langle\rho\partial_\rho u,\rho\partial_\rho u-u/2\rangle-\tfrac{1}{4}(|u|^2-|\Phi|^2).
		\]
		The two terms involving only $\Phi$ cancel pointwise, since $\rho^2|D\Phi|^2=\tfrac{1}{2}|\Phi|^2$ for the Dirac models in Definition~\ref{def:half-frequency-models}. Adding the variation identities and applying Young's inequality to
		$\langle y\cdot D_yu,\rho\partial_\rho u-u/2\rangle$ gives
		\begin{equation}
			\begin{aligned}
				\label{eqn:radial-tilt-variation}
				\tfrac{1}{2}\int (|Du|^2+|D_yu|^2)\psi^2\,dX+\tfrac{1}{2}\int R^{-1}|u|^2\psi\psi'\,dX&\leqslant -2\int \mathrm{Err}(u,\Phi)R^{-1}\psi\psi'\,dX\\&\quad{}+C\int|\rho\partial_\rho u-u/2|^2(\psi')^2\,dX+C\delta.
			\end{aligned}
		\end{equation}
		
		\textbf{Step 2: Tilt-excess along the spine.}
		To estimate the right-hand side of \eqref{eqn:radial-tilt-variation}, we localise at the distance to the spine, as in the graphical description of triple junctions in \cite{Min24}.
		For $s>0$, $y_0\in\mathbb R^{n-2}$ and $0<\kappa\leqslant1/8$, set
		\[
		A_{s,\kappa}(y_0)=\{(z,y):(\rho-s)^2+|y-y_0|^2<(\kappa s)^2\},
		\]
		on which $\rho$ is comparable to $s$. We sometimes suppress $y_0$ from the notation. Choose a cover of $\{1/2<R<3/4\}\setminus P$ by $A_{s_i,1/16}(y_i)$, with the enlarged regions $A_{s_i,1/8}(y_i)$ contained in $B_1$ and having bounded overlap.
		The main strategy is a covering argument: regions with small rescaled $L^2$-excess admit matched sheets by Lemma~\ref{lem:sheet-matching}, while on the remaining regions the Caccioppoli estimate suffices.

		Choose $\varepsilon_0>0$ small enough for the sheet-matching argument in Lemma~\ref{lem:sheet-matching} to apply after rescaling.
		If $s^{-n-1}\int_{A_{s,1/8}}\mathcal G(u,\Phi)^2\,dX<\varepsilon_0$, write $v=u-\varphi$ for the difference of the locally matched sheets. It satisfies a $\delta$-perturbed harmonic equation. Here $|\varphi|^2+s^2|D\varphi|^2\leqslant Cs$ on $A_{s,1/8}$. Test with $\eta^2v$, where $\eta$ depends only on $\rho$ and $y$, equals one on $A_{s,1/16}$, is supported in $A_{s,3/32}$ and satisfies $|D\eta|\leqslant C/s$. This gives
		\begin{equation}
			\begin{aligned}
			\int_{A_{s,1/16}}(|v|^2+s^2|Dv|^2)\,dX
			&\leqslant C\int_{A_{s,1/8}}\mathcal G(u,\Phi)^2\,dX+C\delta^2\int_{A_{s,1/8}}(|\varphi|^2+s^2|D\varphi|^2)\,dX\\
			&\leqslant C\int_{A_{s,1/8}}\mathcal G(u,\Phi)^2\,dX+C\delta^2s^{n+1}.
			\end{aligned}
			\label{eqn:matched-caccioppoli}
		\end{equation}
		Writing $z=\rho e^{i\theta}$ and using $\rho\partial_\rho\varphi=\varphi/2$ and $D_y\varphi=0$, we have
		\[
		\begin{aligned}
		\mathrm{Err}(u,\Phi)
		&=\langle\partial_\theta\varphi,\partial_\theta v\rangle-\tfrac{1}{4}\langle\varphi,v\rangle+O(s^2|Dv|^2+|v|^2),\\ \rho\partial_\rho u-u/2&=\rho\partial_\rho v-v/2.
		\end{aligned}
		\]
		The two linear terms have zero angular integral:
		\begin{equation}
			\int_0^{2\pi}\left(\langle\partial_\theta\varphi,\partial_\theta v\rangle-\tfrac{1}{4}\langle\varphi,v\rangle\right)\,d\theta=[\langle\partial_\theta\varphi,v\rangle]_0^{2\pi}-\int_0^{2\pi}\left\langle\partial_\theta^2\varphi+\tfrac{1}{4}\varphi,v\right\rangle\,d\theta=0,
			\label{eqn:angular-cancellation}
		\end{equation}
		since $\varphi$ and $v$ are $2\pi$-antiperiodic and $\partial_\theta^2\varphi=-\varphi/4$. The cutoff weights are independent of $\theta$, so \eqref{eqn:matched-caccioppoli} controls both error terms after first integrating in $\theta$ and then taking absolute values.

		On the region not covered by such annuli, $s^{n+1}\leqslant\varepsilon_0^{-1}\int_{A_{s,1/8}}\mathcal G(u,\Phi)^2\,dX$. The Caccioppoli estimate \eqref{eqn:caccioppoli-estimate}, with a cutoff in $A_{s,1/8}$, and $|u|^2\leqslant C(\mathcal G(u,\Phi)^2+|\Phi|^2)$ give
		\begin{equation}
			\int_{A_{s,1/16}}(|u|^2+|\Phi|^2+s^2|Du|^2+s^2|D\Phi|^2)\,dX\leqslant C\int_{A_{s,1/8}}\mathcal G(u,\Phi)^2\,dX.
			\label{eqn:unmatched-caccioppoli}
		\end{equation}
		This controls the same error terms directly from the definition of $\mathrm{Err}(u,\Phi)$.

		The weights $R^{-1}\psi\psi'$ and $(\psi')^2$ are bounded. Since the enlarged regions have bounded overlap, the sum of their local excesses is at most $CE$. Their volumes are comparable to $s_i^n$, so
		\[
		\delta^2\sum_i s_i^{n+1}
		\leqslant C\delta^2\sum_i\int_{A_{s_i,1/8}(y_i)}\rho\,dX
		\leqslant C\delta^2\int_{B_1}\rho\,dX\leqslant C\delta^2.
		\]
		Summing the nonnegative integrals obtained above therefore bounds the right-hand side of \eqref{eqn:radial-tilt-variation} by $C(E+\delta)$.
		Combining this with \eqref{eqn:weiss-cutoff-identity} and comparability of $dX$ and $dV_g$, we obtain
		\begin{equation}
		    \int_{B_{1/2}}R^{2-n}\left|\partial_R\left(\frac{u}{R^{1/2}}\right)\right|^2\,dV_g
		    +\int_{B_{1/2}}|D_yu|^2\,dV_g\leqslant C(E+\delta).
		    \label{eqn:centred-radial-tilt}
		\end{equation}

		\textbf{Step 3: Changing the centre.}
		Fix $p=(\xi,\zeta)\in\mathcal B_u\cap B_{9/16}$. Translation invariance along $P$ and $|D\Phi|^2\leqslant C/|z|$ give
		\begin{equation}
		    \mathcal G(\Phi(X-p),\Phi(X))^2
		    \leqslant C|\xi|^2\int_0^1|z-t\xi|^{-1}\,dt.
		    \label{eqn:model-translation-bound}
		\end{equation}
		Thus $\int_{B_1}\mathcal G(u(X),\Phi(X-p))^2\,dV_g(X)\leqslant C(E+|\xi|^2)$, with $|\xi|\to0$ as $E+\delta\to0$ by Lemma~\ref{lem:sheet-matching}.

		Choose geodesic coordinates $Y$ centred at $p$ whose coordinate basis at $p$ is obtained by orthonormalising the original coordinate basis. Extend the original bundle frame at $p$ by parallel transport along radial geodesics. The coordinate change has the form
		\[
		X=p+(I+O(\delta))Y+O(\delta R_p^2),
		\]
		and the bundle frame change is $I+O(\delta R_p)$.
		Write $\Phi=\{\pm(az^{1/2}+b\bar z^{1/2})\}$. Parallel transport the oriented normal frame and the coefficients $a,b$ from $0$ to $p$ along the radial geodesic. This defines a cylindrical Dirac model $\Psi$ at $p$, since parallel transport preserves the Clifford relations.

		Comparing in the same bundle frame, for corresponding coordinates $X,Y$ with $X\notin p+P$, we have
		\[
		\mathcal G(\Phi(X-p),\Psi(Y))^2
		\leqslant C\delta^2\left(\frac{R_p^2}{\operatorname{dist}(X,p+P)}+R_p^3\right)
		\leqslant C\delta^2\frac{R_p^2}{\operatorname{dist}(X,p+P)}.
		\]
		The $O(\delta^2R_p^2/\operatorname{dist}(X,p+P))$ error comes from the $O(\delta R_p)$ change in the variable $z$ in the model formula \eqref{eqn:half-frequency-model}; the $O(\delta^2R_p^3)$ error comes from the bundle frame change. Integrating in polar coordinates for $X-p$ gives
		\begin{equation}
		    \int_{B_{3/8}}\frac{\mathcal G(\Phi(X-p),\Psi(Y))^2}{R_p(X)^{n+1-\sigma}}\,dY
		    \leqslant C\delta^2\int_0^1 t^{\sigma-1}\,dt
		    \leqslant C_\sigma\delta^2.
		    \label{eqn:recentred-model-error}
		\end{equation}

		It follows that $\int_{B_{3/8}}\mathcal G(u,\Psi)^2\,dY\leqslant C(E+|\xi|^2+\delta^2)$. Applying Steps 1 and 2 with a cutoff supported in $B_{3/8}$ and equal to one on $B_{5/16}$ gives
		\begin{equation}
		    \int_{B_{5/16}}R_p^{2-n}
		    \left|\partial_{R_p}\left(\frac{u}{R_p^{1/2}}\right)\right|^2\,dY
		    \leqslant C(E+|\xi|^2+\delta).
		    \label{eqn:recentred-radial-bound}
		\end{equation}
		Choose a radial $\eta\in C_c^\infty(B_{5/16})$ with $0\leqslant\eta\leqslant1$, $\eta=1$ on $B_{9/32}$, and $|D\eta|\leqslant C$. Integration by parts, with nonpositive inner boundary term, Young's inequality and the homogeneity of $\Psi$ give
		\[
		\begin{aligned}
		\int_{B_{9/32}}\frac{\mathcal G(u,\Psi)^2}{R_p^{n+1-\sigma}}\,dY
		&\leqslant C_\sigma\int_{B_{5/16}}R_p^{2-n+\sigma}
		\left|\partial_{R_p}\left(\frac{\eta\mathcal G(u,\Psi)}{R_p^{1/2}}\right)\right|^2\,dY\\
		&\leqslant C_\sigma\int_{B_{5/16}}\left(
		R_p^{2-n}\left|\partial_{R_p}\left(\frac{u}{R_p^{1/2}}\right)\right|^2
		+\mathcal G(u,\Psi)^2\right)\,dY\\
		&\leqslant C_\sigma(E+|\xi|^2+\delta).
		\end{aligned}
		\]
		Since $B_{1/4}(p)$ lies inside $B_{9/32}$ in the new coordinates, \eqref{eqn:recentred-model-error} and volume comparability give, in the original chart,
		\begin{equation}
		    \int_{B_{1/4}(p)}
		    \frac{\mathcal G(u(X),\Phi(X-p))^2}{R_p(X)^{n+1-\sigma}}\,dV_g(X)
		    \leqslant C_\sigma(E+|\xi|^2+\delta).
		    \label{eqn:weighted-predisplacement-bound}
		\end{equation}

		\textbf{Step 4: Estimating the displacement $|\xi|=\operatorname{dist}(p,\mathbb R^{n-2})$.}
		By the explicit formula \eqref{eqn:half-frequency-model} for $\Phi$ and uniform comparability of geodesic and coordinate balls, if $|\xi|\leqslant s\ll\frac{1}{8}$,
		\begin{equation}
			\int_{B_{s}(p)}\mathcal G(\Phi(X),\Phi(X-p))^2\,dV_g(X)
			\gtrsim  \int_{B_s(p)}\frac{|\xi|^2}{s}\,dV_g(X)\gtrsim s^{n-1}|\xi|^2.
			\label{eqn:model-shift-lower}
		\end{equation}
		Taking \(\sigma=1/2\) in \eqref{eqn:weighted-predisplacement-bound}, we have
		\[
		\begin{aligned}
		cs^{n-1}|\xi|^2
		&\lesssim \int_{B_s(p)}(\mathcal{G}(u(X),\Phi(X))^2+\mathcal{G}(u(X),\Phi(X-p))^2)\,dV_g(X)\\
		&\lesssim E+s^{n+1/2}(E+|\xi|^2+\delta).
		\end{aligned}
		\]
		Thus $|\xi|^2\leqslant C(E+\delta)$. Together with \eqref{eqn:centred-radial-tilt} and \eqref{eqn:weighted-predisplacement-bound}, this proves \eqref{eqn:branch-radial-tilt} and bounds the first integral in \eqref{eqn:branch-weighted-excess}.
		
		For the second integral in \eqref{eqn:branch-weighted-excess}, the triangle inequality reduces the remaining estimate to
		\(\int_{B_{1/2}}\frac{\mathcal G(\Phi(X-p),\Phi(X))^2}{R_p(X)^{n-1-\sigma}}\,dV_g(X)\). By \eqref{eqn:model-translation-bound} and comparability of $dV_g$ and $dX$, it suffices to show
		$$\int_{B_{1/2}}\frac{dX}{|z-t\xi|R_p^{n-1-\sigma}}\leqslant C_\sigma,\quad \forall t\in [0,1].$$
		In fact, when $n\geqslant 3$, using $R_p\gtrsim |z-\xi|,|y-\zeta|$,
		$$\int_{B_{1/2}}\frac{dX}{|z-t\xi|R_p^{n-1-\sigma}}\lesssim \left(\int_{|z|\lesssim 1}\frac{dz}{|z-t\xi||z-\xi|^{1-\sigma/2}}\right)\left(\int_{|y-\zeta|\lesssim1}\frac{dy}{|y-\zeta|^{n-2-\sigma/2}}\right)<\infty.$$
		When $n=2$, it directly follows from $|z|^{-2+\sigma}\in L_{\operatorname{loc}}^1$.
	\end{proof}

	\begin{corollary}[$L^2$ non-concentration near the spine]
		\label{cor:spine-nonconcentration}
		For each sufficiently small \(\tau>0\), if the height excess satisfies \(E+\delta<\varepsilon(\tau)\), we have
		\begin{equation}
			\int_{B_{1/2}} \frac{\mathcal G(u,\Phi)^2}{\max\{\rho,\tau\}^{1-\sigma}} \leqslant C_\sigma(E+\delta) .
			\label{eqn:spine-nonconcentration}
		\end{equation}
		The constant $C_\sigma$ is independent of \(\tau\).
	\end{corollary}
	
	\begin{proof}
		At the prescribed scale \(\tau\), by Lemma~\ref{lem:no-hole-principle} we know that for every \(y\in P\cap B_{1/2}\), the normal slice contains a branch point \(p_y\in \mathcal{B}_u\) with \(|p_y-y|<\tau/4\). For \(\tau\leqslant s\leqslant \frac{1}{4}\), applying
		\eqref{eqn:branch-weighted-excess} with exponent \(\sigma/2\), we find
		\[
		\int_{B_s(y)\cap B_{1/2}}\mathcal G(u,\Phi)^2 \leqslant C_\sigma s^{n-1-\sigma/2}(E+\delta).
		\]
		Covering the \(s/2\)-tube of the spine by at most \(Cs^{-(n-2)}\) such balls, we get \(M(s)\leqslant C_\sigma s^{1-\sigma/2}(E+\delta)\) for all $s\geqslant \tau$, where \(M(t):=\int_{B_{1/2}\cap\{\rho<t\}}\mathcal G(u,\Phi)^2\).
		
		Thus,
		\[
		\begin{split}
			\int_{B_{1/2}}\frac{\mathcal G(u,\Phi)^2}{\max\{\rho,\tau\}^{1-\sigma}} \lesssim \underbrace{C_\sigma \tau^{-1+\sigma}M(\tau)}_{\rho<\tau} +\underbrace{C_\sigma\int_\tau^{1/4}s^{-2+\sigma}M(s)\,ds}_{\tau<\rho<\frac{1}{4}}+\underbrace{CE}_{\rho \sim \mathcal{O}(1)}\lesssim C_\sigma(E+\delta).
		\end{split}
		\]
		This completes the proof. 
	\end{proof}
	
	\subsection{Harmonic blow-ups}
	\label{sec:harmonic-blowups}
	Let $(u_j,\Phi_j)$ be as above, with their model spines identified with $P$. Assume $\Phi_j\to\Phi_0=\{\pm\varphi\}$, $E_j:=\mathcal E(u_j,\Phi_j;0,1)\to0$ and $\delta_j/E_j\to0$. We study the blow-up limit of the renormalised differences
	\[
	v_j=\frac{u_j-\varphi_j}{\sqrt{E_j}} \quad\text{on }B_{1/2}\cap\{\rho>\tau_j\}.
	\]
	Here $\tau_j\downarrow0$ is chosen so that the sheets of $u_j$ and $\Phi_j$ are uniquely matched on this region by Lemma~\ref{lem:sheet-matching}. For $\lambda=\lambda_1+i\lambda_2$, write $\lambda\cdot D_z=\lambda_1\partial_{x_1}+\lambda_2\partial_{x_2}$.
	
	\begin{proposition}
		\label{prop:limiting-error-estimates}
		After taking a subsequence, $v_j\to w$ strongly in $L^2(B_{1/2})$. It is smooth and satisfies $\mathcal D_0w=0$ on $B_{1/2}\setminus P$. Moreover,
		\begin{equation}
		\int_{B_{1/4}}R^{2-n}|\partial_R(w/R^{1/2})|^2\leqslant C.
		\label{eqn:blowup-radial-bound}
		\end{equation}
		For each $|y_0|<1/8$, there is $\lambda(y_0)\in\mathbb C$ with $|\lambda(y_0)|\leqslant C$ such that for every $0<\sigma<1$,
		\begin{equation}
		\int_{B_{1/8}(0,y_0)} \frac{|w-\lambda(y_0)\cdot D_z\varphi|^2}{|X-(0,y_0)|^{n+1-\sigma}} \leqslant C_\sigma.
		\label{eqn:blowup-weighted-bound}
		\end{equation}
	\end{proposition}
	
	\begin{proof}
		We first prove the strong $L^2$-convergence and then derive the limiting equation for $w$. Fix $0<\tau<1/8$ and $0<\sigma<1$. Summing \eqref{eqn:matched-caccioppoli} away from $P$ and applying Corollary~\ref{cor:spine-nonconcentration} give, for all sufficiently large $j$,
		\[
		\begin{aligned}
		\int_{B_{1/2}\cap\{\rho>\tau\}}|Dv_j|^2
		&\leqslant C\tau^{-2}(1+\delta_j/E_j),\\
		\frac{1}{E_j}\int_{B_{1/2}\cap\{\rho<\tau\}}\mathcal G(u_j,\Phi_j)^2
		&\leqslant C_\sigma\tau^{1-\sigma}(1+\delta_j/E_j).
		\end{aligned}
		\]
		Rellich compactness away from $P$ and the second estimate give the claimed convergence by first letting $j\to\infty$ and then $\tau\downarrow0$, with $\int_{B_{1/2}}|w|^2\leqslant1$. The source terms in the Dirac equations for $v_j$ are bounded on compact subsets away from $P$ by $C\delta_j/\sqrt{E_j}\to0$. Passing to the limit and applying interior elliptic regularity gives a smooth solution of $\mathcal D_0w=0$ on $B_{1/2}\setminus P$.
		
		Now we pass the $L^2$-estimates (\cf Proposition \ref{prop:branch-point-estimates}) to establish the one for $w$. The $1/2$-homogeneity of $\varphi_j$ gives $\partial_R(v_j/R^{1/2})=E_j^{-1/2}\partial_R(u_j/R^{1/2})$. Thus \eqref{eqn:blowup-radial-bound} follows from \eqref{eqn:branch-radial-tilt} by weak lower semicontinuity away from $P$, followed by an exhaustion argument.

		Fix $|y_0|<1/8$. The no-hole principle gives branch points $p_j=(\xi_j,y_0)$ approaching $(0,y_0)$, and \eqref{eqn:branch-radial-tilt} gives $|\xi_j|^2\leqslant C(E_j+\delta_j)$. For this fixed $y_0$, pass to a further subsequence so that
		\[
		-\frac{\xi_j}{\sqrt{E_j}}\longrightarrow\lambda(y_0),\qquad
		\frac{u_j-\varphi_j(\,\cdot-p_j)}{\sqrt{E_j}}
		\longrightarrow w-\lambda(y_0)\cdot D_z\varphi
		\]
		locally away from $P$, with $|\lambda(y_0)|\leqslant C$. Divide the bound for the first integral in \eqref{eqn:branch-weighted-excess} by $E_j$ and apply Fatou's lemma on $B_{1/8}(0,y_0)$ to obtain \eqref{eqn:blowup-weighted-bound}.
	\end{proof}
	
	\section{Excess decay and \texorpdfstring{$C^{1,\alpha}$}{C1-alpha}-regularity}
	\label{sec:excess-branch-regularity}
	
	In this section, we prove the excess decay via Campanato iteration and finally achieve the $C^{1,\alpha}$ regularity of the frequency $\frac{1}{2}$ branching set for the given $\mathbb Z_2$-harmonic spinor.

	\subsection{The leading term of the linear error}
	\label{sec:leading-linear-term}

	\begin{lemma}[Excess decay of the blow-up]
		\label{lem:linear-excess-decay}
		Let $w$ be a limit in Proposition~\ref{prop:limiting-error-estimates}, with comparison models converging to $\Phi_0$. Then there exists
		$\ell\in T_{\Phi_0}\mathcal C_{\mathcal D}$ such that
		\begin{equation}
			\|\ell\|_{L^2(B_{1/2})}^2 +s^{-n-3}\int_{B_s}|w-\ell|^2 \leqslant C\int_{B_{1/2}}|w|^2, \qquad \forall 0<s<\frac{1}{32}.
			\label{eqn:linear-excess-decay}
		\end{equation}
	\end{lemma}
	\begin{proof}
		We fix $\sigma=1/2$ in Proposition \ref{prop:limiting-error-estimates}. In polar coordinates $(\rho,\theta)$ in the normal direction, $w$ is $4\pi$-periodic in $\theta$, and we decompose it as
		$$w=\sum_{\nu}w_\nu=\sum_\nu(a_\nu(\rho,y) \cos(\nu\theta)+b_\nu(\rho,y)\sin(\nu\theta)),\quad \nu=\frac{1}{2},\frac{3}{2},\cdots.$$
		Each coefficient $f_\nu=a_\nu$ or $b_\nu$ solves \begin{equation}
			\partial_\rho^2f_\nu+\rho^{-1}\partial_\rho f_\nu
			-\nu^2\rho^{-2}f_\nu+\Delta_yf_\nu=0.
			\label{eqn:angular-mode-equation}
		\end{equation}
		
		\smallskip
		\textbf{Case 1: The frequency $\nu=1/2$.}
		Write $g(\rho,y):=\rho^{1/2}w_{1/2}$, which takes values in the real vector space $\mathscr V=\{a\cos(\theta/2)+b\sin(\theta/2):a,b\in E_0\}$ with the real $L^2(0,2\pi)$-inner product. Equation \eqref{eqn:angular-mode-equation} gives $(\partial_\rho^2+\Delta_y)g=0$.
		Set $t_i(\theta)=\rho^{1/2}\partial_i\varphi$ for $i=1,2$; these are independent of $\rho$. Decompose $g=g_\parallel+g_\perp$ relative to $\mathscr T=\operatorname{span}_{\mathbb R}\{t_1,t_2\}$ and its orthogonal complement in $\mathscr V$. We write $B_r^+$ for the half-ball in the variables $(\rho,y)$.
		
		Since $\rho^{1/2}\lambda(y_0)\cdot D_z\varphi\in\mathscr T$, the weighted estimate \eqref{eqn:blowup-weighted-bound} gives, for $|y_0|<1/8$ and $0<s<1/16$,
		\[
		\int_{B_s^+(0,y_0)}\|g_\perp\|_{\mathscr V}^2\,d\rho\,dy
		\leqslant C s^{n+1-\sigma}.
		\]
		The mean-value inequality in a ball of radius $\rho/2$, contained in $B_{2\rho}^+(0,y)$, now implies $|g_\perp(\rho,y)|\leqslant C\rho^{1-\sigma/2}$ for $0<\rho<1/32$ and $|y|<1/8$.
		Thus $g_\perp$ has zero continuous boundary values and its odd reflection is harmonic. Interior estimates for this reflection give
		\[
		\|g_\perp\|_{C^k(B_{1/16}^+)}
		\leqslant C_k\|g_\perp\|_{L^2(B_{1/8}^+)}
		\leqslant C_k\|w\|_{L^2(B_{1/2})},\qquad k\geqslant0.
		\]
		
		To control $g_\parallel$, use the Dirac equation. With $J=\mathfrak c(e_1)\mathfrak c(e_2)$, its polar form is
		\[
		\partial_\rho w=\rho^{-1}J\partial_\theta w
		+\sum_{\alpha=1}^{n-2}\mathfrak c(\partial_\rho)\mathfrak c(e_{\alpha+2})\partial_{y_\alpha}w.
		\]
		The model Dirac equation gives $t_2=-Jt_1$ and $J\partial_\theta t_i=-t_i/2$.
		The operator $\mathfrak c(\partial_\rho)\mathfrak c(e_{\alpha+2})$ is skew-adjoint and anticommutes with $J$, so $\langle \mathfrak c(\partial_\rho)\mathfrak c(e_{\alpha+2})t_j,t_i\rangle=0$ for $i,j=1,2$.
		Moreover, $\mathfrak c(\partial_\rho)\mathfrak c(e_{\alpha+2})t_i$ has only angular frequencies $\pm1/2$, by the explicit form \eqref{eqn:half-frequency-model} of $\Phi$.
		Since $J^*=-J$ and both $w$ and $t_i$ change sign after one turn, integration by parts gives
		\[
		\rho^{-1/2}\int_0^{2\pi}\langle J\partial_\theta w,t_i\rangle\,d\theta
		=\rho^{-1/2}\int_0^{2\pi}\langle w,J\partial_\theta t_i\rangle\,d\theta
		=-\frac{1}{2\rho}\langle g,t_i\rangle_{\mathscr V}.
		\]
		This cancels the term $\frac{1}{2\rho}\langle g,t_i\rangle_{\mathscr V}$ arising when we differentiate the factor $\rho^{1/2}$ in $g$. Hence
		\[
		\partial_\rho\langle g,t_i\rangle_{\mathscr V}=\sum_{\alpha=1}^{n-2}\partial_{y_\alpha}
		\int_0^{2\pi}\langle \mathfrak c(\partial_\rho)\mathfrak c(e_{\alpha+2})g_\perp,t_i\rangle\,d\theta.
		\]
		The right-hand side extends smoothly to $\rho=0$ and vanishes there, since $g_\perp$ has zero boundary values. As $t_1,t_2$ form a basis of $\mathscr T$, integration from a fixed positive value of $\rho$ shows that $g_\parallel$ extends smoothly with $\partial_\rho g_\parallel(0,y)=0$. Its even reflection is therefore harmonic. Interior estimates give
		\[
		\|g\|_{C^2(B_{1/32}^+)}\leqslant C\|w\|_{L^2(B_{1/2})},
		\qquad g(0,y)=g_\parallel(0,y)\in\mathscr T.
		\]
		
		The radial estimate \eqref{eqn:blowup-radial-bound} forces $g(0,0)=0$, since a nonzero constant term in $g$ would make the integral diverge. Thus, near $(0,0)$, the Taylor expansion of $g$ takes the form
		\[
		g(\rho,y)=\rho B(\theta)+\sum_{\alpha,i}c_{\alpha i}y_\alpha (\rho^{1/2}\partial_{i}\varphi)
		+O((\rho^2+|y|^2)\|w\|_{L^2}),\quad B(\theta)\in \mathscr V.
		\]
		Then we define
		$\ell:=\rho^{1/2}B(\theta)+\sum_{\alpha,i}c_{\alpha i}y_\alpha\partial_i\varphi$. For $0<s<1/32$, the remainder estimate is given by
		$$\int_{B_s}|w_{1/2}-\ell|^2\lesssim\|w\|_{L^2(B_{1/2})}^2 \int_{\rho,|y|\lesssim s} (\rho^2+|y|^2)^2\rho^{-1}\rho d\rho dy\lesssim s^{n+3}\|w\|_{L^2(B_{1/2})}^2.$$
		
		\smallskip
		\textbf{Case 2: The higher angular frequencies $\nu \geqslant3/2$.}
		For $\nu\geqslant3/2$ set $g_\nu=\rho^{-\nu}f_\nu$. Equation
		\eqref{eqn:angular-mode-equation} becomes
		\[
		\partial_\rho^2g_\nu+(2\nu+1)\rho^{-1}\partial_\rho g_\nu +\Delta_yg_\nu=0.
		\]
		Regarding $g_\nu$ as the function $g_\nu(|x|,y)$ of $x\in\mathbb R^{2\nu+2}$, we see that $g_\nu$ is harmonic away from the codimension $2\nu+2\geqslant 5$ axis $\{x=0\}$ with $L^2$ integral $\int|f_\nu|^2\rho\,d\rho\,dy$, up to constants. A standard approximation argument implies that $g_\nu$ is weakly harmonic across the axis, and hence harmonic by Weyl's lemma. The monotonicity formula for harmonic functions implies that
		\[
		\int_{B_s}|w_\nu|^2\leqslant (4s)^{n+2\nu}\int_{B_{1/4}}|w_\nu|^2, \qquad 0<s<1/4.
		\]
		
		Since $4s<1$ and $2\nu\geqslant3$, Parseval's identity implies
		\[
		\sum_{\nu\geqslant3/2}\int_{B_s}|w_\nu|^2\leqslant(4s)^{n+3}\sum_{\nu\geqslant3/2}\int_{B_{1/4}}|w_\nu|^2\leqslant C s^{n+3}\|w\|_{L^2(B_{1/2})}^2.
		\]
		The coefficients of $\ell$ are bounded by the preceding $C^1$ 
		estimate for $g$, so $\|\ell\|_{L^2(B_{1/2})}\leqslant C\|w\|_{L^2(B_{1/2})}$ as well.

		Finally $s^{-1/2}w(s\,\cdot)\to\ell$ in $L^2$ on fixed balls as $s\downarrow0$. 
        Passing to the limit in the Dirac equation away from $P$ gives $\mathcal{D}_0\ell=0$ and hence $\ell\in T_{\Phi_0}\mathcal C_{\mathcal D}$ by Lemma \ref{lem:dirac-model-variations}. This completes the proof. 
	\end{proof}
	
	\subsection{Campanato iteration}
	\label{sec:campanato-iteration}
	
	Given Proposition~\ref{prop:limiting-error-estimates} and Lemma~\ref{lem:linear-excess-decay}, we are ready to run the Campanato iteration scheme. Fix any H\"older exponent $\beta \in (0,1)$.
	
	\begin{proposition}
		\label{prop:one-step-decay}
		There are $\theta\in(0,1/8)$, $\varepsilon>0$ and $C<\infty$ with the following property. Suppose that $p\in\mathcal B_u$ with $N_u(p,0)\geqslant\frac{1}{2}$, and a frequency-$\frac{1}{2}$ comparison Dirac model $\Phi$ in a compact non-degenerate model chart, satisfy the smallness of $L^2$ height excess, \ie
		$$\mathcal{E}(u,\Phi;p,r)+r^2\leqslant \varepsilon.$$
        Then there exists a nearby Dirac model $\Phi'$ such that
		\begin{align}
			\mathcal{E}(u,\Phi';p,\theta r)&\leqslant \tfrac{1}{2}\theta^{2\beta}\mathcal{E}(u,\Phi;p, r)+Cr^2,
			\label{eqn:one-step-decay}\\
			d_{\mathcal C}(\Phi',\Phi)&\leqslant C(\mathcal{E}(u,\Phi;p, r)^{1/2}+r).
			\label{eqn:model-parameter-increment}
		\end{align}
	\end{proposition}
	\begin{proof}
		Choose $\theta\in(0,1/32)$ such that $C\theta^2\leqslant \theta^{2\beta}/8$ for the constant $C$ in Lemma~\ref{lem:linear-excess-decay}. The estimates in Proposition~\ref{prop:branch-point-estimates} require only $N_u(p,0)\geqslant\frac{1}{2}$; the comparison models remain homogeneous of degree $\frac{1}{2}$.
		
		Suppose the conclusion fails. For each $j\geqslant1$, choose
		$(u_j,\Phi_j,p_j,r_j)$ with $E_j+r_j^2\leqslant 1/j$, where $E_j:=\mathcal{E}(u_j,\Phi_j;p_j,r_j)$, such that for every nearby Dirac model $\Phi'_j$, at least one of \eqref{eqn:one-step-decay} and \eqref{eqn:model-parameter-increment} fails with $C=j$. Taking $\Phi'_j=\Phi_j$, we know
		$$\mathcal{E}(u_j,\Phi_j;p_j,\theta r_j)>\frac{1}{2}\theta^{2\beta}E_j+jr_j^2.$$
		On the other hand, comparing the volume, we know
		$\mathcal{E}(u_j,\Phi_j;p_j,\theta r_j)\leqslant \theta^{-n-1}E_j.$
		It follows that $\frac{r_j^2}{E_j}\to 0$. Rescale to the unit ball and pass to a subsequence with $\Phi_j\to\Phi_0$. Proposition~\ref{prop:limiting-error-estimates} and Lemma~\ref{lem:linear-excess-decay} give $\ell\in T_{\Phi_0}\mathcal C_{\mathcal D}$ with excess decay estimates.
		Choose $\ell_j\in T_{\Phi_j}\mathcal C_{\mathcal D}$ converging to $\ell$. Lemma~\ref{lem:dirac-model-variations}, with parameter $\sqrt{E_j}$, gives nearby Dirac models $\Phi_j'$ with $d_{\mathcal C}(\Phi_j',\Phi_j)\leqslant C\sqrt{E_j}$. The expansion \eqref{eqn:model-variation-expansion} and strong $L^2$-convergence, together with Corollary~\ref{cor:spine-nonconcentration}, imply
		\[
		\limsup_{j\to\infty}\frac{\mathcal E(u_j,\Phi_j';0,\theta)}{E_j}\leqslant\theta^{-n-1}\int_{B_\theta}|w-\ell|^2 \leqslant C\theta^2.
		\]
		By the choice of $\theta$, both required estimates hold for large $j$, a contradiction.
	\end{proof}
	
	At a fixed centre $p\in\mathcal B_u$ with $N(p,0)=\frac{1}{2}$, Theorem~\ref{thm:spinor-blowups} gives a nonzero tangent model $\Phi$ and arbitrarily small scales $r_0$ at which the height-normalised spinor is close to $\Phi$. To start the iteration with the rescaling used here, multiply $u$ once by $(r_0/h(r_0))^{1/2}$, where $h$ is defined in Subsection~\ref{sec:frequency-control}, and retain the name $u$. Its $r_0^{-1/2}$ rescaling is precisely the height-normalised spinor, so $\mathcal E(u,\Phi;p,r_0)+r_0^2$ is sufficiently small. This multiplicative constant remains fixed throughout the iteration; it can be undone at the end. Set $r_k=\theta^k r_0$. Applying Proposition~\ref{prop:one-step-decay} iteratively, we obtain nearby flat Dirac models $\Phi_k$ with excesses $E_k:=\mathcal E(u,\Phi_k;p,r_k)$ satisfying
	\begin{equation}
		E_k\leqslant C(E_0+r_0^2)\theta^{2\beta k},\qquad\sum_{j\geqslant k}d_{\mathcal C}(\Phi_{j+1},\Phi_j) \leqslant C(E_0^{1/2}+r_0)\theta^{\beta k},\quad \forall k\geqslant 1.
		\label{eqn:iteration-model-drift}
	\end{equation}
	Thus $\{\Phi_j\}$ is Cauchy,  converging to a nondegenerate Dirac model $\Phi_p$ with $L^2$-excess estimate at all small scales, \begin{equation}
		\mathcal E(u,\Phi_p;p,s)\leqslant C(E_0+r_0^2)(s/r_0)^{2\beta},
		\qquad 0<s\leqslant r_0.
		\label{eqn:fixed-centre-decay}
	\end{equation}
	In particular, the tangent map $\Phi_p$ is unique. 
	
	The following lemma compares the parameters of nearby models with their $L^2$-distance.
	
	\begin{lemma}
		\label{lem:nearby-model-comparison}
		Let $p,p'\in\mathcal B_u$ with $d_g(p,p')<r/2$, $r\ll1$, have frequency at least $1/2$, and let $\Phi,\Phi'$ be flat Dirac comparison models centred at these points respectively, in a fixed nondegenerate model chart. Then
		\begin{equation}
			d_{\mathcal C}(\Phi,\Phi')^2 +r^{-2}|\Pi_{P^\perp}\exp_p^{-1}(p')|^2  \leqslant C( \mathcal E(u,\Phi;p,r) +\mathcal E(u,\Phi';p',r)+r^2),
			\label{eqn:nearby-model-comparison}
		\end{equation}
		provided that the right-hand side is sufficiently small, where we compare $\Phi$ and $\Phi'$ by parallel transport, and $P,P'$ stand for their spines respectively. 
	\end{lemma}
	
	\begin{proof}
		Let $\Psi$ be the model at $p'$ obtained by parallel transport of the normal frame and coefficients of $\Phi$. Proposition~\ref{prop:branch-point-estimates} bounds the normal displacement $|\Pi_{P^\perp}\exp_p^{-1}(p')|$. Together with the translation and model comparisons in Step~3 of its proof, this gives
		\[
		r^{-2}|\Pi_{P^\perp}\exp_p^{-1}(p')|^2 +r^{-n-1}\int_{B_{r/4}(p')}\mathcal G(\Phi,\Psi)^2\,dV_g \leqslant C(\mathcal E(u,\Phi;p,r)+r^2).
		\]
		Since $B_{r/4}(p')\subset B_r(p)$, the triangle inequality gives $\mathcal E(u,\Psi;p',r/4)\leqslant C(\mathcal E(u,\Phi;p,r)+r^2)$.

		In $d_{\mathcal C}(\Phi,\Phi')$, the parameters of $\Phi$ are compared at $p'$, where they are represented by $\Psi$. For models with the same centre, the explicit formula \eqref{eqn:half-frequency-model} and the fixed model gauge control parameter distance by normalised $L^2$-distance. Applying the triangle inequality through $u$ gives
		\[
		\begin{aligned}
		d_{\mathcal C}(\Phi,\Phi')^2
		&\leqslant C(\mathcal E(u,\Psi;p',r/4)+\mathcal E(u,\Phi';p',r/4))\\
		&\leqslant C(\mathcal E(u,\Phi;p,r)+\mathcal E(u,\Phi';p',r)+r^2).
		\end{aligned}
		\]
		Together these prove \eqref{eqn:nearby-model-comparison}.
	\end{proof}
	
	\subsection{Nearby centres and the branch graph}
	\label{sec:nearby-branch-graph}

	Fix a branch point $p_0$ with frequency $\frac{1}{2}$; we aim to transfer the estimates
	$\mathcal E(u,\Phi_{p_0};p_0,s)\leqslant Cs^{2\beta},s\ll 1$ to nearby
	branch points. 
	
	\begin{proposition}
		\label{prop:nearby-branch-graph}
		There are a small neighbourhood $U$ of $p_0$, radius $r_0\ll 1$, and $C<\infty$ such that every $p\in\mathcal B_u\cap U$ has frequency $1/2$ and a unique
		nondegenerate Dirac model $\Phi_p$, with
		\begin{equation}
			\mathcal E(u,\Phi_p;p,s)\leqslant Cs^{2\beta},
			\qquad \forall 0<s<r_0.
			\label{eqn:uniform-excess-decay}
		\end{equation}
		Moreover, for $p\neq p'\in \mathcal{B}_u\cap U$, 
		\begin{equation}
			d_{\mathcal C}(\Phi_p,\Phi_{p'})
			+\frac{|\Pi_{P_p^\perp}\exp_p^{-1}(p')|}{d_g(p,p')}
			\leqslant C d_g(p,p')^\beta.
			\label{eqn:two-point-comparison}
		\end{equation}
		In particular, $\mathcal B_u\cap U$ is a single $C^{1,\beta}$-graph over the spine $P_{p_0}$ with $T_p\mathcal{B}_u=P_p$. 
	\end{proposition}
	\begin{proof}
		Choosing $r_0$ sufficiently small, we can make $E_0+r_0^2$ arbitrarily small by the fixed-centre Campanato iteration. Now for $p\in\mathcal B_u\cap B_{r_0/2}(p_0)$, applying Proposition~\ref{prop:branch-point-estimates}, we have
		\[
		r_0^{-2}|\xi|^2\lesssim (E_0+r_0^2).
		\]
		Let $\Psi$ be obtained by parallel transport of the oriented normal frame and coefficients of $\Phi_{p_0}$ from $p_0$ to $p$, as in Step 3 of Proposition~\ref{prop:branch-point-estimates}. Then
		$$\mathcal E(u,\Psi;p,r_0/4)+(r_0/4)^2
		\lesssim (E_0+r_0^{-2}|\xi|^2+r_0^2)
		\lesssim C(E_0+r_0^2).$$
		Since $N_u(p,0)\geqslant1/2$, Proposition~\ref{prop:one-step-decay} applies. Then Campanato iteration implies a unique tangent Dirac model $\Phi_p$ with excess decay estimates
		\[
		\mathcal E(u,\Phi_p;p,s)
		\leqslant C(E_0+r_0^2)(s/r_0)^{2\beta},\qquad 0<s<r_0/4,
		\]
		where the constant $C$ is independent of $p$. Therefore,  $N(p,0)=\frac{1}{2}$ by Theorem~\ref{thm:spinor-blowups}. To compare nearby centres $p,p'$ with distance $d=d_g(p,p')$, Lemma~\ref{lem:nearby-model-comparison} (with $r=10d$)
		implies
		\eqref{eqn:two-point-comparison}:
		\[
		d_{\mathcal C}(\Phi_p,\Phi_{p'})^2
		+d^{-2}|\Pi_{P_p^\perp}\exp_p^{-1}(p')|^2
		\lesssim (d^{2\beta}+d^{2\beta}+d^2)\lesssim d^{2\beta}.
		\]
		The $C^{1,\beta}$-graphical representation then follows from the no-hole principle and \eqref{eqn:two-point-comparison}. 
	\end{proof}
	
	\subsection{Tangent estimates and conclusion}
	\label{sec:tangent-estimates}

	For later use in Part II, we record in the following corollary certain quantitative bounds of the height and tilt excess.  
	\begin{corollary}
		\label{cor:quantitative-tangent-estimates}
		For a sufficiently small neighbourhood $U$, there exists $\alpha\in(0,\beta)$ such that for all small scales $r\ll 1$,
		\begin{align}
			&r^{-n-1}\int_{B_r(p)}\mathcal G(u,\Phi_p)^2\leqslant Cr^{2\alpha},\qquad \sup_{B_{r/2}(p)}\mathcal G(u,\Phi_p)\leqslant Cr^{1/2+\alpha}, \label{eqn:tangent-height-estimates}\\& r^{1-n}\int_{B_{r/2}(p)}
			\mathcal G(\nabla u,\nabla\Phi_p)^2\leqslant Cr^{2\alpha}. \label{eqn:tangent-gradient-estimate}
		\end{align}
		
		For every fixed $\lambda\in(0,1/4)$ and integer $k\geqslant0$,
		locally matched sheets satisfy
		\begin{equation}
			\sup_{\substack{x\in B_{r/2}(p)\\
					\operatorname{dist}(x,\mathcal B_u)>\lambda r}}
			r^{k-1/2}|\nabla^k(u-\varphi_p)(x)|
			\leqslant C_{\lambda,k}r^\alpha,\quad \forall r\ll 1.
			\label{eqn:matched-higher-derivatives}
		\end{equation}
	\end{corollary}
	\begin{proof}[Sketch of the Proof]
		We choose
		\[
		0<\alpha<\min\left\{\frac{\beta}{3},
		\frac{\beta\gamma}{n+\gamma}\right\},
		\]
		where $\gamma$ is as in Proposition~\ref{prop:holder-regularity}, which gives local $C^{\gamma/2}$ estimates for $u$ and $\Phi_p$. By interpolation with the $L^2$ excess decay estimates, \eqref{eqn:tangent-height-estimates} follows. For the gradient estimate \eqref{eqn:tangent-gradient-estimate}, the branch locus is (locally) contained in a $Cr^{1+\beta}$-tube over the spine $P_p$ in $B_r(p)$. Away from this tube at scale $s r$ (with $s\gg r^{\beta}$ and $s^{1/2}\gtrsim r^{\beta\gamma/(n+\gamma)}$), we can select locally matched sheets by Lemma \ref{lem:sheet-matching}. Since $\mathcal G(Du,D\Phi_p)\leqslant|D(u-\varphi_p)|$ for these choices of signs, the interior estimates give
		$$r^{1-n}\int_{\rho>s r}\mathcal{G}(Du,D\Phi_p)^2\lesssim s^{-2}r^{-n-1}\int_{\rho>s r}\mathcal{G}(u,\Phi_p)^2+r^2\lesssim s^{-2}r^{2\beta}+r^2.$$
		We cover the tube by at most $Cs^{2-n}$ balls of radius $Csr$ centred on $\mathcal B_u$. The frequency estimates for $u$ and the explicit formula for $\Phi_p$ give
		\begin{equation*}
			r^{1-n}\int_{\rho<s r}\mathcal{G}(Du,D\Phi_p)^2\lesssim r^{1-n}s^{2-n}(sr)^{n-1}\lesssim s.
		\end{equation*}
		Taking $s=r^{2\alpha}$, we obtain \eqref{eqn:tangent-gradient-estimate}. At points whose distance from the spine is at least $\lambda r$, the higher order estimate \eqref{eqn:matched-higher-derivatives} follows from the interior elliptic estimates. 
	\end{proof}
	Theorem~\ref{thm:main-tangent-estimates} follows directly from Proposition \ref{prop:nearby-branch-graph} and Corollary \ref{cor:quantitative-tangent-estimates}.

\part{Higher regularity of the branching locus}
\label{part:higher-regularity}

In Part~II, we prove that the branching locus is smooth (resp. real-analytic) for smooth (resp. real-analytic) ambient data. The input from Part~I is Theorem~\ref{thm:main-tangent-estimates}: near a frequency-$\frac{1}{2}$ branch point, the branching locus is a $C^{1,\alpha}$-graph, and the spinor is approximated by non-degenerate tangent models at rate $r^{1/2+\alpha}$.

The partial Legendre transform (cf. Section~\ref{sec:legendre-coordinates}) sends the branching locus to $\{\xi=0\}$ and reduces its regularity to the tangential regularity of the transformed spinor $U$. Then the field $U$ satisfies a quasilinear elliptic system with inherited variational structure, degenerating along $\{\xi=0\}$. This degeneracy is captured by Grushin-type weighted spaces. Moreover, tangential derivatives $D_y^\beta U$ satisfy a linearised equation of the schematic form $$\mathscr B_U(Z,Y)=\mathscr F(Y),$$ with sources determined by lower-order derivatives. Higher regularity therefore reduces to a linear
theory for $\mathscr B_U$ near $\{\xi=0\}$, which we develop in three steps.
\begin{itemize}
\item \emph{The frozen operator} (Section~\ref{sec:frozen-operator}). By reduction to the square model $U_0=(\xi^2,0)$, we prove coercivity of the second variation at a transformed tangent model and establish the excess decay for weak solutions of the corresponding linear equation.
\item \emph{Excess decay} (Section~\ref{sec:linearised-polynomial-approximation}). Treating $\mathscr B_U$ as a perturbation of the frozen operator yields a Campanato-type decay to a comparison polynomial, where an extra polynomial correction (using Lemma~\ref{lem:polynomial-corrections}) is required to absorb the principal part of the source term.
\item \emph{Induction on the order of derivatives}
(Section~\ref{sec:tangential-regularity}). We upgrade the $L^2$-excess decay of Section~\ref{sec:linearised-polynomial-approximation} to weighted $C^{1,\alpha}$-estimates, which control the sources of the tangentially differentiated equations. Induction on the order of differentiation then yields smoothness and real-analyticity.
\end{itemize}

\section{Partial Legendre transform and weighted spaces}
\label{sec:legendre-coordinates}
We set up the partial Legendre transform in Section~\ref{sec:partial-legendre} and the weighted spaces adapted to the degeneracy at $\{\xi=0\}$ in Section~\ref{sec:weighted-spaces}. Section~\ref{sec:legendre-initial} translates the estimates of Part~I into these coordinates, and Section~\ref{sec:transformed-equation} derives the transformed Euler--Lagrange equation.

\subsection{Partial Legendre transform}
\label{sec:partial-legendre}
Throughout Part~II, we fix a frequency-$\frac{1}{2}$ branch point $(0,y_0)$ and work in the ball $B_{r_0}(0,y_0)$ such that Theorem~\ref{thm:main-tangent-estimates} applies. For sufficiently small $r_0\ll 1$, all constants in the estimate of Part II will be uniform over this ball. 

We use coordinates $X=(z,y)\in\C\times\R^{n-2}$ in which $P_0=\{z=0\}$ is the spine at $(0,y_0)$, and the exponential gauge centred at $(0,y_0)$, which identifies the fibres of $E$ with the fibre $E_0$ at $(0,y_0)$. By Theorem~\ref{thm:main-tangent-estimates},
\[
 \mathcal B_u\cap B_{r_0}(0,y_0)=\{p(y)=(g(y),y)\},\qquad
 g\in C^{1,\alpha},\quad g(y_0)=0,\quad Dg(y_0)=0.
\]
For $p\in\mathcal B_u\cap B_{r_0}(0,y_0)$, let $\Phi_p$ be its tangent model with spine $P_p=T_p\mathcal{B}_u$. As in Remark~\ref{rmk:model-gauge}, we write $\Phi_p=\pm T_p(\Pi_{P_p^\perp}(X-p))^{1/2}$ with
$T_p(w)=a_pw+b_p\bar w$, where $\Pi_{P_p^\perp}$ is the orthogonal projection onto $P_p^\perp\cong\C$. Fix a real-linear projection $\pi:E_0\to\C$ such that $\pi T_{p(y_0)}$ is invertible with positive determinant; by continuity, so is $\pi T_p$. Let $\pi^\perp$ be the complementary orthogonal projection, and write $\xi$ and $V$ for the components of a spinor in $E_0=\operatorname{im}\pi\oplus\operatorname{im}\pi^\perp$.

By Theorem~\ref{thm:main-tangent-estimates}, the map $(z,y)\mapsto (\xi,y):=(\pi u(z,y),y)$ is invertible, and we write its inverse and the spinor in the new coordinates as
\begin{equation}
 X(\xi,y)=(\varphi(\xi,y),y),\qquad \widehat u(\xi,y):=u(X(\xi,y))=(\xi,V(\xi,y)).
 \label{eqn:legendre-coordinate-maps}
\end{equation}
We call $U:=(\varphi,V)$ the \emph{transformed spinor}. Since $\xi\mapsto-\xi$ exchanges the two sheets of the original $\mathbb Z_2$-harmonic spinor $u$, $\varphi$ is even and $V$ is odd in $\xi$. Moreover, $\varphi(0,y)=g(y)$ and $V(0,y)=0$, so the branching locus becomes $\{\xi=0\}$, and its regularity is then encoded in the tangential regularity of $\varphi$ along
$\{\xi=0\}$.

Near each spine point $(0,y)$, we compare $U$ with its transformed tangent model at $p(y)$, which is a polynomial with coefficients frozen at $y$.

\begin{example}[The transformed tangent model]
\label{ex:transformed-model}
Fix $p=p(y)$, write $X=(z,y')$, and apply the same construction to $\Phi_p$, \ie set $\xi=\pi\Phi_p$, so that $\Pi_{P_p^\perp}(X-p)=((\pi T_p)^{-1}\xi)^2$. Since $\mathcal B_u$ is $C^{1,\alpha}$, the restriction $\Pi_{P_p^\perp}|_{P_0^\perp}:P_0^\perp\to P_p^\perp$ is invertible.
Decomposing
\[
 X-p=\underbrace{(Dg(y)(y'-y),\,y'-y)}_{\in P_p}+(z-g(y)-Dg(y)(y'-y),\,0)
\]
and applying $\Pi_{P_p^\perp}$, we  express the original normal coordinate $z$ in terms of $(\xi,y')$:
\begin{equation}
 z=\varphi_p(\xi,y') =g(y)+Dg(y)(y'-y) +\underbrace{(\Pi_{P_p^\perp}|_{P_0^\perp})^{-1} ((\pi T_p)^{-1}\xi)^2}_{=:q(y,\xi)}.
 \label{eqn:transformed-model-1}
\end{equation}
The remaining $V$-component of the tangent model is
\begin{equation}
 V_p(\xi,y')=\pi^\perp T_p(\pi T_p)^{-1}\xi=:M(y)\xi .
 \label{eqn:transformed-model-2}
\end{equation}
By non-degeneracy, the homogeneous quadratic polynomial $q(y,\cdot)$ satisfies $c|\xi|\leqslant|D_\xi q(y,\xi)|\leqslant C|\xi|$, uniformly in $y$. We call $U_p=(\varphi_p,V_p)$ the \emph{transformed model} at $p$. Its components $\varphi_p$ and $V_p$ are even and odd in $\xi$, respectively, with $\varphi_p(0,y')=g(y)+Dg(y)(y'-y)$ and $V_p(0,y')=0$. 

In the variables $(\xi,y')$, it is an affine function of $y'$ plus a homogeneous polynomial in $\xi$, whose coefficients $g(y),Dg(y),q(y,\cdot),M(y)$ are frozen at $y$. We shall use $U_{p}$ only for fixed $y$, as a comparison model near $(0,y)$, and postpone the study of the dependence of the frozen coefficients on $y$ to Section~\ref{sec:tangential-regularity}.
\end{example}

Since $|\Pi_{P_p^\perp}(X-p)|\sim|\xi|^2$, the rescaling $X\mapsto r^2X$ in the original coordinates corresponds to the anisotropic rescaling $(\xi,y)\mapsto(r\xi,r^2y)$ in the new coordinates.

Let $f(X,u,Du)$ be the $\mathscr J$-energy density of Section~\ref{sec:minimality-compactness}. For a field $U=(\varphi,V)$ with $\det D_\xi\varphi>0$ away from $\{\xi=0\}$, define the partial Legendre transform by \eqref{eqn:legendre-coordinate-maps}. Changing variables in $\mathscr J$ implies the transformed density
\begin{equation}
 \hat f(\xi,y,U,DU)=(\det D_\xi\varphi)\, f(X,\widehat u,D\widehat u\,(DX)^{-1}),\qquad \xi\neq0.
 \label{eqn:exact-density}
\end{equation}
\begin{example}[The Euclidean case]
\label{ex:euclidean-density}
For the Dirichlet density $\frac{1}{2}|Du|^2$,
\begin{equation}
 \hat f_0=\frac{\det D_\xi\varphi}{2}|D\widehat u\,(DX)^{-1}|^2,
 \qquad D\widehat u\,(DX)^{-1}
 =\begin{pmatrix}
 (D_\xi\varphi)^{-1}&-(D_\xi\varphi)^{-1}D_y\varphi\\ D_\xi V(D_\xi\varphi)^{-1}&D_yV-D_\xi V(D_\xi\varphi)^{-1}D_y\varphi
 \end{pmatrix}.
 \label{eqn:euclidean-density}
\end{equation}
\end{example}

\subsection{Weighted spaces}
\label{sec:weighted-spaces}
We use the Euclidean structure in $(\xi,y)$ for all integrals and weak equations. The degeneracy at $\{\xi=0\}$ is encoded in the weighted derivatives and norms below. We write $\rho=|\xi|$ and work on the anisotropic cylinders
\[
 Q_r(y_0)=B_r^2\times B_{r^2}^{n-2}(y_0),\qquad  |Q_r(y_0)|=r^{N_G}|Q_1|,\qquad N_G=2n-2.
\]

A field $Z=(Z_\varphi,Z_V)$ is \emph{admissible} if $Z_\varphi$ is even
and $Z_V$ is odd in $\xi$. Its \emph{Grushin (weighted) derivatives} are defined by
\begin{equation}
 \gradG Z=(\rho^{-1}D_\xi Z_\varphi,\ D_yZ_\varphi,\
           D_\xi Z_V,\ \rho D_yZ_V),
 \label{eqn:weighted-derivative-blocks}
\end{equation}
which are of degree zero under the
rescaling
\[
 (\mathbf S_rZ)(\xi,y)=(r^{-2}Z_\varphi,\,r^{-1}Z_V)(r\xi,\,y_0+r^2y).
\]
The $L^2$-height excess on $Q_r(y_0)$ is defined by
\begin{equation}
 \mathscr E_r(Z) = r^{-N_G-4}\int_{Q_r(y_0)}|Z_\varphi|^2 +r^{-N_G-2}\int_{Q_r(y_0)}|Z_V|^2,
 \label{eqn:excess}
\end{equation}
so that $\mathscr E_r(Z)=\mathscr E_1(\mathbf S_rZ)$ and $\int_{Q_1}|\gradG\mathbf S_rZ|^2=r^{-N_G}\int_{Q_r(y_0)}|\gradG Z|^2$. On $\Omega=Q_r(y_0)$, we define the adapted Sobolev norms as
\begin{equation}
 \|Z\|_{L^2_G(\Omega)}^2=\int_\Omega(|Z_\varphi|^2+\rho^2|Z_V|^2), \qquad \|Z\|_{W_G^{1,2}(\Omega)}^2 =\int_\Omega(|Z_\varphi|^2+|Z_V|^2+|\gradG Z|^2).
 \label{eqn:weighted-energy-norm}
\end{equation}
Write $W_G^{1,2}(\Omega)$ for the completion of the smooth admissible fields under the norm above, and $W_{G,0}^{1,2}(\Omega)$ for the closure of those compactly supported ones. Here, the trace $Z_\varphi(0,y)$ need not vanish. 

We shall study the second variation of the transformed $\mathscr J$-energy $\int\hat f(\xi,y,U,DU)$. For this, we introduce the following
\begin{equation}
\begin{aligned}
 \mathbf A=(\mathbf A^1,\mathbf A^0) =(D_{DU}\hat f,D_U\hat f), \qquad \mathbf A_UZ
 =\left.\frac{d}{dt}\right|_{t=0} \mathbf A(\xi,y,U+tZ,DU+tDZ),
\end{aligned}
\label{eqn:linearised-coefficients}
\end{equation}
with the associated bilinear form
\begin{equation}
 \mathscr B_U(Z,Y)
 =\int_\Omega((\mathbf A_UZ)^1:DY+(\mathbf A_UZ)^0\cdot Y),
 \label{eqn:linearised-weak-form}
\end{equation}
for admissible $Z\in W_G^{1,2}(\Omega)$ and smooth admissible $Y$ compactly supported in $\Omega$. Here $\mathbf A^1$ and $\mathbf A^0$ have the same
structure as $DU$ and $U$, respectively, with the block decomposition
$\mathbf A^1:DY=\mathbf A^1_{\varphi,\xi}:D_\xi Y_\varphi+\mathbf A^1_{\varphi,y}:D_yY_\varphi +\mathbf A^1_{V,\xi}:D_\xi Y_V+\mathbf A^1_{V,y}:D_yY_V$
and $\mathbf A^0\cdot Y=\mathbf A^0_\varphi\cdot Y_\varphi
+\mathbf A^0_V\cdot Y_V$. The weights dual to \eqref{eqn:weighted-derivative-blocks} are
\[
 \mathbf A^1_G
 =(\rho\mathbf A^1_{\varphi,\xi},\ \mathbf A^1_{\varphi,y},\ \mathbf A^1_{V,\xi},\ \rho^{-1}\mathbf A^1_{V,y}),  \qquad\text{so that}\qquad
 \mathbf A^1:DY=\mathbf A^1_G\cdot\gradG Y.
\] For such coefficient fields $\mathbf{A}=(\mathbf A^1,\mathbf A^0)$, we use script letters for the corresponding functionals and introduce the dual weighted (semi-)norm
\begin{equation}
    \label{eqn:functionals}
     \mathscr A(Y)=\int_\Omega(\mathbf A^1:DY+\mathbf A^0\cdot Y),
 \qquad \|\mathbf A\|_{(W_G^{1,2}(\Omega))^*} =\sup_{\substack{Y\in C_c^\infty(\Omega)\text{ admissible}\\ \|\gradG Y\|_{L^2(\Omega)}\leqslant1}}|\mathscr A(Y)|.
\end{equation}
 
 The following elementary inequalities are used throughout Part~II.
\begin{lemma}[Weighted inequalities]
\label{lem:weighted-elementary-inequalities}
There hold: 
\begin{enumerate}[label=\textnormal{(\roman*)}]
\item Every admissible field $Z$ on $Q_r(y_0)$ satisfies
\begin{equation}
 \int\rho^{-2}|Z_V|^2\leqslant\int|D_\xi Z_V|^2 .
 \label{eqn:hardy-trace}
\end{equation}
If further $Z$ is compactly supported in $Q_r(y_0)$, then
\begin{equation}
 \|Z_\varphi\|_{L^2}\leqslant Cr^2\|\gradG Z\|_{L^2},\qquad
 \|Z_\varphi(0,\cdot)\|_{L^2_y}\leqslant Cr\|\gradG Z\|_{L^2}.
 \label{eqn:weighted-poincare-inequality}
\end{equation}
\item If $|\mathbf A^1_G|\leqslant C_1$, $|\mathbf A^0_\varphi|\leqslant C_0$ and $|\mathbf A^0_V|\leqslant C_0\rho$ on $Q_r(y_0)$, then
\begin{equation}
 \|\mathbf A\|_{(W_G^{1,2}(Q_r(y_0)))^*}\leqslant Cr^{N_G/2}(C_1+r^2C_0).
 \label{eqn:source-dual}
\end{equation}
\end{enumerate}
In particular, $\|Z\|_{W_G^{1,2}}^2$ is uniformly equivalent to $\int(|Z_\varphi|^2+|\gradG Z|^2)$, and to $\int|\gradG Z|^2$ on $W_{G,0}^{1,2}$; the trace $Z\mapsto Z_\varphi(0,\cdot)$ is bounded from $W_G^{1,2}(Q_r(y_0))$ to $L^2(B_{r^2}^{n-2}(y_0))$.
\end{lemma}

\begin{proof}
(i) On each normal circle, an odd function has only odd Fourier modes, so Wirtinger's
inequality implies
\[
 \int_0^{2\pi}|Z_V|^2\,d\theta\leqslant\int_0^{2\pi}|\partial_\theta Z_V|^2\,d\theta
 \leqslant\rho^2\int_0^{2\pi}|D_\xi Z_V|^2\,d\theta .
\]
Integrating in $\rho\,d\rho\,dy$ proves \eqref{eqn:hardy-trace}. The first bound in \eqref{eqn:weighted-poincare-inequality} is the standard Poincar\'e inequality in $y$ on $B_{r^2}^{n-2}(y_0)$. To prove $Z_\varphi(0,\cdot)\in L_y^2$, we integrate along normal rays from $\rho=r$, where $Z$ vanishes:
\[
 |Z_\varphi(0,y)|^2
 =\left|\int_0^r\partial_\rho Z_\varphi\,d\rho\right|^2
 \leqslant\frac{r^2}{2}\int_0^r\rho^{-1}|\partial_\rho Z_\varphi|^2\,d\rho,
\]
and average in $\theta$ and integrate in $y$. Then (ii) follows from (i) and Cauchy--Schwarz.
\end{proof}
\smallskip
\noindent\textbf{Weighted H\"older norms.}
Although the equations degenerate along \(\{\xi=0\}\), they are uniformly elliptic on balls whose radii are proportional to their distance from this set, with ellipticity constants independent of that distance. More precisely, fix \(c_0\in(0,1/100)\). For each \((\xi,y)\in\Omega\) with \(\rho=|\xi|>0\), we define $Q_{\xi,y}=(B_{c_0\rho}^2(\xi)\times B_{c_0^2\rho^2}^{n-2}(y))\cap\Omega$, on which $|\xi|$ is comparable to $\rho$.  

We measure fields $Z$, degree-zero quantities $F$ (such as the blocks of
$\gradG Z$), and coefficient fields $\mathbf A$ by adapted H\"older norms
\begin{equation}
\begin{aligned}
 \|Z\|_{C_G^{k,\alpha}(\Omega)}
 &=\sup\|(\rho^{-2}Z_\varphi,\ \rho^{-1}Z_V)
       \|_{C^{k,\alpha}(Q_{\xi,y},g_\rho)}, \\
 \|F\|_{C_G^{k,\alpha}(\Omega)}
 &=\sup\|F\|_{C^{k,\alpha}(Q_{\xi,y},g_\rho)},\\
 \|\mathbf A\|_{C_G^{k,\alpha}(\Omega)}
 &=\sup\|(\mathbf A_G^{1},\
       \rho^2\mathbf A^0_\varphi,\
       \rho\mathbf A^0_V)
       \|_{C^{k,\alpha}(Q_{\xi,y},g_\rho)}.
\end{aligned}
\label{eqn:weighted-holder-field}
\end{equation}
Here each supremum is over $(\xi,y)\in\Omega$ with $\rho=|\xi|>0$, and the $C^{k,\alpha}$-norms are taken with respect to the \textit{rescaled product metric} $g_\rho=\rho^{-2}g_{\xi}+\rho^{-4}g_{y}$.

\smallskip
\noindent\textbf{Comparison polynomials.}
The transformed models \eqref{eqn:transformed-model-1}--\eqref{eqn:transformed-model-2} lie in the finite-dimensional space with norms
\begin{equation}
\begin{aligned}
    & \mathcal P=\{P=(a+B(y-y_0)+q(\xi),\ M\xi)\},\\
    & \|\nabla P\|=|B|+|q|+|M|,\qquad \|P\|_r=r^{-2}|a|+\|\nabla P\|,
\end{aligned} \label{eqn:polynomial-norm}
\end{equation}
where $a\in\C$, $B:\R^{n-2}\to\C$ and $M:\C\to\operatorname{im}\pi^\perp$ are real-linear, and $q$ is a homogeneous quadratic polynomial. Denote the norm of the coefficients of $q$ by $|q|$. Since $\gradG P=(\rho^{-1}D_\xi q,B,M,0)$ and $\mathbf S_rP=(r^{-2}a+By+q(\xi),M\xi)$,
\begin{equation}
 \|\nabla P\|\sim\|\gradG P\|_{L^2(Q_1(y_0))},\qquad \|P\|_r\sim\mathscr E_r(P)^{1/2},\qquad \forall P\in\mathcal P.
 \label{eqn:polynomial-norm-equivalence}
\end{equation}
For a field $Z$ defined near $(0,y_0)$, we call $P\in\mathcal P$ the \emph{comparison polynomial} of $Z$ at $(0,y_0)$ if $\mathscr E_r(Z-P)\to0$ as $r\to0$. It is unique if it exists: by \eqref{eqn:polynomial-norm-equivalence}, two such polynomials satisfy
$\|P-P'\|_r\to0$, so all coefficients of $P-P'$ vanish.

\subsection{Initial estimates}
\label{sec:legendre-initial}
The quantitative estimates of Part~I translate into the following estimates in Legendre coordinates:
\begin{proposition}[Initial estimates in Legendre coordinates]
\label{prop:legendre-initial-estimates}
For $r_0\ll 1$,
\begin{enumerate}
[label=\textnormal{(\roman*)}]
\item The partial Legendre transform $(z,y)\mapsto(\pi u(z,y),y)$ is well-defined.
\item For every $(\xi,y)$ with $0<\rho=|\xi|\leqslant r_0$ and $|y-y_0|\leqslant r_0$,
\begin{equation}
 \|(\rho^{-2}(\varphi-\varphi_{p(y)}),\ \rho^{-1}(V-V_{p(y)}))\|_{C^{1,\alpha}(Q_{\xi,y},g_\rho)} \leqslant C\rho^\alpha .
 \label{eqn:initial-holder-estimate}
\end{equation}
\item For every branch point $p_0=p(y_0)$ and $0<r\leqslant r_0$,
\begin{equation}
 \sup_{Q_r(y_0)\setminus\{\xi=0\}} (r^{-2}|\varphi-\varphi_{p_0}|+r^{-1}|V-V_{p_0}|    +|\gradG(U-U_{p_0})|)\leqslant Cr^\alpha .
 \label{eqn:initial-bounds}
\end{equation}
In particular, the comparison polynomial of $U$ at $(0,y)$ is the transformed model $U_{p(y)}$, with  estimates
\begin{equation}
 \mathscr E_r(U-U_{p_0})\leqslant Cr^{2\alpha}.
 \label{eqn:initial-excess}
\end{equation}

\end{enumerate}
\end{proposition}

\subsection{The transformed equation}
\label{sec:transformed-equation}
We now determine the structure of the transformed density \eqref{eqn:exact-density} in the weighted variables, and derive the weak
equation.
\begin{lemma}[The transformed equation]
\label{lem:transformed-equation}
Consider fields $U=(\varphi,V)$ with \emph{uniformly bounded normalised data}, \ie $(X,\widehat u/\rho,\gradG U)$ is uniformly bounded on $\{\xi\neq0\}$ and $\rho^{-1}D_\xi\varphi$ is uniformly invertible.
\begin{enumerate}[label=\textnormal{(\roman*)}]
\item The transformed density has the form
\begin{equation}
 \hat f=\hat f_{\mathrm{pr}}(X,\gradG U)+\rho^2\operatorname{Err}(\rho,X,\widehat u/\rho,\gradG U),
 \label{eqn:density-powers}
\end{equation}
where $\hat f_{\mathrm{pr}}$ is independent of $\widehat u/\rho$. Both $\hat f_{\mathrm{pr}}$ and $\operatorname{Err}$ are smooth (resp. real-analytic) functions of $\rho$ and the normalised data for smooth (resp. real-analytic) ambient data.
\item The coefficient field $\mathbf A$ in \eqref{eqn:linearised-coefficients} satisfies
\begin{equation}
 |\mathbf A^1_G|+|\mathbf A^0_\varphi|+\rho^{-1}|\mathbf A^0_V|\leqslant C.
 \label{eqn:weighted-coefficient-bounds}
\end{equation}
\item The transformed spinor $U$ and its transformed models $U_p$ have uniformly
bounded normalised data, with bounds independent of $p$, and $U$ satisfies the weak
equation
\begin{equation}
 \int(\mathbf A^1:DY+\mathbf A^0\cdot Y)=0
 \label{eqn:full-weak-equation}
\end{equation}
for every smooth admissible $Y$ compactly supported in $Q_{r_0}(y_0)$.
\end{enumerate}
\end{lemma}
\begin{proof}
(i) In a smooth trivialisation, $f(X,u,Du)$ is the sum of a quadratic
form in $Du$, a bilinear term in $(u,Du)$, and a quadratic form in $u$,
with smooth coefficients in $X$ (\cf the $\mathscr J$-energy \eqref{eqn:jacobi-functional}).
Writing $D_\xi\varphi=\rho\cdot\rho^{-1}D_\xi\varphi$ in
\eqref{eqn:euclidean-density}, each block of $Du=D\widehat u\,(DX)^{-1}$
is $\rho^{-1}$ times a smooth function of $\gradG U$, and
$\det D_\xi\varphi=\rho^2\det(\rho^{-1}D_\xi\varphi)$. Then writing
$\widehat u=\rho\cdot\widehat u/\rho$, the first term contributes
$\hat f_{\mathrm{pr}}(X,\gradG U)$, the transform of the quadratic form in $Du$ frozen at
$X$ (for the Euclidean metric, $\hat f_{\mathrm{pr}}$ is the density $\hat f_0$ in \eqref{eqn:euclidean-density}, a smooth function of $\gradG U$), with $O(\rho^2)$-remainders from the other parts. 

(ii) Since $\hat f$ depends on $DU$ only through $\gradG U$, we have $\mathbf A^1_G=D_{\gradG U}\hat f$, which is bounded. The field $\varphi$ enters only through $X$, so $\mathbf A^0_\varphi$ is bounded; $V$ enters only through $\widehat u/\rho$ in the $O(\rho^2)$-remainder of \eqref{eqn:density-powers}, so $\mathbf A^0_V=O(\rho)$.

(iii) The first assertion follows from the initial bounds \eqref{eqn:initial-bounds}, since $\rho^{-1}D_\xi\varphi_p=\rho^{-1}D_\xi q$ is uniformly invertible. To derive the weak equation, consider in the original coordinates the combined domain and target variation
\[
 u_t(\Psi_t(z,y))=u(z,y)+t(0,\,Y_V(\pi u,y)),\qquad
 \Psi_t(z,y)=(z+tY_\varphi(\pi u,y),\,y).
\]
We check that this is an admissible variation in the sense of Remark~\ref{rmk:rescaled-geometry}:
\begin{itemize}
\item (\textit{inner variation}) Since $Y_\varphi$ is even in $\xi$, $\Psi_t$ is single-valued. Moreover, $D_\xi Y_\varphi=O(\rho)$ and $D_X\xi=O(\rho^{-1})$ guarantees that for $|t|\ll 1$ small, $\Psi_t$ is a bi-Lipschitz homeomorphism, equal to the identity outside a compact set.
\item (\textit{outer variation}) Since $Y_V$ is odd, $u_t$ is a two-valued section.
\end{itemize}
The variation preserves $\pi u$, so $(\xi,y)=(\pi u_t\circ \Psi_t,y)=(\pi u,y)$ are Legendre coordinates for $u_t$, in which the transformed field is $U_t=U+tY$. Therefore, the inner and outer variation identities of Remark~\ref{rmk:rescaled-geometry} apply to $u_t$, and by \eqref{eqn:exact-density}, the local $\mathscr J$-minimality of $u$ yields
\[
 0=\frac{d}{dt}\Big|_{t=0}\int\hat f(\xi,y,U+tY,DU+tDY)
 =\int(\mathbf A^1:DY+\mathbf A^0\cdot Y).
 \qedhere
\]
\end{proof}
Let $\mathbf A_{\mathrm{pr}}(X,DU)=(D_{DU}\hat f_{\mathrm{pr}}(X,\gradG U),0)$ be the coefficient field of the principal part of \eqref{eqn:density-powers}.

\section{The frozen operator}
\label{sec:frozen-operator}
In Section~\ref{sec:linearised-polynomial-approximation}, we shall treat the second variation $\mathscr B_U$ as a perturbation of the second variation at a transformed model, with coefficients frozen at its centre. In this section we study this frozen second variation: we prove its coercivity and an excess decay for its weak solutions (\cf Proposition~\ref{prop:frozen-operator}), and construct the polynomial corrections (\cf Lemma~\ref{lem:polynomial-corrections}) used in the Campanato estimate (\cf Theorem~\ref{thm:campanato}) to absorb the $O(1)$ principal parts of the source term.

Fix a branch point $p_0=p(y_0)$. Freezing the principal part of $\mathscr J$ at $p_0$ yields the density
\[
 f_{p_0}(Du)=\tfrac{1}{2}\sqrt{\det g_{ij}(p_0)}\,|Du|^2_{g_{ij}(p_0),h(p_0)},
\]
where $g_{ij}(p_0),h(p_0)$ are the Riemannian metric and the spinor bundle metric at $p_0$. We write $\hat f_{p_0}(DU)$ for its transform \eqref{eqn:exact-density}; it depends only on $DU$, and $\hat f_{p_0}(DU)=\hat f_{\mathrm{pr}}(p_0,\gradG U)$. By Theorem~\ref{thm:spinor-blowups}, $\Phi_{p_0}$ minimises the energy of $f_{p_0}$, so $U_{p_0}$ is a critical point of $\int\hat f_{p_0}$. The frozen linearised coefficient field and second variation are
\begin{equation}
\begin{aligned}
 \mathbf A_{\mathrm{pr},p_0}Z&=(D^2\hat f_{p_0}(DU_{p_0})(DZ),\,0)=D_{DU}\mathbf A_{\mathrm{pr}}(p_0,DU_{p_0})(DZ),\\
 \mathscr B_{p_0}(Z,Y)&=(\mathscr A_{p_0}Z)(Y)
 =\int D^2\hat f_{p_0}(DU_{p_0})(DZ,DY).
 \label{eqn:frozen-hessian}
 \end{aligned}
\end{equation}

\begin{proposition}[The frozen operator]
\label{prop:frozen-operator}
There are constants $c,C>0$ and $\theta_0\in(0,1)$, uniform in $p_0$,
with the following properties.
\begin{enumerate}[label=\textnormal{(\roman*)}]
\item \textnormal{(Coercivity)} For every $Z\in W_{G,0}^{1,2}$,
\begin{equation}
 c\|\gradG Z\|_{L^2}^2\leqslant\mathscr B_{p_0}(Z,Z)
 \leqslant C\|\gradG Z\|_{L^2}^2.
 \label{eqn:linearised-coercivity}
\end{equation}
\item \textnormal{(Excess decay)} Let $Z\in W^{1,2}_{G,\mathrm{loc}}(Q_1(y_0))$
satisfy $\mathscr E_1(Z)<\infty$ and $\mathscr B_{p_0}(Z,Y)=0$ for all
compactly supported admissible $Y$. Then there is a comparison polynomial
$H\in\mathcal P$ (see \eqref{eqn:polynomial-norm}) with $\mathscr A_{p_0}H=0$ such that
\begin{equation}
 \|H\|_1\leqslant C\mathscr E_1(Z)^{1/2},\qquad
 \mathscr E_\theta(Z-H)\leqslant C\theta^4\mathscr E_1(Z),
 \qquad 0<\theta<\theta_0.
 \label{eqn:frozen-excess-decay}
\end{equation}
\end{enumerate}
\end{proposition}

The proof goes through the rest of this section: we reduce to the square model in Section~\ref{sec:square-model}, and prove (i) in Section~\ref{sec:frozen-coercivity} and (ii) in Section~\ref{sec:frozen-decay}. The polynomial corrections are treated in Section~\ref{sec:polynomial-corrections}.

\subsection{Reduction to the square model}
\label{sec:square-model}
\begin{example}[The square model]
\label{ex:square-model}
Let $\hat f_0$ be the Euclidean density \eqref{eqn:euclidean-density} and
$\mathscr B_0$ the second variation at $U_0$. Then
\begin{equation}
 \mathscr B_0(Z,Z)=\int(\rho^{-2}|\partial_{\bar\xi}Z_\varphi|^2
 +|D_yZ_\varphi|^2+|D_\xi Z_V|^2+4\rho^2|D_yZ_V|^2),
 \label{eqn:canonical-frozen-form}
\end{equation}
where $\partial_{\xi}=\frac{1}{2}(\partial_{\xi_1}-\sqrt{-1}\,\partial_{\xi_2})$ and
$\partial_{\bar \xi}=\frac{1}{2}(\partial_{\xi_1}+\sqrt{-1}\,\partial_{\xi_2})$ for
$\xi=\xi_1+\sqrt{-1}\,\xi_2\in\mathbb C$. Indeed, since
$|D_\xi\varphi|^2=2(|\partial_\xi\varphi|^2+|\partial_{\bar\xi}\varphi|^2)$ and
$\det D_\xi\varphi=|\partial_\xi\varphi|^2-|\partial_{\bar\xi}\varphi|^2$, we find
(\cf the top-left block in \eqref{eqn:euclidean-density})
\[
 \frac{\det D_\xi\varphi}{2}|(D_\xi\varphi)^{-1}|^2
 =\frac{|D_\xi\varphi|^2}{2\det D_\xi\varphi}
 =1+\frac{2|\partial_{\bar\xi}\varphi|^2}{|\partial_\xi\varphi|^2-|\partial_{\bar\xi}\varphi|^2}.
\]
Along the variation $U_0+tZ$, we have $\partial_\xi\varphi=2\xi+t\partial_\xi Z_\varphi$ and
$\partial_{\bar\xi}\varphi=t\partial_{\bar\xi}Z_\varphi$, and substituting into
\eqref{eqn:euclidean-density} yields
\[
 \hat f_0(DU_0+tDZ)=1+\frac{t^2}{2}(\rho^{-2}|\partial_{\bar\xi}Z_\varphi|^2
 +|D_yZ_\varphi|^2+|D_\xi Z_V|^2+4\rho^2|D_yZ_V|^2)
 +O(t^3).
\]
\end{example}

\begin{lemma}[Reduction to the square model]
\label{lem:frozen-form-conjugacy}
There exist an isometry $(z,y)\mapsto(z',y')$ of the base, sending $p_0$ to $(0,y_0)$ and $P_{p_0}$ to $\{z'=0\}$, and a homothety of the fibre $E_0$ with factor $\kappa^{-1/2}=|T_{p_0}(1)|^{-1}$, which together transform $U_{p_0}$ into the square model $U_0$. Let $(\xi',y')=\Psi(\xi,y)$ be the induced change of Legendre coordinates, and $Z\mapsto Z'$ be the linearisation at $U_{p_0}$ of the induced map on fields. Then:
\begin{enumerate}[label=\textnormal{(\roman*)}]
\item $\mathscr B_{p_0}(Z,Y)=\kappa\,\mathscr B_0(Z',Y')$ for admissible $Z$ and compactly supported admissible $Y$;
\item $\|\gradG Z\|_{L^2}$, $\|Z\|_{L^2_G}$ and $\mathscr E_r(Z)$ are uniformly comparable to the corresponding quantities of $Z'$ on the image cylinders;
\item $Z\mapsto Z'$ preserves admissibility, compact support, $\mathcal P$ (see \eqref{eqn:polynomial-norm}) and fields of the form $(q(\xi),0)$; $\Psi$ commutes with anisotropic dilations, and $Q_{cr}(y_0)\subset\Psi(Q_r(y_0))\subset Q_{Cr}(y_0)$ with uniform $c,C>0$.
\end{enumerate}
\end{lemma}

\begin{proof}
The tangent model is determined by the plane $P_{p_0}$ and the real-linear map $T=T_{p_0}:\C\to E_0$, $T(w)=a_{p_0}w+b_{p_0}\bar w$. We normalise the former by an isometry of the base and the latter by a homothety of the fibre, and read off the induced \textit{non-linear} change of Legendre coordinates. 

Since $a_{p_0}w$ and $b_{p_0}\bar w$ lie in orthogonal eigenspaces of $J$ (\cf Definition~\ref{def:half-frequency-models}), we have $|T(w)|^2=\kappa|w|^2$. Let $G$ be a linear isometry from $(\R^n,g(p_0))$ to Euclidean space sending $P_{p_0}$ to $\{z'=0\}$, so that $z'=\Pi_{P_{p_0}^\perp}(X-p_0)$ under the identification $P_{p_0}^\perp\cong\C$ of Section~\ref{sec:partial-legendre}, and set $X'=G(X-p_0)+(0,y_0)$. Let $Q$ be an orthogonal map of $E_0$ with $QT(w)=(\kappa^{1/2}w,0)$ in the decomposition $(\pi,\pi^\perp)$, and set $(\xi',V')=H(\xi,V)$ with $H=\kappa^{-1/2}Q$. Then the model becomes $\pm((z')^{1/2},0)$, \ie $z'=(\xi')^2$ and $V'=0$ in the new Legendre coordinates.

In block form with respect to $(z,y)$ and $(\pi,\pi^\perp)$, this reads
\[
 G_{11}Dg(y_0)+G_{12}=0,\qquad H_{11}+H_{12}M=(\pi T)^{-1},
\]
\[
H_{21}+H_{22}M=0,\qquad \text{and} \qquad G_{11}q(y_0,\xi)=((\pi T)^{-1}\xi)^2,
\]
where $M=M(y_0)$. For a field $U$, the new Legendre coordinates are $\Psi_U=(H_{11}\xi+H_{12}V,\ y_0+G_{21}(\varphi-g(y_0))+G_{22}(y-y_0))$, and at $U=U_{p_0}$,
\begin{equation}
 \Psi(\xi,y)=((\pi T)^{-1}\xi,\ y_0+(G_{21}Dg(y_0)+G_{22})(y-y_0)+G_{21}q(y_0,\xi)).
 \label{eqn:full-transform}
\end{equation}
Linearising the induced map $U\mapsto U'$ at $U_{p_0}$ yields
\[
 Z'\circ\Psi=\begin{pmatrix}G_{11}&-2\xi'H_{12}\\0&H_{22}\end{pmatrix}Z .
\]

We can now prove (i)--(iii). Since $G$ and $\kappa^{1/2}H$ are isometries, $f_{p_0}(Du)=\kappa |\sqrt{\det g(p_0)}| f_0(Du')$, and changing variables gives $\int_\Omega\hat f_{p_0}(DU)=\kappa\int_{\Psi_U(\Omega)}\hat f_0(DU')$ for fields $U$ near $U_{p_0}$. Taking the second variation at $U_{p_0}$ proves (i).

By \eqref{eqn:full-transform}, $\Psi$ is triangular with constant Jacobian, preserves parity and commutes with the anisotropic dilations. Since $\xi'$ and $Z_V$ are odd, the matrix above preserves admissibility, $\mathcal P$ and fields of the form $(q(\xi),0)$. This proves (iii).

For (ii), we have $\rho'\sim\rho$, and the only new terms are $2\xi'H_{12}Z_V=O(\rho|Z_V|)$ in $Z'_\varphi$ by the chain rule, which is compatible with the weights in $\|\cdot\|_{L^2_G}$ and $\mathscr E_r$, and $\rho^{-1}D_\xi(2\xi'H_{12}Z_V)=O(\rho^{-1}|Z_V|+|D_\xi Z_V|)$ in $\gradG Z'$, whose $L^2$-norm is bounded by $\|D_\xi Z_V\|_{L^2}$ by \eqref{eqn:hardy-trace}.
\end{proof}

\subsection{Coercivity}
\label{sec:frozen-coercivity}
To prove Proposition ~\ref{prop:frozen-operator}\textnormal{(i)}, it remains to bound $\int\rho^{-2}|\partial_{\bar\xi}Z_\varphi|^2$ from below. We start with a Hardy-type estimate for $Z_\varphi$. For such an even $Z_\varphi$, regarded as $\C$-valued, we write
\begin{equation}
 f(\xi,y)=\frac{Z_\varphi(0,y)-Z_\varphi(\xi,y)}{2\xi}.
 \label{eqn:normal-difference-quotient}
\end{equation}
Then $f$ is odd in $\xi$, with $f=O(\rho)$ and $D_\xi f=O(1)$, and
$\partial_{\bar\xi}f=-\partial_{\bar\xi}Z_\varphi/(2\xi)$, so
\begin{equation}
 \rho^{-2}|\partial_{\bar\xi}Z_\varphi|^2=4|\partial_{\bar\xi}f|^2 =|D_\xi f|^2-2\det D_\xi f .
 \label{eqn:difference-quotient-identity}
\end{equation}
We then prove the a Hardy-type inequality for the even component $Z_\varphi$.
\begin{lemma}\label{lem:even-hardy}
Every smooth $Z_\varphi$ on $Q_r(y_0)$, even in $\xi$, satisfies
\begin{equation}
 \int_{Q_r(y_0)}\frac{|Z_\varphi(\xi,y)-Z_\varphi(0,y)|^2}{\rho^4} \leqslant\int_{Q_r(y_0)}\frac{|\partial_\rho Z_\varphi|^2}{\rho^2},
 \qquad \|D_\xi f\|_{L^2}+r^{-1}\|f\|_{L^2}\leqslant C\|\rho^{-1}D_\xi Z_\varphi\|_{L^2}.
 \label{eqn:radial-hardy}
\end{equation}
\end{lemma}
\begin{proof}
Note that $Z_\varphi-Z_\varphi(0,y)=O(\rho^2)$ by evenness. For fixed $(\theta,y)$ and $0<s<r$, integration by parts and Young's inequality give
\begin{align*}
     \int_0^s\frac{|Z_\varphi-Z_\varphi(0,y)|^2}{\rho^3}\,d\rho
 & \leqslant \int_0^s\frac{\operatorname{Re}\langle Z_\varphi-Z_\varphi(0,y),
 \partial_\rho Z_\varphi\rangle}{\rho^2}\,d\rho \\
 & \leqslant\frac{1}{2}\int_0^s\frac{|Z_\varphi-Z_\varphi(0,y)|^2}{\rho^3}\,d\rho
 +\frac{1}{2}\int_0^s\frac{|\partial_\rho Z_\varphi|^2}{\rho}\,d\rho .
\end{align*}
Integrating in $(\theta,y)$ and letting $s\to r$ proves the first inequality. The second follows from it and $|D_\xi f|+r^{-1}|f|\lesssim\rho^{-1}|D_\xi Z_\varphi| +\rho^{-2}|Z_\varphi-Z_\varphi(0,y)|$.
\end{proof}

\begin{proof}[Proof of Proposition~\ref{prop:frozen-operator}\textnormal{(i)}]
By Lemma~\ref{lem:frozen-form-conjugacy}, it suffices to consider the square model. The upper bound is clear from \eqref{eqn:canonical-frozen-form}. For the lower bound, we aim to bound $\int\rho^{-2}|D_\xi Z_\varphi|^2$ by $\int\rho^{-2}|\partial_{\bar\xi}Z_\varphi|^2$ for smooth compactly supported
admissible $Z$. Fix $y$ and let $f$ be as in \eqref{eqn:normal-difference-quotient}; outside the support, $f=Z_\varphi(0,y)/(2\xi)=O(\rho^{-1})$. The determinant in \eqref{eqn:difference-quotient-identity} integrates to zero: on an annulus $\{\varepsilon<\rho<s\}$,
\[
 \int\det D_\xi f\,d\xi =\frac{1}{2}\int_0^{2\pi}\operatorname{Im}(\bar f\,\partial_\theta f)\,d\theta \Big|_{\rho=\varepsilon}^{\rho=s} =O(\varepsilon^2)+O(s^{-2}).
\]
Hence $\int_{\R^2}|D_\xi f|^2=\int_{\R^2}\rho^{-2}|\partial_{\bar\xi}Z_\varphi|^2$. Since $D_\xi Z_\varphi=-2(f\,d\xi+\xi D_\xi f)$ and $f$ is odd, the weighted inequality in \eqref{eqn:hardy-trace} yield
\[
 \int_{\R^2}\rho^{-2}|D_\xi Z_\varphi|^2 \lesssim\int_{\R^2}(|D_\xi f|^2+\rho^{-2}|f|^2) \lesssim\int_{\R^2}|D_\xi f|^2 =\int_{\R^2}\rho^{-2}|\partial_{\bar\xi}Z_\varphi|^2 .
\]
Integrating in $y$ completes the proof.
\end{proof}

\subsection{Excess decay}
\label{sec:frozen-decay}
\begin{proof}[Proof of Proposition~\ref{prop:frozen-operator}\textnormal{(ii)}]
By Lemma~\ref{lem:frozen-form-conjugacy}, we may assume $U_{p_0}=U_0$, and write $Q_a=Q_a(y_0)$. The strategy (essentially the same as in Part I) is to prove that $Z$ is smooth up to $\{\xi=0\}$, with bounded $C^{m,\alpha}$-norm, and take $H$ to be its Taylor polynomial in $\mathcal P$ (see \eqref{eqn:polynomial-norm}). We divide this proof in three steps. 

\medskip
\textbf{Step 1: Tangential regularity.}
For $0<a<b<1$, let $\eta=\eta_1(\rho)\eta_2(y)$ be a cutoff of $Q_a\subset Q_b$ with $\eta_1=1$ near $\rho=0$. Since the coefficients of $\mathscr B_0$ do not depend on $y$, we have the exact identity
\[
 \mathscr B_0(\eta Z,\eta Z)=\mathscr B_0(Z,\eta^2Z)
 +\int D^2\hat f_0(DU_0)(Z\otimes D\eta,Z\otimes D\eta).
\]
The first term vanishes. Coercivity (\cf Proposition \ref{prop:frozen-operator} (i)), the square-model form \eqref{eqn:canonical-frozen-form} and $D_\xi\eta=0$ near $\rho=0$ imply the weighted Caccioppoli inequality
\begin{equation}
 \int_{Q_a}|\gradG Z|^2\leqslant C_{a,b}\|Z\|_{L^2_G(Q_b)}^2 \leqslant C_{a,b}\mathscr E_1(Z).
 \label{eqn:frozen-local-energy}
\end{equation}
The equation is invariant under translations in $y$, so tangential difference quotients $\delta_{t,\alpha}Z(\xi,y)
   =\frac{Z(\xi,y+t\partial_{y_\alpha})-Z(\xi,y)}{t}$ of $Z$ are also weak solutions with $L^2_G$-norms bounded by $\|\gradG Z\|_{L^2}$. Applying \eqref{eqn:frozen-local-energy} to $\delta_{t,\alpha}Z$, we get
\begin{equation}
 \sum_{|\gamma|\leqslant m}(\|D_y^\gamma Z\|_{L^2(Q_a)}  +\|\gradG D_y^\gamma Z\|_{L^2(Q_a)}) \leqslant C_{m,a}\mathscr E_1(Z)^{1/2}.
 \label{eqn:frozen-tangential}
\end{equation}

\textbf{Step 2: Regularity across $\{\xi=0\}$.}
Let $f$ be the normal difference quotient \eqref{eqn:normal-difference-quotient} of $Z_\varphi$. Lemma~\ref{lem:even-hardy} and the trace bound of Lemma~\ref{lem:weighted-elementary-inequalities}, applied to cutoffs of $D_y^\gamma Z$, give
\begin{equation}
 \sum_{|\gamma|\leqslant m}(\|D_y^\gamma Z_\varphi(0,\cdot)\|_{L^2(B_{a^2}^{n-2}(y_0))} +\|D_y^\gamma f\|_{L^2(Q_a)}+\|D_\xi D_y^\gamma f\|_{L^2(Q_a)}) \leqslant C_{m,a}\mathscr E_1(Z)^{1/2}.
 \label{eqn:tangent-estimate-f}
\end{equation}
Testing the equation with $Y=(-2\xi\psi,0)$, with $\psi$ odd, we find that $f$ and $Z_V$ are weak solutions of
\begin{equation}
 \Delta_\xi f+4\rho^2\Delta_yf=2\bar\xi\,\Delta_yZ_\varphi(0,y),\qquad \Delta_\xi Z_V+4\rho^2\Delta_yZ_V=0,
 \label{eqn:quotient-equations}
\end{equation}
where odd test functions suffice by oddness of $f$ and $Z_V$. On each normal slice, the standard $W^{2,2}$-estimates imply
\[
 \|D_\xi^2f\|_{L^2(B_a^2)} \leqslant C_{a,b} (\|\Delta_\xi f\|_{L^2(B_b^2)} +\|f\|_{L^2(B_b^2)}).
\]
Using the equation and integrating in $y$, we obtain
\[
 \|D_\xi^2f\|_{L^2(Q_a)}\leqslant C_{a,b}( \|\Delta_yf\|_{L^2(Q_b)} +\|\Delta_y Z_\varphi(0,\cdot)\|_{L^2(B_{b^2}^{n-2}(y_0))} +\|f\|_{L^2(Q_b)})\leqslant C_{a,b}\mathscr E_1(Z)^{1/2},
\]
by \eqref{eqn:tangent-estimate-f}. The same bound holds for $D_\xi^2D_y^\gamma f$ and their higher order weighted derivatives. The same estimates apply to $Z_V$. Thus, by Sobolev embedding, we find
\[
 \|Z_\varphi\|_{C^m(Q_a)}+\|Z_V\|_{C^m(Q_a)}
 \leqslant C_{m,a}\mathscr{E}_1(Z)^{1/2}.
\]

\textbf{Step 3: Approximation by a kernel element.}
By parity, the Taylor polynomial of $Z$ at $(0,y_0)$ in $\mathcal P$ (see \eqref{eqn:polynomial-norm}) is
\[
 H=(Z_\varphi(0,y_0)+D_yZ_\varphi(0,y_0)(y-y_0)+\tfrac{1}{2}D_\xi^2Z_\varphi(0,y_0)(\xi,\xi),\ D_\xi Z_V(0,y_0)\xi),
\]
and Step~2 gives, on $Q_{\theta_0}$ with $\theta_0\ll 1$,
\[
 |Z_\varphi-H_\varphi|\lesssim(\rho^4+\rho^2|y-y_0|+|y-y_0|^2)\mathscr E_1(Z)^{1/2},  \qquad |Z_V-H_V|\lesssim(\rho^3+\rho|y-y_0|)\mathscr E_1(Z)^{1/2}.
\]
Integrating over $Q_\theta$ with $\theta\leqslant \theta_0$ proves \eqref{eqn:frozen-excess-decay}. It remains to show $\mathscr A_0H=0$. The same expansion gives $\|\gradG\mathbf S_r(Z-H)\|_{L^2(Q_1)}\lesssim r^2\mathscr E_1(Z)^{1/2}$. Since $\mathscr B_0$ depends only on derivatives, $D(\mathbf S_rH)=DH$ for the rescaling $\mathbf S_r$ centred at $y_0$, and the equation is invariant under $\mathbf S_r$, we find
\[
 |\mathscr B_0(H,Y)|=|\mathscr B_0(\mathbf S_r(H-Z),Y)|  \lesssim r^2\mathscr E_1(Z)^{1/2}\|\gradG Y\|_{L^2}
\]
for all small $r$ and all compactly supported admissible $Y$. Taking $r\to0$ completes the proof.
\end{proof}

\subsection{Polynomial corrections}
\label{sec:polynomial-corrections}
For $P\in\mathcal P$ (see \eqref{eqn:polynomial-norm}), the functional $\mathscr A_{p_0}P$ acts on compactly supported admissible test fields. Let $\mathscr A_{p_0}(\mathcal P)$ be the image of $\mathcal P$ under $\mathscr A_{p_0}$, and $\ker\mathscr A_{p_0}=\{P\in\mathcal P:\mathscr A_{p_0}P=0\}$. In Section~\ref{sec:linearised-polynomial-approximation}, the principal parts of the sources lie in the finite-dimensional space $\mathscr A_{p_0}(\mathcal P)$.
\begin{lemma}[Polynomial corrections]
\label{lem:polynomial-corrections}
There is a linear right inverse $\mathbf R_{p_0}:\mathscr A_{p_0}(\mathcal P)\to\mathcal P$ such that each $\mathbf R_{p_0}\mathscr F$ is of the form $(q(\xi),0)$ and
\begin{equation}
 \mathscr A_{p_0}\mathbf R_{p_0}\mathscr F=\mathscr F,\qquad  \|\mathbf R_{p_0}\mathscr F\|_1\leqslant C\|\mathscr F\|_{(W_G^{1,2}(Q_1(y_0)))^*},
 \qquad \mathscr F\in\mathscr A_{p_0}(\mathcal P),
 \label{eqn:right-inverse}
\end{equation}
with $C$ uniform in $p_0$.
\end{lemma}
\begin{proof}
By Lemma~\ref{lem:frozen-form-conjugacy}, it suffices to consider the square model. For $P=(a+B(y-y_0)+q(\xi),M\xi)$, the affine part and $M\xi$ contribute zero to $\mathscr B_0(P,Y)$ by \eqref{eqn:canonical-frozen-form} and integration by parts, so $\mathscr A_0P=\mathscr A_0(q,0)$. Thus $\mathscr A_0$ maps the finite-dimensional space of fields $(q(\xi),0)$ onto $\mathscr A_0(\mathcal P)$, and any linear right inverse of this map satisfies \eqref{eqn:right-inverse}, since all norms on a finite-dimensional space are equivalent.
\end{proof}

\section{Excess decay for the second variation}
\label{sec:linearised-polynomial-approximation}
In this section we prove a general  Campanato-type excess decay result for weak solutions $Z$ of the linearised equation
\[
 \mathscr B_U(Z,Y)=\mathscr F(Y),
\]
whose source term $\mathscr F$ with $\mathbf F \in (W_G^{1,2}(Q_r(y_0)))^*$ (see \eqref{eqn:functionals} for definition) has a principal part $\mathscr F_{\mathrm{pr}}\in\mathscr A_{p_0}(\mathcal P)$ with size $O(1)$ and a remainder of order $r^{N_G/2+\alpha}$ at scale $r$ (\cf Theorem~\ref{thm:campanato}). In Section~\ref{sec:tangential-regularity}, $Z$ will be taken as tangential derivatives $D_y^\beta U$, and $\mathscr{F}$ will be specified by \eqref{eqn:proper-source}, consisting of lower-order tangential derivatives of the transformed spinor $U$, and $\mathscr F_{\mathrm{pr}}$ will be lower-order derivatives frozen at $p_0$.

The idea is to treat $\mathscr B_U$ as a perturbation of the frozen form $\mathscr B_{p_0}$ of Section~\ref{sec:frozen-operator}. By Proposition~\ref{prop:legendre-initial-estimates}, the coefficients of $\mathscr B_U$ differ from the frozen ones by $O(r^\alpha)$ on $Q_r(y_0)$. If only the difference $\mathscr B_U-\mathscr B_{p_0}$ were present, the usual Campanato iteration would apply: at each scale, we approximate $Z$ by a polynomial $P_r$, compare $Z-P_r$ with a solution of the frozen equation, and update $P_r$ by the kernel element given by Proposition~\ref{prop:frozen-operator}(ii). However, the density \eqref{eqn:density-powers} \textit{depends on the original position $X=(\varphi,y)$}. The constant term $a$ of a kernel element \textit{does move the centre $p_0$ in the normal $z$-direction}, and thus contributes a principal part $\mathscr N_{p_0}a$ to the source (see \eqref{eqn:normal-variation} below), of order $|a|=O(1)$, which is not small at scale $r$. To remove this systematic error caused by freezing, together with the given principal part $\mathscr F_{\mathrm{pr}}$, we further correct the polynomials by Lemma~\ref{lem:polynomial-corrections}, so that they lie in the affine class \eqref{eqn:polynomial-comparison-class}. The resulting argument is a variant of the Campanato iteration.

Throughout this section, we fix a branch point $p_0=p(y_0)$. By Proposition~\ref{prop:legendre-initial-estimates}(iii) and $\varphi_{p_0}-g(y_0)=O(r^2)$ on $Q_r(y_0)$, for $0<r\leqslant r_0$, we have
\begin{equation}
 |\gradG(U-U_{p_0})|\leqslant Cr^\alpha,\qquad |X-p_0|\leqslant Cr^2
 \qquad\text{on }Q_r(y_0)\setminus\{\xi=0\}.
 \label{eqn:background-deviation}
\end{equation}
We write
\begin{equation}
 \mathbf N_{p_0}a=\frac{d}{dt}\Big|_{t=0}\mathbf A_{\mathrm{pr}}(p_0+t(a,0),DU_{p_0}),\qquad a\in\C,
 \label{eqn:normal-variation}
\end{equation}
that is, the change of the frozen coefficients when the centre moves in the normal direction by $a$ as discussed above, with $DU_{p_0}$ fixed. 

\subsection{Comparison with the frozen operator}
\label{sec:frozen-comparison}
The following lemma collects the estimates comparing $\mathscr B_U$ with
$\mathscr B_{p_0}$ that are used below.
\begin{lemma}[Comparison with the frozen operator]
\label{lem:linearised-operator-comparison}
For $0<r\leqslant r_0\ll 1$, the following hold on $Q_r=Q_r(y_0)$.
\begin{enumerate}[label=\textnormal{(\roman*)}]
\item For $Z\in W_G^{1,2}(Q_r)$ and $Y\in W_{G,0}^{1,2}(Q_r)$,
\begin{equation}
 |\mathscr B_U(Z,Y)-\mathscr B_{p_0}(Z,Y)| \leqslant Cr^\alpha(\|\gradG Z\|_{L^2}+r^{-2}\|Z_\varphi\|_{L^2})\|\gradG Y\|_{L^2}.
 \label{eqn:background-form-comparison}
\end{equation}
In particular, $\mathscr B_U$ is uniformly coercive on $W_{G,0}^{1,2}(Q_r)$.
\item For $P\in\mathcal P$ with constant term $a$, the error field
$\mathbf E=\mathbf A_UP-\mathbf A_{\mathrm{pr},p_0}P-\mathbf N_{p_0}a$ satisfies
\begin{align}
 r^{-N_G/2}\|\mathbf E\|_{(W_G^{1,2}(Q_r))^*}&\leqslant Cr^\alpha(|a|+\|\nabla P\|),
 \label{eqn:polynomial-error}\\
 \|\mathbf E\|_{C_G^{0,\alpha}(Q_r)}&\leqslant Cr^\alpha(|a|+\|\nabla P\|).
 \label{eqn:polynomial-annular-error}
\end{align}
\item For every $a\in\C$ there is $P_a\in\mathcal P$ with constant term $a$ such that $\mathscr N_{p_0}a=-\mathscr A_{p_0}P_a$. In particular,
\begin{equation}
 \mathscr N_{p_0}a\in\mathscr A_{p_0}(\mathcal P),\qquad
 \|\mathbf R_{p_0}\mathscr N_{p_0}a\|_1\leqslant C|a|.
 \label{eqn:normal-variation-image}
\end{equation}
\end{enumerate}
\end{lemma}

\begin{proof}
The main ingredients are the structure \eqref{eqn:density-powers} of the transformed density, the deviation estimate \eqref{eqn:background-deviation}, and the weighted inequalities in Lemma~\ref{lem:weighted-elementary-inequalities}.

(i) Write the second variation as
$\mathscr B_U(Z,Y)=\int(\mathbf a(DZ,DY)+\mathbf b(Z,DY)+\mathbf b(Y,DZ) +\mathbf c(Z,Y))$, with coefficients given by second derivatives of $\hat f$ at $U$. By \eqref{eqn:density-powers}, the principal part is $D^2_{DU}\hat f_{\mathrm{pr}}(X,\gradG U)$ up to $O(\rho^2)$, and
$\mathscr B_{p_0}$ is obtained by evaluating it at $(p_0,\gradG U_{p_0})$. By \eqref{eqn:background-deviation}, the difference of the principal parts is bounded by $Cr^\alpha|\gradG Z||\gradG Y|$. The remaining terms are bounded by
\[
 C(|Z_\varphi|+\rho|Z_V|+\rho^2|\gradG Z|)|\gradG Y|
 +C(|Y_\varphi|+\rho|Y_V|)|\gradG Z|
 +C(|Z_\varphi|+|Z_V|)(|Y_\varphi|+|Y_V|).
\]
By Lemma~\ref{lem:weighted-elementary-inequalities},
$$\|Z_V\|_{L^2}\leqslant r\|\gradG Z\|_{L^2}, \qquad \|Y_\varphi\|_{L^2}+r\|Y_V\|_{L^2}\leqslant Cr^2\|\gradG Y\|_{L^2}.$$
The Cauchy--Schwarz inequality then bounds the remaining terms by
\[
 C(\|Z_\varphi\|_{L^2}+r^2\|\gradG Z\|_{L^2})\|\gradG Y\|_{L^2}=Cr^2(\|\gradG Z\|_{L^2}+r^{-2}\|Z_\varphi\|_{L^2})\|\gradG Y\|_{L^2},
\]
which proves \eqref{eqn:background-form-comparison}. For compactly supported $Z$, we have $r^{-2}\|Z_\varphi\|_{L^2}\leqslant C\|\gradG Z\|_{L^2}$, so coercivity of $\mathscr B_U$ follows from the coercivity \eqref{eqn:linearised-coercivity} of the frozen form.

(ii) Write $P=(a,0)+(P-(a,0))$. Since $DU$ is unchanged by the constant
variation $(a,0)$, $\mathbf A_U(a,0)=\frac{d}{dt}\big|_{t=0}\mathbf A_{\mathrm{pr}}(X+t(a,0),DU)+O(\rho^2|a|)$, where the $O(\rho^2|a|)$-remainder is measured by the weighted $C^0$-norm \eqref{eqn:weighted-holder-field}. Compared with $\mathbf N_{p_0}a$, the coefficients are evaluated at $(X,DU)$ instead of $(p_0,DU_{p_0})$, so \eqref{eqn:background-deviation} gives an error $O(r^\alpha|a|)$. The non-constant part satisfies
$$|\gradG P|\leqslant C\|\nabla P\|,~~|P_\varphi-a|\leqslant Cr^2\|\nabla P\|,~~|P_V|\leqslant\rho\|\nabla P\|,\quad \textit{on } Q_r,$$ so the proof of (i), applied pointwise, implies
$$|\mathbf E^1_G|\leqslant Cr^\alpha(|a|+\|\nabla P\|), \qquad |\mathbf E^0_\varphi|+\rho^{-1}|\mathbf E^0_V|\leqslant C(|a|+\|\nabla P\|).$$
Then the dual error estimate \eqref{eqn:polynomial-error} follows from Lemma~\ref{lem:weighted-elementary-inequalities}(ii). Together with the transformed $C^{1,\alpha}$-estimates \eqref{eqn:initial-holder-estimate} of Proposition~\ref{prop:legendre-initial-estimates}(ii), the same computation shows \eqref{eqn:polynomial-annular-error}.

(iii) This follows essentially from the Leibniz rule. Let $p_t=p_0+t(a,0)$. Using the exponential gauge and parallel transport, we transport $\Phi_{p_0}$ to a smooth family of flat Dirac models $\Phi_t$ centred at $p_t$, each minimising the energy frozen at $p_t$. Their transformed models have the form
$$U_t=(g(y_0)+ta+B_t(y-y_0)+q_t(\xi),\,M_t\xi),\qquad U_0=U_{p_0},$$ and satisfy $$\int\mathbf A^1_{\mathrm{pr}}(p_t,DU_t):DY=0,$$ for every compactly supported admissible $Y$.  Differentiating at $t=0$ implies
$$\mathscr N_{p_0}a=-\mathscr A_{p_0}P_a,\qquad \textit{where }P_a:=\left.\frac{d}{dt}\right|_{t=0}U_t\in\mathcal P,$$ 
whose constant term is $a$. The bound \eqref{eqn:normal-variation-image} follows from the bound \eqref{eqn:right-inverse} of the right inverse in Lemma~\ref{lem:polynomial-corrections}.
\end{proof}

\begin{remark}
\label{rmk:translated-backgrounds}
For $h\in\R^{n-2}$ with $|h|\leqslant r^2/100$, write $U_h(\xi,y)=U(\xi,y+h)$ and $U_{h,t}=(1-t)U+tU_h$ for $t\in[0,1]$. The proof of Lemma~\ref{lem:linearised-operator-comparison}(i) applies verbatim to the following bilinear forms:
\begin{itemize}
\item \textnormal{(Spine-translated bilinear form)} The second variation at $U_h$ of the density translated by $h$ along the spine,
\begin{equation}
\begin{aligned}
     \mathscr B_h(Z,Y)=\frac{d}{ds}\Big|_{s=0}\int
 (&\mathbf A^1(\xi,y+h,U_h+sZ,DU_h+sDZ):DY
 \\&+\mathbf A^0(\xi,y+h,U_h+sZ,DU_h+sDZ)\cdot Y).
\end{aligned}
 \label{eqn:translated-form}
\end{equation}
\item \textnormal{(Averaged bilinear form)} The average along the segment from $U$
to $U_h$,
\begin{equation}
\begin{aligned}
     \overline{\mathscr B}_h(Z,Y)=\frac{d}{ds}\Big|_{s=0}\int_0^1\int
 (& \mathbf A^1(\xi,y+th,U_{h,t}+sZ,DU_{h,t}+sDZ):DY
 \\&\qquad+\mathbf A^0(\xi,y+th,U_{h,t}+sZ,DU_{h,t}+sDZ)\cdot Y)\,dt.
\end{aligned}
 \label{eqn:averaged-form}
\end{equation}
\end{itemize}
\end{remark}

As a consequence, we obtain the following weighted Caccioppoli inequality with source $\mathscr F$. 
\begin{corollary}
\label{cor:linearised-caccioppoli}
Let $\mathscr B$ be $\mathscr B_U$ or one of the forms in Remark~\ref{rmk:translated-backgrounds}. If $Z\in W_G^{1,2}(Q_r(y_0))$ satisfies $\mathscr B(Z,Y)=\mathscr F(Y)$ for every compactly supported admissible $Y$, then
\begin{equation}
 \|\gradG Z\|_{L^2(Q_{r/2}(y_0))}
 \leqslant C(r^{-2}\|Z\|_{L^2_G(Q_r(y_0))}
 +\|\mathbf F\|_{(W_G^{1,2}(Q_r(y_0)))^*}).
 \label{eqn:linearised-caccioppoli}
\end{equation}
\end{corollary}

\begin{proof}
Let $\eta=\eta_1(\rho)\eta_2(y)$ be a cutoff of $Q_{r/2}\subset Q_r$ with $\eta_1=1$ for $\rho\leqslant r/2$. Since $\mathscr B$ is symmetric, $\mathscr B(\eta Z,\eta Z)=\mathscr B(Z,\eta^2Z)+\int\mathbf a(Z\otimes D\eta,Z\otimes D\eta)$, as in Step~1 of the proof of Proposition~\ref{prop:frozen-operator}(ii) without source term. Since $D_\xi\eta$ is supported where $\rho\geqslant r/2$, the weighted bounds on $\mathbf a$ give $|\mathbf a(Z\otimes D\eta,Z\otimes D\eta)|\leqslant Cr^{-4}(|Z_\varphi|^2+\rho^2|Z_V|^2)$. The conclusion follows from the coercivity in Lemma~\ref{lem:linearised-operator-comparison}(i) and $|\mathscr F(\eta^2Z)|\leqslant\|\mathbf F\|_{(W_G^{1,2})^*}\|\gradG(\eta^2Z)\|_{L^2}$.
\end{proof}

\subsection{Campanato iteration}
\label{sec:linearised-campanato}
By Lemma~\ref{lem:polynomial-corrections} and Lemma~\ref{lem:linearised-operator-comparison}(iii), for $\mathscr F\in\mathscr A_{p_0}(\mathcal P)$, the polynomials solving the frozen equation with source $\mathscr F$, corrected for the motion of the centre, form an affine space
\begin{equation}
 \mathcal P(\mathscr F) 
 =\{P\in\mathcal P:\mathscr A_{p_0}P+\mathscr N_{p_0}a_P=\mathscr F\}
 =\mathbf R_{p_0}\mathscr F + \{H-\mathbf R_{p_0}\mathscr N_{p_0}a_H:H\in\ker\mathscr A_{p_0}\}.
 \label{eqn:polynomial-comparison-class}
\end{equation}
\begin{theorem}[Campanato estimate]
\label{thm:campanato}
Let $\mathscr F_{\mathrm{pr}}\in\mathscr A_{p_0}(\mathcal P)$, and suppose that $Z\in W_G^{1,2}(Q_{r_0}(y_0))$ satisfies $\mathscr B_U(Z,Y)=\mathscr F(Y)$ for every compactly supported admissible $Y$, where
\begin{equation}
 \|\mathbf F-\mathbf F_{\mathrm{pr}}\|_{(W_G^{1,2}(Q_r(y_0)))^*} \leqslant Cr^{N_G/2+\alpha},\qquad 0<r\leqslant r_0.
 \label{eqn:source-hyp}
\end{equation}
Then there exists a unique comparison polynomial $P_{Z,y_0}\in \mathcal P(\mathscr F_{\mathrm{pr}})$ of $Z$ at $(0,y_0)$ with estimates,
\begin{equation}
 \mathscr E_r(Z-P_{Z,y_0})\leqslant Cr^{2\alpha},\qquad 0<r\leqslant r_0,
 \label{eqn:campanato-decay}
\end{equation}
where the constant depends linearly on the data at scale $r_0$.
\end{theorem}

\begin{proof}
For $0<r\leqslant r_0$, let $P_r$ minimise $\mathscr E_r(Z-P)$ over $\mathcal P(\mathscr F_{\mathrm{pr}})$, and set $e(r)=\mathscr E_r(Z-P_r)^{1/2}$: this  minimiser exists and is unique, since $\mathscr E_r(Z-\cdot)$ is a coercive strictly convex quadratic functional (\cf \eqref{eqn:polynomial-norm-equivalence}) on the finite-dimensional affine space $\mathcal P(\mathscr F_{\mathrm{pr}})$. As $\mathscr A_{p_0}P_r+\mathscr N_{p_0}a_{P_r}=\mathscr F_{\mathrm{pr}}$, the field $W=Z-P_r$ then satisfies
\begin{equation}
 \mathscr B_U(W,Y)=(\mathscr F-\mathscr F_{\mathrm{pr}})(Y) -\int(\mathbf E^1:DY+\mathbf E^0\cdot Y),
 \label{eqn:campanato-equation}
\end{equation}
where $\mathbf E$ is the error field of $P_r$. By \eqref{eqn:source-hyp} and the error estimate \eqref{eqn:polynomial-error} of Lemma~\ref{lem:linearised-operator-comparison}(ii), the right-hand side has $(W_G^{1,2}(Q_r(y_0)))^*$-norm $O(r^{N_G/2+\alpha}(1+\|P_r\|_{r_0}))$.

Let now $W_0$ solve $\mathscr B_{p_0}(W_0,Y)=0$ on $Q_{r/2}(y_0)$ with $W-W_0\in W_{G,0}^{1,2}(Q_{r/2}(y_0))$, whose existence is guaranteed by the coercivity estimate \eqref{eqn:linearised-coercivity} of the frozen form and the Lax--Milgram theorem. Testing \eqref{eqn:campanato-equation} and the frozen equation with $W-W_0$ and applying weighted Caccioppoli estimates (\cf Corollary~\ref{cor:linearised-caccioppoli}), we find
\[
\begin{aligned}
 \|\gradG(W-W_0)\|_{L^2}^2
 &\lesssim (\mathscr B_{p_0}-\mathscr B_U)(W,W-W_0)
\\& \qquad +(\mathscr F-\mathscr F_{\mathrm{pr}})(W-W_0)
 -\int(\mathbf E^1:D(W-W_0)+\mathbf E^0\cdot(W-W_0))\\
 &\lesssim(r^\alpha(\|\gradG W\|_{L^2(Q_{r/2})}
 +r^{-2}\|W_\varphi\|_{L^2(Q_{r/2})}) \\
 & \qquad +r^{N_G/2+\alpha}(1+\|P_r\|_{r_0}))\|\gradG(W-W_0)\|_{L^2}\\
 &\lesssim r^{N_G/2+\alpha}(e(r)+1+\|P_r\|_{r_0})\|\gradG(W-W_0)\|_{L^2}.
\end{aligned}
\]
By weighted Poincar\'e inequalities, this implies
\begin{equation}
 \mathscr E_{r/2}(W-W_0)^{1/2}\leqslant Cr^{-N_G/2}\|\gradG(W-W_0)\|_{L^2(Q_{r/2})} \leqslant Cr^\alpha(e(r)+1+\|P_r\|_{r_0}).
 \label{eqn:campanato-comparison}
\end{equation}

By Proposition~\ref{prop:frozen-operator}(ii), applied to scale $r/2$, there is $H\in\ker\mathscr A_{p_0}$ with $\|H\|_{r/2}\leqslant C\mathscr E_{r/2}(W_0)^{1/2}$ and $\mathscr E_{\theta r}(W_0-H)\leqslant C\theta^4\mathscr E_{r/2}(W_0)$ for $\theta<\theta_0/2\ll 1$. As $H$ may have a non-zero constant term, we consider the competitor $P_r+H-\mathbf R_{p_0}\mathscr N_{p_0}a_H\in\mathcal P(\mathscr F_{\mathrm{pr}})$. Since the correction $\mathbf R_{p_0}\mathscr N_{p_0}a_H$ is purely quadratic in $\xi$, we know $\mathscr E_{\theta r}(\mathbf R_{p_0}\mathscr N_{p_0}a_H)^{1/2}\lesssim \|\mathbf R_{p_0}\mathscr N_{p_0}a_H\|_1\lesssim |a_H|\lesssim r^2\|H\|_{r/2}$ (\cf Lemma \ref{lem:linearised-operator-comparison} (iii)). Together with
\eqref{eqn:campanato-comparison}, this gives
\[
 e(\theta r)\leqslant C_0\theta^2e(r) + C_\theta r^\alpha(e(r)+1+\|P_r\|_{r_0}),
\]
with $C_0$ independent of $\theta$. Fixing $\theta$ with $C_0\theta^2\leqslant\theta^\alpha/4$, and then $r_0$ with $C_\theta r_0^\alpha\leqslant\theta^\alpha/4$, we obtain
\begin{align}
 e(\theta r)&\leqslant\tfrac{1}{2}\theta^\alpha e(r)+Cr^\alpha(1+\|P_r\|_{r_0}),
 \label{eqn:campanato-one-step-decay}
\end{align}
and
\begin{equation}
\label{eqn:polynomial-coefficient-increments}
\|P_{\theta r}-P_r\|_r\lesssim \mathscr E_{\theta r}(P_{\theta r}-P_r)^{1/2}\leqslant e(\theta r)+\mathscr E_{\theta r}(Z-P_r)^{1/2} \leqslant C(e(r)+r^\alpha(1+\|P_r\|_{r_0})).
\end{equation}

Theorem \ref{thm:campanato} follows from iterating \eqref{eqn:campanato-one-step-decay} and \eqref{eqn:polynomial-coefficient-increments}
along $r_k=\theta^kr_0$ in the standard way. \end{proof}

\section{Tangential smoothness and real analyticity}
\label{sec:tangential-regularity}
We now prove Theorem~\ref{thm:main-higher-regularity} by induction on the order of tangential derivatives. Suppose that $D_y^\gamma U$ has been controlled for all $|\gamma|<k$. Then at order $k$, differentiating the transformed weak equation \eqref{eqn:full-weak-equation} shows that $Z_\beta=D_y^\beta U$ with $|\beta|=k$ solves the linearised equation $\mathscr B_U(Z_\beta,Y)=\mathscr F_\beta(Y)$, whose source \textit{depends only on lower-order derivatives}. To apply Theorem~\ref{thm:campanato}, we need two properties of $\mathscr F_\beta$:
\begin{itemize}
\item \emph{Principal part.} $\mathscr F_\beta$ is close to its principal part $\mathscr F_{\beta,\mathrm{pr},y}$ at every centre $(0,y)$. This requires pointwise and $L^\infty$/H\"older bounds on the lower-order derivatives, not only $L^2$-bounds. We therefore upgrade the $L^2$-excess decay of Theorem~\ref{thm:campanato} to weighted $C^{1,\alpha}$-estimates (Proposition~\ref{prop:weighted-holder-estimate}), using interior Schauder estimates at the scale of the distance to $\{\xi=0\}$.
\item \emph{Compatibility and polynomial corrections.} The principal part lies in $\mathscr A_{p(y)}(\mathcal P)$.
This follows by taking $\beta$-th derivative of the identity expressing that transformed models are critical points of the frozen energy, along a family of models following the branching locus (Proposition~\ref{prop:differentiated-equation}(ii)).
\end{itemize}
The smoothness of the tangential derivatives of $U$ (and hence of the branching locus) shall follow from above discussion. To obtain real-analyticity, we use Kato's trick \cite{Kat96} together with the diffeomorphism invariance of the transformed system and a power-series majorant to obtain factorial bounds on the tangential derivatives (\cf Theorem~\ref{thm:tangential-regularity}).

\subsection{From \texorpdfstring{$L^2$}{L2} excess decay to weighted
\texorpdfstring{$C^{1,\alpha}$}{C1,alpha}-estimates}
\label{sec:weighted-holder-estimates}
Let $\Omega=Q_r(y_0)$ with $0<r\leqslant r_0$, and let $Z$ be a field on $\Omega$ which has a comparison polynomial
\begin{equation}
 P_{Z,y}=(a_Z(y)+B_Z(y)(y'-y)+q_Z(y,\xi),\ M_Z(y)\xi)
 \label{eqn:field-leading-polynomial}
\end{equation}
at every centre $(0,y)\in\Omega$. We define the following scaling-invariant $\mathcal{X}_\alpha$-norm for fields $Z$,
\begin{equation}
\begin{aligned}
 \|Z\|_{\mathcal X_\alpha(\Omega)}
 &:=r^{-N_G/2-2}\|Z\|_{L^2_G(\Omega)}
   +r^{-N_G/2}\|\gradG Z\|_{L^2(\Omega)}\\&\quad
   +r^{-2}\|a_Z\|_{C^0}
   +\|D_ya_Z\|_{C^0}
   +\|q_Z\|_{C^0}
   +\|M_Z\|_{C^0}\\
 &\quad+r^\alpha(
   [D_ya_Z]_{C^{0,\alpha/2}}
   +[q_Z]_{C^{0,\alpha/2}}
   +[M_Z]_{C^{0,\alpha/2}})+\sup_{y,\,s}
   \left(\frac{r}{s}\right)^\alpha
   \mathscr E_s(Z-P_{Z,y})^{1/2}\\
 &\quad+\sup_{\xi,y}\left(\frac{r}{\rho}\right)^\alpha
   \|(\rho^{-2}(Z-P_{Z,y})_\varphi,\  \rho^{-1}(Z-P_{Z,y})_V)\|_{C^{1,\alpha}(Q_{\xi,y},g_\rho)}.
\end{aligned}
\label{eqn:field-norm}
\end{equation}
Here the norms of coefficients are taken over
$B_{r^2}^{n-2}(y_0)$; the first supremum is over centres $(0,y)\in\Omega$ and radii $0<s\leqslant r$ with $Q_s(y)\subset\Omega$; the second supremum is over $Q_{\xi,y}\subset\Omega$ with $\rho=|\xi|>0$. As $B_Z=D_ya_Z$ (established below), the coefficient bounds and the last term in \eqref{eqn:field-norm} imply
\[
 \sup_{\Omega\setminus\{\xi=0\}} (r^{-2}|Z_\varphi|+\rho^{-1}|Z_V|+|\gradG Z|) \leqslant C\|Z\|_{\mathcal X_\alpha(\Omega)}.
\]

Correspondingly, for a source $\mathbf F$ with prescribed principal coefficient fields $\mathbf F_{\mathrm{pr},y}$ at the centres $(0,y)$, we define the scaling-invariant $\mathcal{Y}_\alpha$-norm
\begin{equation}
\begin{aligned}
 \|\mathbf F\|_{\mathcal Y_\alpha(\Omega)}
 &=r^{-N_G/2}\|\mathbf F\|_{(W_G^{1,2}(\Omega))^*}
   +\sup_y \|\mathbf F_{\mathrm{pr},y}\|_{C_G^{0,\alpha}(Q_1(y))}\\ 
   &\quad+\sup_{y,\,s}\left(\frac{r}{s}\right)^\alpha s^{-N_G/2} \|\mathbf F-\mathbf F_{\mathrm{pr},y}\|_{(W_G^{1,2}(Q_s(y)))^*} +\sup_{\xi,y} \left(\frac{r}{\rho}\right)^\alpha\|\mathbf F-\mathbf F_{\mathrm{pr},y}\|_{C_G^{0,\alpha}(Q_{\xi,y})},
\end{aligned}
\label{eqn:source-norm}
\end{equation}
with the suprema as above, $Q_s(y),Q_{\xi,y}\subset\Omega$.

We remark that the above norms just collect all the norms involved (adapted weighted $W^{1,2}$ and $C^{1,\alpha}$ norms) in the Campanato iteration procedure.  
\begin{proposition}[Weighted $C^{1,\alpha}$-estimates]
\label{prop:weighted-holder-estimate}
Let $0<2r\leqslant r_0$. Suppose that $Z\in W_G^{1,2}(Q_{2r}(y_0))$ satisfies $\mathscr B_U(Z,Y)=\mathscr F(Y)$ for every compactly supported
admissible $Y$, where $\|\mathbf F\|_{\mathcal Y_\alpha(Q_{2r}(y_0))}<\infty$ and $\mathscr F_{\mathrm{pr},y}\in\mathscr A_{p(y)}(\mathcal P)$ at every centre.
Then $Z$ has a comparison polynomial $P_{Z,y}\in\mathcal P(\mathscr F_{\mathrm{pr},y})$ at every centre $(0,y)\in Q_r(y_0)$, with $B_Z=D_ya_Z$, and
\begin{equation}
 \|Z\|_{\mathcal X_\alpha(Q_r(y_0))}\leqslant C(r^{-N_G/2-2}\|Z\|_{L^2_G(Q_{2r}(y_0))} +\|\mathbf F\|_{\mathcal Y_\alpha(Q_{2r}(y_0))} ).
 \label{eqn:local-linear-estimate}
\end{equation}
In particular,
\[
 |D_\xi Z_\varphi(\xi,y)| \leqslant C|\xi|( r^{-N_G/2-2}\|Z\|_{L^2_G(Q_{2r}(y_0))} +\|\mathbf F\|_{\mathcal Y_\alpha(Q_{2r}(y_0))} )
\]
on $Q_r(y_0)\setminus\{\xi=0\}$.

If further $Z\in W_{G,0}^{1,2}(Q_{2r}(y_0))$ is supported in $K\Subset Q_{2r}(y_0)$, then
\begin{equation}
 \|Z\|_{\mathcal X_\alpha(Q_{2r}(y_0))} \leqslant C_K \|\mathbf F\|_{\mathcal Y_\alpha(Q_{2r}(y_0))}.
 \label{eqn:compact-support-estimate}
\end{equation}
\end{proposition}

\begin{proof}
Since the equation is linear, it suffices to prove \eqref{eqn:local-linear-estimate} under the normalisation $r^{-N_G/2-2}\|Z\|_{L^2_G(Q_{2r}(y_0))}+\|\mathbf F\|_{\mathcal Y_\alpha(Q_{2r}(y_0))} \leqslant1$. Throughout the proof, $C$ is independent of $r$.

\textbf{Step 1: Excess decay at each centre.}
Applying the weighted Caccioppoli estimate (Corollary~\ref{cor:linearised-caccioppoli}, with the fixed interior radius ratio $3/4$) to $Q_{2r}(y_0)$, we find
\[
 r^{-N_G/2}\|\gradG Z\|_{L^2(Q_{3r/2}(y_0))} \leqslant C( r^{-N_G/2-2}\|Z\|_{L^2_G(Q_{2r}(y_0))} +r^{-N_G/2}\|\mathbf F\|_{(W_G^{1,2}(Q_{2r}(y_0)))^*} ) \leqslant C.
\]
Together with the normalisation, this bounds the first two terms of \eqref{eqn:field-norm}. For $(0,y)\in Q_r(y_0)$, we have $Q_r(y)\subset Q_{3r/2}(y_0)$, so \eqref{eqn:hardy-trace} implies
\[
 r^{-N_G/2-2}\|Z_\varphi\|_{L^2(Q_r(y))} +r^{-N_G/2-1}\|Z_V\|_{L^2(Q_r(y))}  \leqslant C.
\]
Thus Theorem~\ref{thm:campanato} applies at scale $r$ centred at $(0,y)$, and we obtain a unique comparison polynomial $P_{Z,y}\in\mathcal P(\mathscr F_{\mathrm{pr},y})$ with estimates
\[
 \left(\frac{r}{s}\right)^\alpha \mathscr E_s(Z-P_{Z,y})^{1/2}  +r^{-2}|a_Z(y)|+\|\nabla P_{Z,y}\| \leqslant C, \qquad 0<s\leqslant r.
\]

\medskip
\textbf{Step 2: Estimates away from $\{\xi=0\}$.}
Since $P_{Z,y}\in\mathcal P(\mathscr F_{\mathrm{pr},y})$, the difference $Z-P_{Z,y}$ satisfies
\[
 \mathscr B_U(Z-P_{Z,y},Y)=(\mathscr F-\mathscr F_{\mathrm{pr},y})(Y)
 -\int(\mathbf E^1:DY+\mathbf E^0\cdot Y),
\]
where $\mathbf E$ is the error field of $P_{Z,y}$ at $p(y)$.
Let $(\xi,y)\in Q_r(y_0)$ with $\rho=|\xi|>0$. With respect to the rescaled metric $g_\rho$ of \eqref{eqn:weighted-holder-field}, this is a uniformly elliptic system in divergence form on the annulus $Q_{\xi,y}$ with coefficients uniformly bounded in $C^{0,\alpha}$, by Proposition~\ref{prop:legendre-initial-estimates}(ii) and Lemma~\ref{lem:transformed-equation}.

The source has $C^{0,\alpha}$-norm at most $C(\rho/r)^\alpha$, by \eqref{eqn:source-norm}, the polynomial coefficient bounds from Step~1, and \eqref{eqn:polynomial-annular-error}.
Moreover, Step~1 yields an $L^2$-bound $C(\rho/r)^\alpha$ for the normalised field $
 (\rho^{-2}(Z-P_{Z,y})_\varphi,\ \rho^{-1}(Z-P_{Z,y})_V)$
with respect to the rescaled metric $g_\rho=\rho^{-2}g_\xi+\rho^{-4}g_y$. The interior elliptic estimate therefore implies
\[
 \|( \rho^{-2}(Z-P_{Z,y})_\varphi,\ \rho^{-1}(Z-P_{Z,y})_V  )\|_{C^{1,\alpha}(Q_{\xi,y},g_\rho)}  \leqslant C\left(\frac{\rho}{r}\right)^\alpha,
\]
which bounds the last term of the $\mathcal{X}_\alpha$-norm \eqref{eqn:field-norm}.

\medskip
\textbf{Step 3: Estimates of the comparison polynomials.}
Let $(0,y_1),(0,y_2)\in Q_r(y_0)$ and set
$s=4|y_1-y_2|^{1/2}\leqslant r$. Then
$Q_{s/2}(y_1)\subset Q_s(y_1)\cap Q_s(y_2)$.
Step~1 and the triangle inequality imply
\[
 \mathscr E_{s/2}(P_{Z,y_1}-P_{Z,y_2})^{1/2}
 \leqslant C\left(\frac{s}{r}\right)^\alpha,
\]
with $\mathscr E_{s/2}$ centred at $y_1$.
Expanding $P_{Z,y_2}$ around $y_1$ and using
\eqref{eqn:polynomial-norm-equivalence}, we obtain
\begin{equation}
\begin{aligned}
 |B_Z(y_2)-B_Z(y_1)|
 +|q_Z(y_2)-q_Z(y_1)|
 +|M_Z(y_2)-M_Z(y_1)|
 &\leqslant Cr^{-\alpha}|y_2-y_1|^{\alpha/2},\\
 |a_Z(y_2)-a_Z(y_1)-B_Z(y_1)(y_2-y_1)|
 &\leqslant Cr^{-\alpha}|y_2-y_1|^{1+\alpha/2}.
\end{aligned}
\label{eqn:coefficient-coherence}
\end{equation}
Hence $D_ya_Z=B_Z$. For more distant pairs of points, the coefficient supremum bounds in Step~1 imply the required H\"older bounds.
Consequently,
\[
\begin{aligned}
 &r^{-2}\|a_Z\|_{C^0}
 +\|D_ya_Z\|_{C^0}
 +\|q_Z\|_{C^0}
 +\|M_Z\|_{C^0}
 +r^\alpha(
 [D_ya_Z]_{C^{0,\alpha/2}}
 +[q_Z]_{C^{0,\alpha/2}}
 +[M_Z]_{C^{0,\alpha/2}})
 \leqslant C.
\end{aligned}
\]
Steps~1--3 prove \eqref{eqn:local-linear-estimate}.

\medskip
\textbf{Step 4: Compactly supported fields.}
Suppose that $Z\in W_{G,0}^{1,2}(Q_{2r}(y_0))$ is supported in $K\Subset Q_{2r}(y_0)$. Testing the equation with $Z$, and using coercivity and the weighted Poincar\'e inequalities, we find
\begin{align*}
    r^{-N_G/2-2} \|Z\|_{L^2_G(Q_{2r}(y_0))}+r^{-N_G/2}\|\gradG Z\|_{L^2(Q_{2r}(y_0))}
    & \leqslant Cr^{-N_G/2}  \|\mathbf F\|_{(W_G^{1,2}(Q_{2r}(y_0)))^*} \\
    & \leqslant C  \|\mathbf F\|_{\mathcal Y_\alpha(Q_{2r}(y_0))}.
\end{align*}
Choose a fixed interior neighbourhood of $K$ and cover it by finitely many spine-centred cylinders and annuli $Q_{\xi,y}$. Apply \eqref{eqn:local-linear-estimate} on the spine-centred cylinders and interior Schauder estimates on the annuli. Outside this neighbourhood, the comparison polynomials and the excesses at sufficiently small scales vanish. Patching these estimates together yields \eqref{eqn:compact-support-estimate}.
\end{proof}
\subsection{The differentiated equations}
\label{sec:tangential-source-estimates}
For a multi-index $\beta=(\beta_1,\ldots,\beta_{n-2})$, write $|\beta|=\sum_i\beta_i$, $\beta!=\prod_i\beta_i!$ and $Z_\beta=D_y^\beta U$, and let $e_j$ be the $j$-th unit multi-index. Let $D^q\mathbf A$ be the $q$-th derivative of $\mathbf A$ with respect to $(y,U,DU)$, with $\xi$ fixed, and set
\begin{equation}
 \Xi_\gamma=D_y^\gamma(y,U,DU)=
 \begin{cases}
 (e_j,Z_{e_j},DZ_{e_j}),&\gamma=e_j,\\ (0,Z_\gamma,DZ_\gamma),&|\gamma|\geqslant2.
 \end{cases}
 \label{eqn:chain-rule-input}
\end{equation}
Formally differentiating the transformed weak equation \eqref{eqn:full-weak-equation} $k=|\beta|$ times by the chain rule, we find that $Z_\beta$ satisfies $\mathscr B_U(Z_\beta,Y)=\mathscr F_\beta(Y)$, with $\mathbf F_{e_j}=-\partial_{y_j}\mathbf A$ and
\begin{equation}
 \mathbf F_\beta=-\sum_{q=2}^k\frac{1}{q!}\sum_{\substack{\gamma_1+\cdots+\gamma_q=\beta\\\gamma_i\neq0}}\frac{\beta!}{\gamma_1!\cdots\gamma_q!}\,D^q\mathbf A(\Xi_{\gamma_1},\ldots,\Xi_{\gamma_q}), \qquad k\geqslant2.
 \label{eqn:proper-source}
\end{equation}
The source involves only $Z_\gamma$ with $|\gamma|\leqslant k-1$. At a centre $(0,y)$, the \emph{principal part} of the source is
\begin{equation}
 \mathbf F_{\beta,\mathrm{pr},y}=-\sum_{q=2}^k\frac{1}{q!}\sum_{\substack{\gamma_1+\cdots+\gamma_q=\beta\\\gamma_i\neq0}}\frac{\beta!}{\gamma_1!\cdots\gamma_q!}\,D^q\mathbf A_{\mathrm{pr}}(p(y),DU_{p(y)}) (\Xi^{\mathrm{pr}}_{\gamma_1},\ldots,\Xi^{\mathrm{pr}}_{\gamma_q}),\qquad k\geqslant 2,
 \label{eqn:leading-source}
\end{equation}
where $\Xi_\gamma^{\mathrm{pr}} =(D_y^\gamma y,\ a_{Z_\gamma}(y),\ DP_{Z_\gamma,y})$, in the order of $(y,U,DU)$ as in \eqref{eqn:chain-rule-input} (the $U$-component is $a_{Z_\gamma}(y)$, since $\mathbf A_{\mathrm{pr}}(X,DU)$ depends on $U$ only through $X=(\varphi,y)$), and $\mathbf F_{e_j,\mathrm{pr},y}=-\partial_{y_j}\mathbf A_{\mathrm{pr}}(p(y),DU_{p(y)})$ for $k=1$. On $Q_{\xi,y}$, in the normalised quantities of the weighted H\"older norms \eqref{eqn:weighted-holder-field}, the independence of $\mathbf A_{\mathrm{pr}}$ from the value of $V$ gives
\[
\begin{aligned}
 \mathbf F_\beta-\mathbf F_{\beta,\mathrm{pr},y}
 &=O(\rho^2)+O(|X-p(y)|+|\gradG(U-U_{p(y)})|)\\
 &\quad+O\left(\sum_{1\leqslant|\gamma|<k}(|(Z_\gamma)_\varphi-a_{Z_\gamma}(y)|+|\gradG(Z_\gamma-P_{Z_\gamma,y})|)\right)=O(\rho^\alpha),
\end{aligned}
\]
by the density decomposition \eqref{eqn:density-powers}, the initial estimates of Proposition~\ref{prop:legendre-initial-estimates}(ii), and the lower-order $\mathcal X_\alpha$-estimates.

The following proposition collects the properties of the sources used in the induction.
\begin{proposition}[Source estimates]
\label{prop:differentiated-equation}
Let $\Omega=Q_r(y_0)$, $k\geqslant1$, and suppose that $\|Z_\gamma\|_{\mathcal X_\alpha(\Omega)}<\infty$ for $|\gamma|\leqslant k-1$,
with
\begin{equation}
 a_{Z_\gamma}=D_y^\gamma g,\qquad q_{Z_\gamma}=D_y^\gamma q,\qquad
 M_{Z_\gamma}=D_y^\gamma M .
 \label{eqn:coefficient-identification}
\end{equation}
\begin{enumerate}[label=\textnormal{(\roman*)}]
\item \textnormal{(Source estimate)} For $k=1$, the source $\mathbf F_{e_j}=-\partial_{y_j}\mathbf A$ satisfies $\|\mathbf F_{e_j}\|_{\mathcal Y_\alpha(\Omega)}\leqslant C_1$ by the background estimates. For $k\geqslant2$, $q\geqslant1$ and nonzero $\gamma_1,\ldots,\gamma_q$ with $|\gamma_i|\leqslant k-1$, the term $D^q\mathbf A(\Xi_{\gamma_1},\ldots,\Xi_{\gamma_q})$ with the principal parts defined above satisfies
\begin{equation}
 \|D^q\mathbf A(\Xi_{\gamma_1},\ldots,\Xi_{\gamma_q})\|_{\mathcal Y_\alpha(\Omega)} \leqslant C_q\prod_{i=1}^q(\|Z_{\gamma_i}\|_{\mathcal X_\alpha(\Omega)}
 +\mathbf 1_{\{|\gamma_i|=1\}}), \label{eqn:estimate-on-dq}
\end{equation}
where $C_q\leqslant \Lambda_0\Lambda_1^qq!$ for real-analytic ambient data, with $\Lambda_0,\Lambda_1$ independent of $q$ and the multi-indices. In particular, $\|\mathbf F_\beta\|_{\mathcal Y_\alpha(\Omega)}<\infty$ for $|\beta|=k$.
\item \textnormal{(Validity of polynomial corrections)} For $|\beta|=k$ and every centre $(0,y)$, $\mathscr F_{\beta,\mathrm{pr},y}\in\mathscr A_{p(y)}(\mathcal P)$.
\end{enumerate}
\end{proposition}
We remark that the indicator function in \eqref{eqn:estimate-on-dq} accounts for the entry $e_j$ of $\Xi_{e_j}$.

\begin{proof}
\textbf{Step 1: Pointwise structure and the constants.}
By \eqref{eqn:density-powers}, $D^q\mathbf A(\Xi_{\gamma_1},\ldots,\Xi_{\gamma_q})$ is a multilinear expression in the normalised quantities $(Z_{\gamma_i})_\varphi$, $\gradG Z_{\gamma_i}$, $\rho^{-1}(Z_{\gamma_i})_V$ (and $e_j$), whose coefficients are derivatives of order at most $q+1$ of the energy density $\hat f$, regarded as a function of $\rho$ and the normalised data $(X,\widehat u/\rho,\gradG U)$ (see Lemma~\ref{lem:transformed-equation}). As in the proof of Lemma~\ref{lem:transformed-equation}(ii), the resulting coefficient field satisfies the weighted bounds \eqref{eqn:weighted-coefficient-bounds} with constant $C_q\prod_i(\|Z_{\gamma_i}\|_{\mathcal X_\alpha}+\mathbf 1_{\{|\gamma_i|=1\}})$, where
\begin{equation}
    C_q=C^{q+1}\left(1+\max_{m\leqslant q+3}\sup\|D^m\hat f\|\right)
    \label{eqn:source-constants}
\end{equation}
with $C$ independent of $q$, where $D^m\hat f$ denotes the $m$-th derivative in $(X,\widehat u/\rho,\gradG U)$, and the supremum is over $\Omega\setminus\{\xi=0\}$ along $U$. For real-analytic data, $\hat f$ is real-analytic, so $\sup\|D^m\hat f\|\leqslant\Lambda\Lambda_1^mm!$ for some constants $\Lambda,\Lambda_1>0$. Since $(q+3)!\leqslant4^{q+3}q!$, we obtain $C_q\leqslant\Lambda_0\Lambda_1^qq!$ for sufficiently large $\Lambda_0,\Lambda_1\gg 1$.

\medskip
\textbf{Step 2: Closeness to the principal part.}
Fix a centre $(0,y)$, and let $D^q\mathbf A_{\mathrm{pr}}$ denote the derivative of $\mathbf A_{\mathrm{pr}}$ evaluated at $(p(y),DU_{p(y)})$. The difference of the two multilinear expressions telescopes:
\begin{align}\label{eqn:source-telescoping}
\begin{split}
     D^q\mathbf A(\Xi_{\gamma_1},\ldots,\Xi_{\gamma_q}) - & D^q\mathbf A_{\mathrm{pr}}(\Xi^{\mathrm{pr}}_{\gamma_1},\ldots,\Xi^{\mathrm{pr}}_{\gamma_q}) \\ 
     &=(D^q\mathbf A-D^q\mathbf A_{\mathrm{pr}})(\Xi_{\gamma_1},\ldots,\Xi_{\gamma_q})\\ 
     &\quad+\sum_{i=1}^qD^q\mathbf A_{\mathrm{pr}}(\Xi^{\mathrm{pr}}_{\gamma_1},\ldots, \Xi^{\mathrm{pr}}_{\gamma_{i-1}},\Xi_{\gamma_i}-\Xi^{\mathrm{pr}}_{\gamma_i}, \Xi_{\gamma_{i+1}},\ldots,\Xi_{\gamma_q}).
\end{split} 
\end{align}
On $Q_{\xi,y}$, with the normalised quantities \eqref{eqn:weighted-holder-field} and the $C^{0,\alpha}$-norms with respect to $g_\rho$, we have
\[
 \|\Xi_{\gamma_i}\|+\|\Xi^{\mathrm{pr}}_{\gamma_i}\| \leqslant C(\|Z_{\gamma_i}\|_{\mathcal X_\alpha}+\mathbf 1_{\{|\gamma_i|=1\}})
\]
and
\[
 \|(Z_{\gamma_i})_\varphi-a_{Z_{\gamma_i}}(y)\|_{C^{0,\alpha}} +\|\gradG(Z_{\gamma_i}-P_{Z_{\gamma_i},y})\|_{C^{0,\alpha}}
 \leqslant C\left(\frac{\rho}{r}\right)^\alpha\|Z_{\gamma_i}\|_{\mathcal X_\alpha},
\]
by the last term in \eqref{eqn:field-norm} and
$|(P_{Z_{\gamma_i},y})_\varphi-a_{Z_{\gamma_i}}(y)| \leqslant C\rho^2\|\nabla P_{Z_{\gamma_i},y}\|$.
 Since $\mathbf A_{\mathrm{pr}}$ is independent of the value of $V$, the $V$-value component of $\Xi_{\gamma_i}-\Xi^{\mathrm{pr}}_{\gamma_i}$ contributes zero in each summand of \eqref{eqn:source-telescoping}. Moreover, in the corresponding normalised coefficient blocks,
\[
 \|D^q\mathbf A\|+\|D^q\mathbf A_{\mathrm{pr}}\|
 +\rho^{-\alpha}\|D^q\mathbf A-D^q\mathbf A_{\mathrm{pr}}\|
 \leqslant C_q,
\]
by Proposition~\ref{prop:legendre-initial-estimates}(ii) and \eqref{eqn:density-powers}. Substituting these bounds into \eqref{eqn:source-telescoping} implies
\[
 \|\mathbf F-\mathbf F_{\mathrm{pr},y}\|_{C_G^{0,\alpha}(Q_{\xi,y})} \leqslant C_q\left(\frac{\rho}{r}\right)^\alpha\prod_{i=1}^q  (\|Z_{\gamma_i}\|_{\mathcal X_\alpha}               +\mathbf 1_{\{|\gamma_i|=1\}}),
\]
after enlarging $C$ in the definition \eqref{eqn:source-constants} of $C_q$ to absorb the factor $q$. This bounds the last term in the $\mathcal{Y}_\alpha$-norm \eqref{eqn:source-norm}. On $Q_s(y)$, we decompose $$\mathbf F-\mathbf F_{\mathrm{pr},y}=(\mathbf F-\mathbf F_{\mathrm{pr},y'}) +(\mathbf F_{\mathrm{pr},y'}-\mathbf F_{\mathrm{pr},y})$$ at each point $(\xi',y')$. The first part is $O((s/r)^\alpha)$ by the previous bound, and the second is $O(r^{-\alpha}|y'-y|^{\alpha/2})=O((s/r)^\alpha)$ by the H\"older continuity of the coefficients of the comparison polynomials. Then Cauchy--Schwarz, in the form of the dual estimate \eqref{eqn:source-dual} of Lemma~\ref{lem:weighted-elementary-inequalities}(ii), bounds the deviation $\mathbf F-\mathbf F_{\mathrm{pr},y}$, completing (i).
\medskip

\textbf{Step 3: Validity of polynomial corrections $\mathscr F_{\beta,\mathrm{pr},y}\in\mathscr A_{p(y)}(\mathcal P)$.} This is essentially the Leibniz rule again, now in the form of the (generalised) chain rule. Fix $|\beta|=k$ and a centre $(0,y)$, and write $p=p(y)$. Let
\[
 p_t=(g_t,\,y+t),\qquad g_t=g(y)+\sum_{1\leqslant|\gamma|\leqslant k-1}\frac{t^\gamma}{\gamma!}D^\gamma g(y), \qquad t\in\R^{n-2},\ |t|\ll1,
\]
which agrees with $p(y+t)$ to order $k-1$ at $t=0$, and the $k$-th order derivatives of $g_t$ vanish. Using the exponential gauge and parallel transport as in the proof of Lemma~\ref{lem:linearised-operator-comparison}(iii), we find a smooth family of flat Dirac models $\Phi_t$ centred at $p_t$, each minimising the energy frozen at $p_t$, whose transformed models $U_t$ agree with $U_{p(y+t)}$ to order $k-1$ at $t=0$. They satisfy
\begin{equation} \label{eqn:model-identity}
 \int\mathbf A^1_{\mathrm{pr}}(p_t,DU_t):DY=0,
\end{equation}
compactly supported admissible $Y$ and, by \eqref{eqn:coefficient-identification},
\[
 \partial_t^\gamma(y+t,g_t,DU_t)\big|_{t=0} =(D_y^\gamma y,\ a_{Z_\gamma}(y),\ DP_{Z_\gamma,y}) =\Xi^{\mathrm{pr}}_\gamma,\qquad 1\leqslant|\gamma|\leqslant k-1.
\]
At order $k$, since
differentiation in $t$ preserves parity and the polynomial form of the
models, we have
\[
 a:=\partial_t^\beta g_t\big|_{t=0}=0,\qquad
 P:=\partial_t^\beta U_t\big|_{t=0}\in\mathcal P .
\]
In the source expansion \eqref{eqn:proper-source}, the term with $q=1$ is the first derivative
of $\mathbf A^1_{\mathrm{pr}}$ applied to $\partial_t^\beta(y+t,g_t,DU_t)|_{t=0}=(D_y^\beta y,a,DP)$, and the terms with $q\geqslant2$ only involve the jets $\Xi^{\mathrm{pr}}_{\gamma_i}$ with $|\gamma_i|\leqslant k-1$. Hence, by \eqref{eqn:leading-source},
\[
\begin{aligned}
 \partial_t^\beta\mathbf A^1_{\mathrm{pr}}(p_t,DU_t)\big|_{t=0}
 &=(\mathbf N_pa)^1
 +D_{DU}\mathbf A^1_{\mathrm{pr}}(p,DU_p)(DP)\\
 &\quad+\begin{cases}
 \partial_{y_j}\mathbf A^1_{\mathrm{pr}}(p,DU_p), & \beta=e_j,\\[2pt]
 \displaystyle\sum_{q=2}^k\frac{1}{q!} \sum_{\substack{\gamma_1+\cdots+\gamma_q=\beta\\\gamma_i\neq0}} \frac{\beta!}{\gamma_1!\cdots\gamma_q!}\,  D^q\mathbf A^1_{\mathrm{pr}}(p,DU_p)  (\Xi^{\mathrm{pr}}_{\gamma_1},\ldots,\Xi^{\mathrm{pr}}_{\gamma_q}), & k\geqslant2,
 \end{cases}\\
 &=(\mathbf N_pa)^1+(\mathbf A_{\mathrm{pr},p}P)^1-\mathbf F^1_{\beta,\mathrm{pr},y} =(\mathbf A_{\mathrm{pr},p}P)^1-\mathbf F^1_{\beta,\mathrm{pr},y}.
\end{aligned}
\]
For $|t|\ll1$, the $t$-derivatives of $\mathbf A^1_{\mathrm{pr}}(p_t,DU_t)$ are bounded on $\operatorname{supp}Y$ uniformly in $t$, so by the Lebesgue's dominated convergence theorem we differentiate \eqref{eqn:model-identity} $\beta$ times and obtain (ii):
\[
 \mathscr F_{\beta,\mathrm{pr},y}=\mathscr A_pP\in\mathscr A_p(\mathcal P).
\qedhere\]
\end{proof}

\subsection{Tangential regularity}
\label{sec:tangential-main}
We now combine Propositions~\ref{prop:weighted-holder-estimate} and \ref{prop:differentiated-equation} in an induction on $|\beta|$. For real-analytic data, we make the estimates quantitative in $|\beta|$.

\begin{theorem}[Tangential regularity]
\label{thm:tangential-regularity}
Suppose that the ambient data are smooth (resp. real-analytic), and let $0<r\leqslant r_0/2$.
For every multi-index $\beta\neq0$, the field $Z_\beta=D_y^\beta U$ belongs to $W_G^{1,2}(Q_r(y_0))$ and solves $\mathscr B_U(Z_\beta,Y)=\mathscr F_\beta(Y)$ for every compactly supported admissible $Y$, with the source \eqref{eqn:proper-source}. Moreover, the comparison polynomials of $Z_\beta$ satisfy \eqref{eqn:coefficient-identification}, and
\begin{equation}
 \|D_y^\beta U\|_{\mathcal X_\alpha(Q_r(y_0))}\leqslant C_\beta \qquad(\text{resp. }\leqslant C^{|\beta|+1}|\beta|!),
 \label{eqn:tangential-factorial-bound}
\end{equation}
where $C$ is independent of $\beta$. In particular,
\begin{equation}
 |D_\xi D_y^\beta\varphi(\xi,y)|\leqslant C_\beta|\xi|, \qquad (\xi,y)\in Q_r(y_0)\setminus\{\xi=0\},
 \label{eqn:tangential-normal-bound}
\end{equation}
the branching locus is smooth (resp. real-analytic), and $U$ and the coefficients $g,q,M$ of the transformed models $U_{p(y)}$ depend smoothly (resp. real-analytically) on $y$.
\end{theorem}

\begin{proof}[Proof for smooth data]
We argue by induction on $k=|\beta|$. For $\beta=0$, these are the transformed estimates from Part~I (Proposition~\ref{prop:legendre-initial-estimates}). Each step below shrinks the domain by a fixed factor, from $Q_{4r}(y_0)$ to $Q_r(y_0)$. This is harmless since the assertions are local: we prove them at each order on every compact subset of $Q_{r_0}(y_0)$, by covering it with finitely many such cylinders. So we assume that the assertions hold for all orders below $k$ on $Q_{4r}(y_0)$.

Now, let $|\beta|=k$, choose $j$ with $\beta_j>0$,and set $\gamma=(\beta_1,\ldots,\beta_j-1,\ldots,\beta_{n-2})$. For $0<|h|\leqslant r^2/100$, consider the tangential difference quotient
\[
\delta_{h,j}Z_\gamma(\xi,y)=\frac{Z_\gamma(\xi,y+he_j)-Z_\gamma(\xi,y)}{h}.
\]

\medskip
\textbf{Step 1: $Z_\beta\in W_G^{1,2}(Q_{2r}(y_0))$ solves the weak equation.} For $k=1$, subtracting the transformed weak equation \eqref{eqn:full-weak-equation} from its translate by $he_j$ and applying the Newton--Leibniz formula
along the segment from $U$ to $U_{he_j}$, we have
\begin{align*}
    &\overline{\mathscr B}_{he_j}(\delta_{h,j}U,Y) \\
    & = -\int_0^1\!\!\int(\partial_{y_j}\mathbf A^1(\xi,y+the_j,U_{he_j,t},DU_{he_j,t}):DY +\partial_{y_j}\mathbf A^0(\xi,y+the_j,U_{he_j,t},DU_{he_j,t})\cdot Y)\,dt ,
\end{align*}
where $\overline{\mathscr B}_{he_j}$ is the averaged bilinear form considered in Remark \ref{rmk:translated-backgrounds}.

For $k\geqslant2$, subtracting the equation of $Z_\gamma$ from its translate by $he_j$ implies
\[
 \mathscr B_{he_j}(\delta_{h,j}Z_\gamma,Y)
 =(\delta_{h,j}\mathscr F_\gamma)(Y)
 -\frac{(\mathscr B_{he_j}-\mathscr B_U)(Z_\gamma,Y)}{h},
\]
where $\mathscr B_{he_j}$ is the spine-translated bilinear form considered in Remark \ref{rmk:translated-backgrounds}. Since $\mathbf F_\gamma$ and the coefficients of $\mathscr B_U$ only involve $Z_{\gamma'}$ with $|\gamma'|\leqslant k-2$, their difference quotients only involve $Z_{\gamma'}$ with $|\gamma'|\leqslant k-1$. These are bounded pointwise in the normalised blocks, by the induction hypothesis and the last term in the $\mathcal X_\alpha$-norm \eqref{eqn:field-norm}. Together with the weighted bounds \eqref{eqn:weighted-coefficient-bounds} on the coefficients and the dual estimate \eqref{eqn:source-dual} for weighted sources, the right-hand sides above are bounded in $(W_G^{1,2}(Q_{3r}(y_0)))^*$ uniformly in $h$. Since $\|\delta_{h,j}Z_\gamma\|_{L^2_G(Q_{3r}(y_0))}\leqslant C\|\gradG Z_\gamma\|_{L^2(Q_{4r}(y_0))}$, the weighted Caccioppoli estimates (Corollary~\ref{cor:linearised-caccioppoli}) imply that
\[
 \|\gradG\delta_{h,j}Z_\gamma\|_{L^2(Q_{2r}(y_0))}\leqslant C_k,
 \qquad 0<|h|\leqslant r^2/100 .
\]
Letting $h\to0$, we obtain $\delta_{h,j}Z_\gamma\rightharpoonup Z_\beta$ in $W_G^{1,2}(Q_{2r}(y_0))$, with $\mathscr B_U(Z_\beta,Y)=\mathscr F_\beta(Y)$ on $Q_{2r}(y_0)$ for all compactly supported admissible test fields.

\medskip
\textbf{Step 2: $Z_\beta$ satisfies the hypotheses of Proposition~\ref{prop:weighted-holder-estimate}.} By the induction hypothesis, the fields $Z_{\gamma'}$, $|\gamma'|\leqslant k-1$, have finite
$\mathcal X_\alpha(Q_{2r}(y_0))$-norms and satisfy \eqref{eqn:coefficient-identification}. These are the conditions of Proposition~\ref{prop:differentiated-equation} on $Q_{2r}(y_0)$, so
\[
 \|\mathbf F_\beta\|_{\mathcal Y_\alpha(Q_{2r}(y_0))}<\infty, \qquad  \mathscr F_{\beta,\mathrm{pr},y}\in\mathscr A_{p(y)}(\mathcal P)\quad \text{at every centre }(0,y).
\]
Proposition~\ref{prop:weighted-holder-estimate} therefore implies that $Z_\beta$ has the unique comparison polynomial $P_{Z_\beta,y}\in\mathcal P(\mathscr F_{\beta,\mathrm{pr},y})$ at every centre $(0,y)\in Q_r(y_0)$, and
\[
 \|Z_\beta\|_{\mathcal X_\alpha(Q_r(y_0))} \leqslant C(r^{-N_G/2-2}\|Z_\beta\|_{L^2_G(Q_{2r}(y_0))} +\|\mathbf F_\beta\|_{\mathcal Y_\alpha(Q_{2r}(y_0))})<\infty .
\]
Since $D_\xi D_y^\beta\varphi=D_\xi(Z_\beta)_\varphi$, \eqref{eqn:tangential-normal-bound} follows.

\medskip
\textbf{Step 3: The coefficients of $P_{Z_\beta,y}$ satisfy
\eqref{eqn:coefficient-identification}.} By the last term in the $\mathcal X_\alpha$-norm \eqref{eqn:field-norm}, every field $Z$ of finite $\mathcal X_\alpha$-norm satisfies, uniformly in $y$ and $|\omega|=1$,
\[
 Z_\varphi(\rho\omega,y)=a_Z(y)+\rho^2q_Z(y,\omega)+O(\rho^{2+\alpha}),\qquad
 Z_V(\rho\omega,y)=\rho M_Z(y)\omega+O(\rho^{1+\alpha}).
\]
Applying this to $Z_\gamma$ and $Z_\beta=\partial_{y_j}Z_\gamma$, we find $(a,q,M)_{Z_\beta}=\partial_{y_j}(a,q,M)_{Z_\gamma}=D_y^\beta(g,q,M)$ by \eqref{eqn:coefficient-identification} for $\gamma$.

Steps~1--3 complete the induction. In particular, $g$, $q$ and $M$ are smooth, and so is the branching locus, the graph of $g$.
\end{proof}

For real-analytic data, the difficulty is that each application of the local estimate \eqref{eqn:local-linear-estimate} shrinks the domain. Instead, we consider $\chi^kD_y^\beta U$ for a smooth cut-off function $\chi$ following Kato \cite{Kat96} and use the diffeomorphism invariance of the transformed system. The advantage of this method is that only compactly supported fields are involved, so the estimate \eqref{eqn:compact-support-estimate} applies at every order. The constants are then tracked through Proposition~\ref{prop:differentiated-equation}(i), and a majorant series implies factorial bounds.

\begin{proof}[Proof for real-analytic data]
\textbf{Step 1: Kato's trick and diffeomorphism invariance.} Let $\Omega=Q_{2r}(y_0)$ and $K=\overline{Q_{3r/2}(y_0)}$. Choose a smooth cutoff $\chi=\chi_1(y)\chi_2(\xi)$ supported in $K$, even in $\xi$, equal to one near $\overline{Q_r(y_0)}$, and with $\chi(\xi,y)=\chi(0,y)$ for $|\xi|\leqslant r$. For $t\in\R^{n-2}$ with $|t|\ll1$, the map
\begin{equation}
 \Psi_t(\xi,y)=(\xi,y+\chi(\xi,y)t),\qquad U_t=U\circ\Psi_t,
 \label{eqn:localised-translation-map}
\end{equation}
is a diffeomorphism of $\Omega$ which equals the identity outside $K$ and preserves parity and $\{\xi=0\}$. We have
\begin{equation}
\begin{aligned}
 D\Psi_t
 &=
 \begin{pmatrix}
 I_2 & 0\\
 tD_\xi\chi & I_{n-2}+tD_y\chi
 \end{pmatrix},
 \quad 
 (D\Psi_t)^{-1}
 =
 \begin{pmatrix}
 I_2 & 0\\[2pt]
 -\dfrac{tD_\xi\chi}{1+D_y\chi\,t}
 &
 I_{n-2}-\dfrac{tD_y\chi}{1+D_y\chi\,t}
 \end{pmatrix},
\end{aligned}
\label{eqn:translation-jacobian-inverse}
\end{equation}
and $\det D\Psi_t=1+D_y\chi\,t$
in the Legendre coordinates $(\xi,y)$. Moreover, $U_t$ solves the weak equation
$\int_\Omega\mathbf A_t^1(\xi,y,U_t,DU_t):DY+\mathbf A_t^0(\xi,y,U_t,DU_t)\cdot Y=0$
for every compactly supported admissible $Y$, with coefficients
\begin{equation}
\begin{aligned}
 \mathbf A_t^1(\xi,y,U,DU)
 &=(\det D\Psi_t)\,\mathbf A^1(\xi,y+\chi t,U,DU(D\Psi_t)^{-1})(D\Psi_t)^{-T},\\
 \mathbf A_t^0(\xi,y,U,DU)
 &=(\det D\Psi_t)\,\mathbf A^0(\xi,y+\chi t,U,DU(D\Psi_t)^{-1}).
\end{aligned}
\label{eqn:pulled-back-field}
\end{equation}
For $|\beta|=k\geqslant1$, consider the localised field
\[
 Z_\beta:=\partial_t^\beta U_t\big|_{t=0}=\chi^kD_y^\beta U .
\]
By the smooth case, $Z_\beta\in W_{G,0}^{1,2}(\Omega)$ is supported in $K$ and has finite $\mathcal X_\alpha(\Omega)$-norm. Since $\chi=\chi(0,y)$ for $|\xi|\leqslant r$, its comparison polynomials have coefficients
\[
 a_{Z_\beta}=\chi(0,y)^kD_y^\beta g,\qquad
 q_{Z_\beta}=\chi(0,y)^kD_y^\beta q,\qquad
 M_{Z_\beta}=\chi(0,y)^kD_y^\beta M .
\]

\medskip
\textbf{Step 2: The recursive estimate for
$A_k:=\max_{|\beta|=k}\|Z_\beta\|_{\mathcal X_\alpha(\Omega)}$.} Differentiating the equation of $U_t$ at $t=0$ as in Step~1 of the proof for smooth data, we see that $Z_\beta$ solves $\mathscr B_U(Z_\beta,Y)=\mathscr F_\beta(Y)$, where $\mathbf F_\beta$ is given by \eqref{eqn:proper-source} with $(y,\mathbf A)$ replaced by $(t,\mathbf A_t)$ and $\Xi_\gamma=\partial_t^\gamma(t,U_t,DU_t)|_{t=0}$.
For $|\xi|\leqslant r$, we have $D_\xi\chi=0$, so $\Psi_t$ is induced by $(z,y)\mapsto(z,y+\chi(0,y)t)$ in the original coordinate chart. Hence right multiplication by $(D\Psi_t)^{-1}$ only changes the tangential columns of $DU$, and $\mathbf A_t$ keeps the structure \eqref{eqn:density-powers} of the transformed density. For $|\xi|\geqslant r$, the weighted and ordinary derivatives are uniformly comparable. 
Moreover, \eqref{eqn:translation-jacobian-inverse} and \eqref{eqn:pulled-back-field} are real-analytic in $(t,U,DU)$ on a fixed neighbourhood of the normalised data of $U$.

Applying Proposition~\ref{prop:differentiated-equation} to $\mathbf A_t$, we get
\begin{itemize}
\item for $q\geqslant1$,
\begin{equation}
 \|D^q_{(t,U,DU)}\mathbf A_t\big|_{t=0}(\Xi_{\gamma_1},\ldots,\Xi_{\gamma_q})\|_{\mathcal Y_\alpha(\Omega)} \leqslant\Lambda_0\Lambda_1^qq!\prod_{i=1}^q
 (\|Z_{\gamma_i}\|_{\mathcal X_\alpha(\Omega)}+\mathbf 1_{\{|\gamma_i|=1\}}),
 \label{eqn:analytic-coefficient-bound}
\end{equation}
where $\Lambda_0,\Lambda_1$ depend on $\chi$, but not on $q$ or the multi-indices;
\item the same proof gives $\mathscr F_{\beta,\mathrm{pr},y}\in\mathscr A_{p(y)}(\mathcal P)$. 
\end{itemize}
Since $Z_\beta$ is supported in $K$, the estimate \eqref{eqn:compact-support-estimate} of Proposition~\ref{prop:weighted-holder-estimate} for compactly supported fields gives
\begin{equation}
 A_k\leqslant C_K\max_{|\beta|=k}\|\mathbf F_\beta\|_{\mathcal Y_\alpha(\Omega)},
 \label{eqn:bound-ak-by-source}
\end{equation}
with $C_K$ independent of $k$. Summing \eqref{eqn:analytic-coefficient-bound} over the terms of \eqref{eqn:proper-source}, we obtain $A_1\leqslant C_K\Lambda_0\Lambda_1$ and
\begin{equation}
 A_k\leqslant C_K\Lambda_0\sum_{q=2}^k\Lambda_1^q
 \sum_{\substack{\gamma_1+\cdots+\gamma_q=\beta\\\gamma_i\neq0}}
 \frac{\beta!}{\gamma_1!\cdots\gamma_q!}
 \prod_{i=1}^q(A_{|\gamma_i|}+\mathbf 1_{\{|\gamma_i|=1\}}),
 \qquad k\geqslant2.
 \label{eqn:source-bound}
\end{equation}

\medskip
\textbf{Step 3: Factorial bounds via a majorant series.}
For fixed $j_1+\cdots+j_q=k$, comparing coefficients in
$\prod_i(x_1+\cdots+x_{n-2})^{j_i}=(x_1+\cdots+x_{n-2})^k$ yields the multinomial identity
\[
\sum_{\substack{\gamma_1+\cdots+\gamma_q=\beta\\|\gamma_i|=j_i}}
 \frac{\beta!}{\gamma_1!\cdots\gamma_q!}=\frac{k!}{j_1!\cdots j_q!}.
\]
Thus \eqref{eqn:source-bound} becomes
\begin{equation}
 \frac{A_k}{k!}\leqslant C_K\Lambda_0\sum_{q=2}^k\Lambda_1^q
 \sum_{\substack{j_1+\cdots+j_q=k\\j_i\geqslant1}}
 \prod_{i=1}^q\left(\frac{A_{j_i}}{j_i!}+\mathbf 1_{\{j_i=1\}}\right),
 \qquad k\geqslant2.
 \label{eqn:recurrence-relation}
\end{equation}
Let $b=\max\{1,C_K\Lambda_0\Lambda_1\}\geqslant A_1$, and let $w(s)=\sum_{k\geqslant1}w_ks^k$ be the solution with $w(0)=0$ of
\begin{equation}
 w=bs+C_K\Lambda_0\frac{\Lambda_1^2(s+w)^2}{1-\Lambda_1(s+w)}.
 \label{eqn:scalar-majorant-equation}
\end{equation}
Comparing coefficients of $s^k$, we know the recurrence for $w_k$ is exactly the  equality case of \eqref{eqn:recurrence-relation}, with nonnegative coefficients and $w_1=b\geqslant A_1$. By induction, $A_k/k!\leqslant w_k$. With $c=\Lambda_1+C_K\Lambda_0\Lambda_1^2$, we have the explicit formula for $w(s)$
\[
 w=-s+\frac{1+(b+1)\Lambda_1s
 -\sqrt{(1+(b+1)\Lambda_1s)^2-4(b+1)cs}}{2c},
\]
showing that $w(s)$ is real-analytic near $s=0$; hence $w_k\leqslant C^{k+1}$ for some $C>0$, so $A_k\leqslant C^{k+1}k!$. Since $\chi=1$ near $\overline{Q_r(y_0)}$, we have $Z_\beta=D_y^\beta U$ there, which proves \eqref{eqn:tangential-factorial-bound}. Finally, $|D_y^\beta g|=|a_{Z_\beta}|\leqslant C^{|\beta|+1}|\beta|!$ by \eqref{eqn:coefficient-identification} and the definition \eqref{eqn:field-norm} of the $\mathcal X_\alpha$-norm, and similarly for $q$, $M$ and $U$. This proves real analyticity.
\end{proof}

\bibliographystyle{amsalpha} 
\bibliography{bibliography}
\end{document}